\documentstyle[12pt,leqno]{report}

\def\@chapapp{\vskip -0.25in CHAPTER}

\def\appendix{
 \setcounter{chapter}{0}
 \setcounter{section}{0}
 \setcounter{subsection}{0}
 \def\thechapter{\Alph{chapter}}}
\def\chappendix#1{\chapter*{#1}
\addcontentsline{toc}{chapter}{#1}} 

\newcommand{\seq}{\begin{equation}}                 
\newcommand{\eeq}[1]{\label{#1}\end{equation}
 }

\newtheorem{Theorem}{Theorem}[section]
\newcommand{\sthm}{\begin{Theorem}}         

\newcommand{\ethm}{\end{Theorem}}           

\newtheorem{Corollary}[Theorem]{Corollary}   
\newcommand{\scor}{\begin{Corollary}}       
\newcommand{\ecor}{\end{Corollary}}         

\newtheorem{Lemma}[Theorem]{Lemma}
\newcommand{\slm}{\begin{Lemma}}            

\newcommand{\elm}{\end{Lemma}}              

\newtheorem{Definition}[Theorem]{Definition}
\newcommand{\sde}{\begin{Definition}}            

\newcommand{\ede}{\end{Definition}}              

\newtheorem{Conjecture}[Theorem]{Conjecture}
\newcommand{\sconj}{\begin{Conjecture}}            

\newcommand{\econj}{\end{Conjecture}}              

\newtheorem{Proposition}[Theorem]{Proposition}

\newcommand{\sprop}{\begin{Proposition}}            

\newcommand{\eprop}{\end{Proposition}}              

\newcommand{\smlist}[1]{\begin{list}           
                      {(#1{zzcount})}{\usecounter{zzcount}}}
\newcommand{\elist}{\end{list}}

\newcommand{\longdownarrow}{
\setlength{\unitlength}{1em}
\begin{picture}(0.6,1)
\put(0.3,1){\vector(0,-1){1.5}}
\end{picture}

}

\newcommand{\h}{{\cal H}^1(T^2)\otimes i{\bf R}}
\newcommand{\zz}{L_{1,\delta}^2(\Lambda^0(Y)\otimes i{\bf R})}
\newcommand{\zo}{L_{1,\delta}^2(\Lambda^1(Y)\otimes i{\bf R})}
\newcommand{\oz}{L_{2,\delta}^2(\Lambda^0(Y)\otimes i{\bf R})}
\newcommand{\oo}{L_{2,\delta}^2(\Lambda^1(Y)\otimes i{\bf R})}
\newcommand{\tz}{L_{3,\delta}^2(\Lambda^0(Y)\otimes i{\bf R})}
\newcommand{\w}{L_{2,\delta}^2(W)}
\newcommand{\wo}{L_{1,\delta}^2(W)}
\newcommand{\dd}{d_{T^2}}
\newcommand{\hh}{\cal H}
\newcommand{\A}{\cal A}

\newcommand{\B}{\cal B}
\renewcommand{\S}{\cal S}
\newcommand{\G}{\cal G}
\newcommand{\M}{\cal M}
\newcommand{\E}{\cal L}
\newcommand{\cs}{\cal{CSD}}
\newcommand{\s}{s^{\prime}}

\newcommand{\css}{{\cal{CSD}}^{\prime}}

\newcommand{\K}{\cal K}

\renewcommand{\ker}{Ker}

\begin{document}
\author{Weimin Chen}

\title{ THE SEIBERG-WITTEN THEORY OF HOMOLOGY 3-SPHERES }


\date{August 1998 }





\maketitle


\pagestyle{empty}
\tableofcontents
\pagestyle{plain}

\pagenumbering{arabic}
\setcounter{page}{0}
\chapter*{Introduction}

\addcontentsline{toc}{chapter}{INTRODUCTION}

\quad 
Let $Y$ be an oriented homology 3-sphere, i.e. $H_\ast(Y)=H_\ast(S^3)$.
Equip $Y$ with a Riemannian metric $g_0$. The unique spin structure on $Y$ gives
rise to a (unique) $SU(2)$ vector bundle $W$ on $Y$ such that the oriented volume form of $Y$ acts on $W$ as identity by Clifford multiplication. Consider pairs
$(A,\psi)$ where $A$ is an imaginary valued 1-form on $Y$ and $\psi$ is a 
smooth section of $W$. The 3-dimensional Seiberg-Witten equations for $(A,\psi)$ read as
$$
\left\{\begin{array}{c}
D_{g_0}\psi+A\psi=0 \\ \ast dA+\tau(\psi,\psi)=0.
\end{array} \right.
$$
Here $D_{g_0}$ is the Dirac operator on $Y$ associated to the metric $g_0$ and $\tau(\cdot,\cdot)$ is a certain bilinear form on $\Gamma(W)$ with values in 
the 
space of imaginary valued 1-forms on $Y$. The group of gauge transformations ${\G}(Y)=Map(Y,S^1)$ acts on the pairs
$(A,\psi)$ by the following rule: 
$$
s\cdot (A,\psi)=(A-s^{-1} ds, s\psi) \hspace{2mm} \mbox {for} \hspace{2mm} s\in {\G}(Y).
$$
The Seiberg-Witten moduli space ${\M}(Y)$ is the space of gauge equivalence classes of solutions to the Seiberg-Witten equations (these solutions are called monopoles). It is compact and has virtual dimension zero.

The algebraic count of the elements in ${\M}(Y)$ is called the Seiberg-Witten invariant of $Y$ and is denoted by $\chi(Y)$ throughout.
${\M}(Y)$ can be regarded as the set of critical points of the 
Chern-Simons-Dirac functional and $\chi(Y)$ its Euler characteristic.
 
The first question we consider is whether the Seiberg-Witten invariant $\chi(Y)$ is independent of the data involved in its definition, such as the Riemannian 
metric on $Y$ and the perturbations of the Seiberg-Witten equations.
Unfortunately, the answer to this question turns out to be negative. To be more
precise, suppose that the oriented homology 3-sphere $Y$ bounds a smooth spin
4-manifold $X$ endowed with a Riemannian metric which is a product near $Y$.
We set 
$$
\alpha(Y)=\chi(Y)-(\mbox{index} D_X+\frac{1}{8}\mbox{Sign}(X)),
$$
where $D_X$ is the Dirac operator on $X$ defined with the APS global boundary
condition (\cite {APS}) and Sign$(X)$ is the signature of $X$. 

In Chapter 1,
we give a rigorous definition of $\chi(Y)$ and $\alpha(Y)$ and prove the following theorem.

{\bf Theorem A}{\it

Let $Y$ be an oriented homology 3-sphere. Then 
\begin{enumerate}
\item $\alpha(Y)$ is a topological invariant of $Y$, and 
$\alpha(Y)+\alpha(-Y)=0$.
\item $\alpha(Y)\equiv \mu(Y) \pmod 2$, where $\mu(Y)$ is the Rohlin invariant of $Y$.
\end{enumerate}
}

The Casson's invariant satisfies both these properties. Thus
this result strongly supports the recent conjecture of Kronheimer
and Mrowka (\cite{KM2}) that $\alpha(Y)$ equals Casson's
invariant of $Y$.

In order to define $\chi(Y)$, we need to
consider the following perturbations of the Seiberg-Witten equations:
$$
\left\{\begin{array}{c}
D_g\psi + A\psi + f\psi=0 \\ \ast dA+\tau(\psi,\psi)+\mu=0,
\end{array} \right. 
$$
where $g$ is a perturbation of the metric $g_0$,
$f$ is a real valued smooth function on $Y$ and $\mu$ is a small, co-closed, imaginary valued 1-form on $Y$.
The topological invariance of $\alpha(Y)$ is roughly saying that the space of pairs $(g,f)$ has a chamber structure and the
Seiberg-Witten invariant  $\chi(Y)$ depends only on the chamber of the perturbed Dirac operator $D_g+f$ (assuming the perturbation $\mu$ is small). In \cite{H},
Hitchin studied a family of metrics on $S^3$ which shows that the Dirac
operator associated to this family of metrics has infinitely many
different chambers. Using Hitchin's observation, we show that even for the simplest 3-manifold, $S^3$, the Seiberg-Witten invariant $\chi(S^3)$ takes infinitely many different values.

The rest of this thesis is devoted to understanding the Seiberg-Witten
invariant $\chi(Y)$ in the following geometric setting. Assume that $Y$ is decomposed into a union of two submanifolds $Y_1$ and $Y_2$ by an embedded torus $T^2$ where $Y_2$ is diffeomorphic to $D^2\times S^1$. We put a Riemannian metric on $Y$ such that a collar neighborhood of $T^2$ is isometric to $(-1,1)\times {\bf R}/2\pi{\bf Z}\times {\bf R}/2\pi{\bf Z}$ and $Y_2$ carries a metric whose scalar curvature is
non-negative and somewhere positive. By inserting cylinders $[0,2L+1]\times T^2$, we obtain a family of stretched versions $Y_L$ of $Y$.
Our goal is to express the Seiberg-Witten invariant $\chi(Y_L)$ in terms of $Y_1$ and $Y_2$ when the neck is sufficiently long.
We regard $Y_L$ as a result of cutting and pasting of two
cylindrical end manifolds obtained by attaching infinite cylinders to $Y_1$
and $Y_2$ (still denoted by $Y_1$ and $Y_2$ for simplicity). It turns out that
the (finite energy) Seiberg-Witten moduli spaces of the cylindrical end manifolds $Y_1$ and $Y_2$ are
generically 1-dimensional manifolds which are immersed into the space of 
equivalence classes of flat
$U(1)$ connections on $T^2$ via a map which sends a finite energy monopole to 
its limiting value at the infinity of the cylindrical end. After fixing
orientations, these moduli spaces define an ``intersection'' number
$\#{\S} (Y_1,Y_2)$, which we prove equals to the Seiberg-Witten invariant $\chi(Y_L)$ when the length of the neck is large enough. This result is refered to as the gluing formula of $\chi$.

{\bf Theorem B}{\it

For large enough $L$, $\chi(Y_L)=\#{\S} (Y_1,Y_2)$.
}

In Chapter 2, we set up the Fredholm theory for Seiberg-Witten equations on cylindrical end 3-manifolds. The issue of
perturbation and transversality, and analytic properties of the finite energy monopoles
such as exponential decay estimates and ``compactness'' are discussed.
The gluing formula is proved in Chapter 3.

Two technical results needed in Chapters 2 and 3 are included as Appendice
A and B.

Part of this thesis has appeared in the Proceedings of 5th G\"{o}kova
Geometry-Topology Conference (1996) (\cite {C1},\cite{C2}).
\vskip .1in
{ \large\bf Acknowledgment:}
This thesis grew out of an unsuccessful endeavor searching for a 
homology bordism invariant lifting of 
the Rohlin invariant of an oriented homology 3-sphere via the Seiberg-Witten theory.  The existence of such an invariant would imply that there are no
$Z_2$ torsion elements with
non-zero Rohlin invariant in the 3-dimensional homology bordism group, which in turn would imply that not every higher
dimensional topological manifold is simplicially triangulable. I am very grateful to my
thesis advisor Professor Selman Akbulut for suggesting that I work on this problem and for sharing with me his ideas of using gauge theory. I have benefited greatly from numerous discussions with him in the past three years.
I would not have been able to survive the hard work in these years
without his patience, encouragement, and support.

I would also like to thank other members of my Thesis Committee,
Professors Ron Fintushel, Tom Parker, Jon Wolfson and Zhengfang Zhou, as well as other faculty members in the department for their generosity in educating me,
their interest in my work, their help, and the useful conversations
with them.
Thanks also go to Professor Wei-Eihn Kuan, Director of Graduate Studies, for
his generosity and support in these years, and to the Department of Mathematics at Michigan
State University for her financial assistance during my participation in
various
conferences. While working on this thesis, I was awarded with a Summer
Acceleration Fellowship in 1996 and a Dissertation Completion Fellowship in
1997 from the Graduate School.

During the years of my graduate study at MSU, I was fortunate to meet several outstanding fellow graduate students. I benefited greatly from interactions with them. I especially wish to thank Liviu Nicolaescu who has taught me my first lessons in elliptic PDE, and Slava Matveyev who has taught me a great deal of geometric topology and whose generous help has made my life a lot easier.

Finally thanks to my family and friends, especially my wife Zhaorong, without
whose love, understanding, and support, my Ph.D would not have been possible.

\chapter{Topological Invariance}

\section{Seiberg-Witten theory in dimension 3}

Let $Y$ be an oriented homology 3-sphere equipped with a Riemannian metric $g$ (many facts stated in this section hold for general 3-manifolds). There exists 
a unique $SU(2)$ vector bundle $W_0$ over $Y$ as a Clifford module of the Clifford algebra bundle $Cl(TY)\otimes_{\bf R}{\bf C}$ such that the oriented volume form on $Y$ acts as identity on $W_0$. Let $W=W_0\otimes L$, where $L$ is the trivial complex line bundle over $Y$. $W$ is a $U(2)$ vector bundle.

Let $(e^1,e^2,e^3)$ be an oriented local orthonormal basis of $T^{\ast}Y$. This 
gives rise to a local unitary basis of $W_0$ and $W$, within which the 
Clifford multiplication is given by the following matrices:
$$
c(e^1)=\left(\begin{array}{clcr}
i & 0\\
0 & -i
\end{array} \right),
c(e^2)=\left(\begin{array}{clcr}
0 & -1\\
1 & 0
\end{array} \right),
c(e^3)=\left(\begin{array}{clcr}
0 & i\\
i & 0
\end{array} \right).
$$
Let $\psi=(z,w)$, $\phi=(u,v)$, $\psi, \phi\in W$, we define
$$\tau(\psi,\phi)=
{1\over 2}\left(\begin{array}{clcr}
{Re}(z\bar{u}-w\bar{v}) & z\bar{v}+\bar{w}u\\
\bar{z}v+w\bar{u}       & -{Re}(z\bar{u}-w\bar{v})
\end{array} \right).
$$
It is straightforward to show

\slm
$i\tau(\psi,\phi)={1\over 2}
({Re}(z\bar{u}-w\bar{v})(e^1)+{Im}(z\bar{v}+\bar{w}u)(e^2)+
{Re}(z\bar{v}+\bar{w}u)(e^3))$, so $\tau(\psi,\phi)\in \Lambda^1(Y)\otimes i
{\bf R}$. Moreover, we have
$$\langle ie\cdot\psi,\phi \rangle_{Re}=-2\langle e,i\tau(\psi,\phi)\rangle$$
for any $e\in \Lambda^1(Y)$, and $|\tau(\psi,\psi)|^2={1\over 4}|\psi|^4$.
\elm

The Levi-Civita connection of the Riemannian metric $g$ lifts to
a connection on $W_0$. Coupled with a $U(1)$ connection $A$ on the complex line bundle $L$, the 
Dirac operator $D_A$: $\Gamma(W)\longrightarrow \Gamma(W)$ is given in a local
frame by
$$
D_A=\sum_{j=1}^{3}e^j\cdot(\nabla_{e_j}+iA_j).
$$

Let ${\cal A}={\cal C}\times \Gamma(W)$ where ${\cal C}$ is the space of smooth $U(1)$ connections on $L$. The gauge group ${\cal G}={Map}(Y,S^1)$ acts on ${\cal A}$ by $s\cdot(A,\psi)=(A-s^{-1}ds,s\psi)$, $s\in{\cal G}$, $(A,\psi)\in{\cal A}$. Note that $\pi_0({\cal G})=H^1(Y,{\bf Z})=0$. Each element in ${\cal G}$ can be written as $e^f$ with $f\in\Gamma(\Lambda^{0}(Y)\otimes i{\bf R})$ determined up to a constant $2\pi ik$, $k\in\bf Z$.
So ${\cal G}=K({\bf Z},1)$. Let ${\cal B}={\cal A}/{\cal G}$. The action of
${\cal G}$ is free on the subset ${\cal A}^{\ast}={\cal A}\setminus\{\psi\equiv 0\}$,
and with stabilizer $S^1$ on the rest. Hence ${\cal B}^{\ast}={\cal A}^{\ast}/{\cal G}$
is homotopic to $CP^{\infty}$.

We shall work within the context of Sobolev spaces and Banach manifolds.
By fixing a trivialization of $L$, ${\cal C}$ can be identified with 
$\Omega^1(Y)\otimes i{{\bf R}}$, the space of imaginary valued 1-forms on $Y$. Define ${\cal A}_1^2=L_1^2(\Lambda^1(Y)\otimes i{{\bf R}})\times
L_1^2(W_0)$, ${\cal G}_1^2=\{L_2^2$ maps from $Y$ to $S^1 \}$. For simplicity,
we still use the old symbols to denote the Sobolev objects.

\slm
${\cal B}^{\ast}$ is a Banach manifold whose tangent space at $(A,\psi)$ is
$$
T{\cal B}^{\ast}_{(A,\psi)}=\{(a,\phi)\in{\cal A}| -d^{\ast}a + i\langle i\psi,\phi
\rangle_{Re}=0\}.
$$
\elm

\noindent{\bf Proof:} Standard arguments. The key point is that the operator 
$d^{\ast}d + |\psi|^2$ is invertible if $\psi$ is not identically zero. See 
\cite{FU}.
\hfill ${\quad \Box}$

\noindent{\bf Remark:} A neighborhood of $[(A,0)]$ in ${\cal B}$ is diffeomorphic to ${U}/S^1$, where
${U}=\{(a,\phi)\in{\cal A}| d^{\ast}a=0, \|(a,\phi)\|<\delta \}$.

There is a natural ${\bf Z}_4$ action $\sigma$ on ${\cal A}$ given by
$\sigma(A,\psi)=(-A,J\psi)$, where $J$ is the quaternion structure on $W_0$.
The action $\sigma$ descends to an involution on ${\cal B}$ and acts freely on ${\cal B}^{\ast}$.

The Chern-Simons-Dirac functional on ${\cal A}$ is defined by
$${{\cal {CSD}}}(A,\psi)=
-\frac{1}{2}\int_{Y}A\wedge dA +\frac{1}{2}\int_{Y}
\langle\psi,D_A\psi\rangle_{g Re}{Vol}_g,
$$
which is gauge invariant and descends to ${\cal B}$. It is also
$\sigma$-invariant.
The gradient of ${{\cal {CSD}}}$ at $(A,\psi)$ is given by
$$
s(A,\psi)=
(\ast dA + \tau(\psi,\psi),D_{A}\psi).
$$
It can be regarded as a `weak'
tangent vector field on ${\cal B}^{\ast}$ in the sense that it is not in
$T{\cal B}^{\ast}$ but in its $L^2$ completion ${{\cal L}}$, i.e., ${{\cal L}}_{(A,\psi)}=
\{(a,\phi)\in L^2 | -d^{\ast}a+i\langle i\psi,\phi\rangle_{Re}=0 \}$.

The covariant derivative $\nabla s$ is given by
$$\nabla s_{(A,\psi)}(a,\phi)=
(\ast da+2\tau(\psi,\phi)-df(\phi),D_{A}\phi+a\psi+f(\phi)\psi)
$$
where $f(\phi)$ is the unique solution to the equation $(d^{\ast}d+|\psi|^2)f
=i\langle iD_{A}\psi,\phi\rangle_{Re}$. As in \cite{T2}, we have

\slm
$\nabla s_{(A,\psi)}$ defines a closed, essentially selfadjoint, Fredholm operator 
on ${{\cal L}}_{(A,\psi)}$, and its eigenvectors form an $L^2$-complete orthonormal basis
for ${{\cal L}}_{(A,\psi)}$. The domain of $\nabla s_{(A,\psi)}$ is the $L^2_1$-Sobolev
space completion of ${{\cal L}}_{(A,\psi)}$. The eigenvalues form a discrete subset of the real line which has no accumulation points, and which is unbounded in both
directions. Each eigenvalue has finite multiplicity.
\elm

The 3-dimensional Seiberg-Witten moduli space ${{\cal M}}$ is the set of critical
points of ${{\cal {CSD}}}$ on ${\cal B}$, i.e. the  equivalence classes of solutions to
the Seiberg-Witten equations
$$
\left\{\begin{array}{c}
\ast dA + \tau(\psi,\psi)=0 \\ D_{A}\psi=0.
\end{array} \right.
$$
Let $[\theta]$ denote the unique reducible solution $[(0,0)]$. Then the moduli
space of irreducible solutions is ${{\cal M}}^{\ast}={{\cal M}}\setminus [\theta]$. As in \cite{KM1}, we have

\slm
The moduli space ${{\cal M}}$ can be represented by smooth sections and it is compact.
\elm

In order to define the Seiberg-Witten invariant, i.e. the 
Euler characteristic of ${{\cal {CSD}}}$, we need 
suitable perturbations of ${{\cal {CSD}}}$.

\sde
A perturbation ${\cal {CSD}}^\prime$ of ${{\cal {CSD}}}$ is admissible if:
\begin{enumerate}
\item The critical points of ${\cal {CSD}}^\prime$ in ${\cal B}^{\ast}$ are non-degenerate, i.e.
$\nabla{s^\prime}_{[(A,\psi)]}$ is invertible at ${[(A,\psi)]}\in {\cal B}^{\ast}\bigcap {s^\prime}^{-1}(0)$.
\item The Dirac operator at the reducible $[\theta]$ is invertible so that $[\theta]$ is isolated. 
\end{enumerate}
Here $s^\prime$ is the gradient of ${\cal {CSD}}^\prime$ and $\nabla{s^\prime}$ is the covariant
derivative of $s^\prime$. The Dirac operator at $[\theta]$ will be clear when we specify the perturbation.
\ede

An admissible perturbation has only finitely many isolated critical points in ${\cal B}^{\ast}$. This is because the reducible $[\theta]$ is isolated so that ${{\M}}^\ast$ is compact. Each irreducible critical point is assigned a sign by the mod $2$ spectral flow of $\nabla{s^\prime}$. Since $\pi_1({\cal B}^{\ast})=0$, the spectral flow does not depend on the path chosen. See 
\cite{T2}.

We will consider two classes of admissible perturbations. The first class is
$\sigma$-invariant. First we need to perturb the Dirac operator so that it is invertible and still quaternionic. These perturbations take the form of $D_g +f$ where $g$
stands for the metric and $f$ is a smooth real valued function on $Y$. The perturbed Chern-Simons-Dirac functional takes the form of
$$
{{\cal {CSD}}^\prime}(A,\psi)={{\cal {CSD}}}(A,\psi)+\frac{1}{2}\int_{Y}f|\psi|_g^2{Vol}_g +u,
$$
where $u$ is some functional on $\B$ which will be constructed in Section 1.5. The corresponding Dirac operator at the reducible $[\theta]$ is $D_g+f$.
For convenience, we set
$$
{{\cal {CSD}}}_f (A,\psi)={{\cal {CSD}}}(A,\psi)+\frac{1}{2}\int_{Y}f|\psi|_g^2{Vol}_g.
$$ 

The following proposition is proved in Section 1.4, in which $Met$ stands for the space of metrics.

\sprop
Let $Y$ be a closed oriented 3-manifold.
For a generic pair $(g,f)\in {Met}\times C^k(Y)$, the perturbed Dirac 
operator ${D}_g + f$ is invertible. Moreover, any two such regular
pairs $(g_0,f_0)$ and $(g_1,f_1)$ can be connected by a generic path
$(g_t,f_t)$ such that the perturbed Dirac operators ${D}_{g_t}+f_t$
are invertible except for $t_i\in(0,1)$ with 
Ker$({D}_{g_{t_i}}+f_{t_i})=\bf H$, $i=1,2,\ldots,n$.
Let $\lambda_t,\psi_t$ be the eigenvalue and eigenvector near $t_i$,
i.e. $({D}_{g_t}+f_t)\psi_t=\lambda_t \psi_t$ with $\lambda_{t_i}=0$ and 
$\|\psi_t\|_{L^2}=1$,
we have
$$\frac{d\lambda_t}{dt}(t_i)=
\int_{Y}\langle \frac{d}{dt}({D}_{g_t}+f_t)(t_i)(\psi_{t_i}),\psi_{t_i}
\rangle_{Re}\neq 0.
$$
As a corollary, the spectral flow of ${D}_{g_t}+f_t$ at $t_i$ is $\pm 4$ for 
$i=1,2,\ldots,n$.
\eprop

The next proposition concerning the existence of $\sigma$-invariant
admissible perturbations is proved in Section 1.5.

\sprop
Fix a regular pair $(g,f)$ so that the reducible $[\theta]$ is isolated.
There exist $\sigma$-invariant admissible perturbations of ${{\cal {CSD}}}_f$ which are supported in the complement of
$[\theta]$ and the non-degenerate critical points of ${{\cal {CSD}}}_f$. Any two such admissible perturbations can be connected by a path supported in the complement of $[\theta]$.
\eprop

The second class of admissible perturbations of ${{\cal {CSD}}}$ has the form 
of
$$
{{\cal {CSD}}^\prime}_{\mu}(A,\psi)={{\cal {CSD}}}(A,\psi)-\int_{Y} A\wedge \ast\mu
$$
where $\mu$ is a generic imaginary valued co-closed 1-form.
The gradient of ${{\cal {CSD}}^\prime}_{\mu}$ at $(A,\psi)$ is 
$$
{s^\prime}_{\mu}(A,\psi)=(\ast dA + \tau(\psi,\psi) + \mu, D_{A}\psi).
$$
${{\cal {CSD}}^\prime}_{\mu}$ has a unique reducible critical point $[\theta_{\mu}]=[(a_{\mu},0)]$ where
$a_\mu$ is the unique solution to the equations $\ast da_\mu+\mu=0$ and $d^{\ast}a_\mu=0$.
The covariant derivative $\nabla{s^\prime}_\mu$ is given by
$$\nabla{s^\prime}_{\mu,{(A,\psi)}}(a,\phi)=(\ast da+2\tau(\psi,\phi)-df(\phi),
D_{A}\phi+a\psi+f(\phi)\psi)
$$
where $f(\phi)$ is the unique solution to the equation 
$$
(d^{\ast}d+|\psi|^2)f=i\langle iD_{A}\psi,\phi\rangle_{Re}.
$$ 
The corresponding Dirac operator at $[\theta_\mu]$ is $D_\mu=D+a_\mu$.

\sprop
For a generic $\mu$, ${{\cal {CSD}}^\prime}_\mu$ is admissible.
Moreover, any two such regular $\mu_0$ and $\mu_1$ can be connected by a path
$\mu_t$, $t\in[0,1]$, such that
\begin{enumerate}
\item ${s^\prime}_{\mu_t}$ is transversal to the zero section of the Hilbert
bundle ${{\cal L}}$ over ${\cal B}^{\ast}\times [0,1]$.
\item $D_{\mu_t}$ is invertible for all but finitely many points $t_i\in
(0,1)$ with Ker$D_{\mu_{t_i}}={\bf C}$. Moreover, if $\lambda_t$ and $\psi_t$
are the eigenvalue and eigenvector of $D_{\mu_t}$ near $t_i$, i.e.
$D_{\mu_t} {\psi_t}=\lambda_t \psi_t$ with $\|\psi_t\|_{L^2}=1$ and $\lambda_{t_i}=0$, then
$$
\frac{d\lambda_t}{dt}(t_i)=\int_{Y}\langle \frac{d}{dt}(D_{\mu_t})(t_i)
\psi_{t_i},\psi_{t_i}\rangle_{Re}\neq 0.
$$
In particular, the spectral flow of $D_{\mu_t}$ (as complex linear operators) at $t_i$ is equal to $\pm 1$.
\end{enumerate}
\eprop

\noindent{\bf Proof:} 
The ``universal'' gradient $s(\mu,A,\psi)=(\ast dA +\tau(\psi,\psi) +\mu,D_{A}\psi)$ is a
section of the Hilbert bundle ${{\cal L}}$ over Ker $d^{\ast}\times {\cal A}^{\ast}$
which is transversal to the zero section. So $s^{-1}(0)$
is a Banach manifold, and so is $s^{-1}(0)/{\cal G}$. The 
projection $P$: $s^{-1}(0)/{\cal G}\rightarrow$ Ker $d^{\ast}$ is a
Fredholm map of index $0$. So for a generic $\mu$, $\nabla{s^\prime}_\mu$ is
invertible at ${s^\prime}_\mu^{-1}(0)$, and  any two such regular $\mu_0$ and $\mu_1$ can be connected by a path
$\mu_t$, $t\in[0,1]$, such that ${s^\prime}_{\mu_t}$ is transversal to the zero section of the Hilbert
bundle ${{\cal L}}$ over ${\cal B}^{\ast}\times [0,1]$.

Consider the real Hilbert bundle ${\cal E}$ over Ker $d^{\ast}\times (L^2_1(W_0)
\setminus\{0\})$ given by ${\cal E}_{(a,\psi)}=\{\phi\in L^2(W_0)|\phi$ is 
orthogonal to $i\psi\}$. Then $L(a,\psi)=D\psi+a\psi$ is a section of ${\cal E}$ which is transversal to the zero section. Therefore
$L^{-1}(0)$ is a Banach manifold. The projection $\Pi$:
$L^{-1}(0) \rightarrow$ Ker $d^{\ast}$ is a Fredholm map of index $1$. Since
$D_a=D + a$ is complex linear, by Sard-Smale theorem, for a generic $a\in $ Ker $d^{\ast}$,
$\Pi^{-1}(a)$ is empty, i.e. $D_a$ is invertible. Two such regular $a_0$ and
$a_1$ can be connected by a path $a_t$ which is transversal to $\Pi$. We can
take an analytic path $a_t$ so that for all but finitely many points $t_i$,
$D_{a_t}$ is invertible and Ker$D_{a_{t_i}}={\bf C}$ by index counting.
If $D_{a_t} {\psi_t}=\lambda_t \psi_t$ with $\|\psi_t\|_{L^2}=1$ and $\lambda_{t_i}=0$, then
$$
\frac{d\lambda_t}{dt}(t_i)=\int_{Y}\langle \frac{d}{dt}(D_{a_t})(t_i)
\psi_{t_i},\psi_{t_i}\rangle_{Re}.
$$
Since $a_t$ is transversal to the projection $\Pi$, 
$\int_{Y}\langle \frac{d}{dt}(D_{a_t})(t_i)\psi_{t_i},\psi_{t_i}\rangle_{Re}
\neq 0$.
\hfill ${\quad \Box}$

\noindent{\bf Remark:} The same conclusions hold if we also allow the metrics to change.

\section{The definition of $\chi$ and $\alpha$}

Fix an admissible perturbation ${{\cal {CSD}}^\prime}$ of ${{\cal {CSD}}}$ with gradient ${s^\prime}$. Denote the Dirac operator at the reducible critical point $[\theta]$ by $D^\prime$. Let ${{\cal M}}^{\ast} =\{{[(A,\psi)]}\in {\cal B}^{\ast}| {s^\prime} (A,\psi)=0\}$.
We define for $\beta_j\in{\cal M}^\ast$,
$$
\chi^j=\sum_{\beta_i\in {{\cal M}}^{\ast}}(-1)^{SF(\beta_j,\beta_i)}
$$
where $SF(\beta_j,\beta_i)$ is the spectral flow between $\nabla{s^\prime}_{\beta_j}$ and $\nabla{s^\prime}_{\beta_i}$. As in \cite{T2},
it is easy to show that $|\chi^j|$ is independent of the choice of $\beta_j$. In order to give a sign to $|\chi^j|$, we need to fix a sign near
the reducible critical point $[\theta]$.

At $(A,\psi)\in {\cal A}^{\ast}$, we have a short exact sequence
$$
0\longrightarrow T{\cal G}_{id} \stackrel{d_{(A,\psi)}}{\longrightarrow} T{\cal A}^{\ast}\stackrel{\pi_\ast}{
\longrightarrow} T{\cal B}^{\ast} \longrightarrow 0
$$
where $d_{(A,\psi)}(f)=(-df,f\psi)$ and $\pi: {\cal A}^{\ast} \rightarrow{\cal B}^{\ast}$.
This enables us to extend any endomorphism of $T{\cal B}^{\ast}$ to a ${\cal G}$-equivariant one of $T{\cal A}\oplus T{\cal G}_{id}$. An endomorphism $L$ of 
$T{\cal B}^{\ast}$ is extended to 
$${\cal K}^\prime_L=\left(\begin{array}{clcr}
L & 0  & 0\\
0 & 0  & d_{(A,\psi)}\\
0 & d^{\ast}_{(A,\psi)} & 0
\end{array} \right),
$$
an endomorphism of $T{\cal A}\oplus T{\cal G}_{id}=T{\cal B}^{\ast}\oplus {Im}(d_{(A,\psi)})\oplus T{\cal G}_{id}$. ${\cal K}^\prime_L$ is self-adjoint if and only if $L$ is. For $L=\nabla{s^\prime}_\mu$,
we use ${\cal K}^\prime$ for ${\cal K}^\prime_L$.

At $(A,\psi)\in {\cal A}$, we define a self-adjoint endomorphism of $T{\cal A}\oplus T{\cal G}_{id}$: 
$$
{\cal K}_{(A,\psi)}(a,\phi,f)=(\ast da+2\tau(\psi,\phi)-df, D_{A}\phi+a\psi+f\psi,
-d^{\ast}a+i\langle i\psi,\phi\rangle_{Re})
$$ 
or
$$
{\cal K}_{(A,\psi)}=\left(\begin{array}{clcr}
D_A  &  \psi\cdot  & \psi\cdot\\
2\tau(\psi,\cdot) & \ast d &  -d\\
i\langle i\psi,\cdot\rangle_{Re} & -d^{\ast} & 0
\end{array} \right).
$$
As in \cite{T2}, we have

\slm
For smooth $(A,\psi)\in{\cal A}$, ${\cal K}_{(A,\psi)}$ extends to 
$L^2(\Lambda^1(Y)\otimes i{{\bf R}}\oplus W_0\oplus \Lambda^0(Y)\otimes i{{\bf R}})$
as a closed, essentially selfadjoint, Fredholm operator. It has discrete
spectrum with no accumulation points, and each eigenvalue has finite 
multiplicity. The spectrum is unbounded from above and below. The same holds
for ${\cal K}^\prime_{(A,\psi)}$ if $(A,\psi)\in{\cal A}^{\ast}$. Moreover, one can 
replace $\nabla{s^\prime}$ by ${\cal K}$ for the purpose of computing the spectral
flow.
\elm

For any $(a,\phi)\in{\cal A}^{\ast}$, we need to study the small eigenvalues of ${\cal K}_{t}(a,\phi)={\cal K}_{0}+tC(a,\phi)$ as $t\rightarrow 0$
where
$${\cal K}_{0}=\left(\begin{array}{clcr}
D^\prime & 0 & 0\\
0 & \ast d & -d\\
0 & -d^{\ast} & 0
\end{array} \right), \mbox{ and } 
C(a,\phi)=\left(\begin{array}{clcr}
a & \phi\cdot & \phi\cdot\\
2\tau(\phi,\cdot) & 0 & 0\\
i\langle i\phi,\cdot\rangle_{Re} & 0 & 0
\end{array} \right).
$$
Here $D^\prime$ is the Dirac operator at $[\theta]$ which is invertible. 
${\cal K}_{0}$ has only one zero
eigenvector which is the constant function $i$. ${\cal K}_{t}(a,\phi)$ is expected to have exactly one small eigenvalue $\lambda_t$ which is analytic in $t$ as $t \rightarrow 0$. See \cite{K}.

\slm
$\dot{\lambda_t}(0)=0$, $\ddot{\lambda_t}(0)=-2 \int_{Y}\langle D^\prime
\tilde{\phi},\tilde{\phi}\rangle_{Re}$
where $\tilde{\phi}=(D^\prime)^{-1}(i\phi)$.
\elm

\noindent{\bf Proof:} For simplicity let $K_t={\cal K}_{t}(a,\phi)$, $C=C(a,\phi)$.
Suppose $(K_t-\lambda_t)f_t=0$ where $\|f_t\|=1$, $f_0=i$. By
differentiating the equation, we have 
$$
(C-\dot{\lambda_t})f_t +(K_t-\lambda_t)\dot{f_t}=0.
$$
So $\dot{\lambda_t}=(C(f_t),f_t)$, and $\dot{\lambda_t}(0)=(C(i),i)=(i\phi,i)=0$. $K_0(\dot{f_t}(0))=-C(f_0)=-i\phi$.
Let $\tilde{\phi}=(D^\prime)^{-1}(i\phi)$, then $\ddot{\lambda_t}(0)=
(C(\dot{f_t}(0)),f_0)+(C(f_0),\dot{f_t}(0))=-2\int_{Y}\langle D^\prime
\tilde{\phi},\tilde{\phi}\rangle_{Re}$.

\hfill ${\quad \Box}$

\scor
For a generic $\phi$, $\ddot{\lambda_t}(0)\neq 0$. $\lambda_t \sim \lambda t^2$
where $\lambda=-\int_{Y}\langle D^\prime
\tilde{\phi},\tilde{\phi}\rangle_{Re}$
and $\tilde{\phi}=(D^\prime)^{-1}(i\phi)$.
\ecor

For $\beta_j\in {{\cal M}}^{\ast}$, we define 
$$
\mbox{sign}(\beta_j)=-\mbox{sign}(\int_{Y}\langle D^\prime
\tilde{\phi},\tilde{\phi}\rangle_{Re})\cdot (-1)^{SF(\beta_j,\phi)}
$$
for a generic $\phi$, where $SF(\beta_j,\phi)$ is the spectral flow between
${\cal K}_{\beta_j}$ and ${\cal K}_{t}(a,\phi)$ for small $t$.

\sde
$\chi=\mbox{sign}(\beta_j)\cdot \chi^j$.
\ede

It is easy to see that
sign$(\beta_j)$ is independent of $(a,\phi)$, and $\chi$ is independent of
$\beta_j$ as in \cite{T2}.

\slm
$\chi(Y)=-\chi(-Y)$, and $\chi \equiv 0 \pmod 2$ if ${{\cal {CSD}}^\prime}$
is a $\sigma$-invariant admissible perturbation.
\elm

\noindent{\bf Proof:} $W_0$ still can serve for $-Y$ if we change the Clifford multiplication
by a factor of $-1$. Under this change, ${{\cal {CSD}}^\prime}(Y)=-{{\cal {CSD}}^\prime}(-Y)$,
$\nabla{s^\prime}(Y)=-\nabla{s^\prime}(-Y)$, ${{\cal M}}(Y)={{\cal M}}(-Y)$,
and $\int_{Y}\langle D^\prime
\tilde{\phi},\tilde{\phi}\rangle_{Re}=-\int_{-Y}\langle D^\prime
\tilde{\phi},\tilde{\phi}\rangle_{Re}$. So $\chi(Y)=-\chi(-Y)$.
The other statement is obvious.
\hfill ${\quad \Box}$

Let $X$ be a smooth compact spin 4-manifold with $\partial X=Y$. 
Equip $X$ with a Riemannian metric such that a neighborhood of $Y$ is isometric to $(-1,0]\times Y$. Suppose $D_X$ is a perturbed Dirac operator on $X$ which takes the form 
$$
c(dt)(\frac{d}{dt}+D^\prime)
$$
near the boundary $Y$. Here $D^\prime$ is the Dirac operator at $[\theta]$ for
an admissible perturbation of the Chern-Simons-Dirac functional and takes the form of $D_g+f+a$ where $a$
is a co-closed imaginary valued 1-form, $g$ stands for the metric and $f$ is a
smooth real valued function on $Y$. $D^\prime$ is invertible.
Index$D_X$ is the $L^2$ index if an infinite cylinder is attached to $X$,
or the index of $D_X$ satisfying the APS global boundary
condition. 

\slm (\cite{APS})
$\mbox{Index}D_X + \frac{1}{8}\mbox{Sign}(X)$ is independent of $X$, and
$$
(\mbox{Index}D_X^1 + \frac{1}{8}\mbox{Sign}(X))-(\mbox{Index}D_X^2 + \frac{1}{8}\mbox{Sign}(X))=-SF(D^\prime_1,D^\prime_2),
$$
where $D_X^i=c(dt)(\frac{d}{dt}+D^\prime_i)$ near $Y$.
In the case that $a=0$ and $(g,f)$ is a regular pair, $\mbox{Index}D_X + \frac{1}{8}\mbox{Sign}(X)\equiv \mu(Y) \pmod 2$ where $\mu(Y)$ is the Rohlin invariant.  $\mbox{Index}D_X + \frac{1}{8}\mbox{Sign}(X)$ changes by a factor of $-1$ if the orientation of $Y$ is reversed.
\elm

\sde 
For any admissible perturbation, define
$$
\alpha=\chi-(\mbox{Index}D_X + \frac{1}{8}\mbox{Sign}(X)).
$$
Here $D_X$ takes the form of $c(dt)(\frac{d}{dt}+D^\prime)$ near
$Y$ where $D^\prime$ is the Dirac operator at the reducible critical point $[\theta]$ associated to the admissible perturbation.
\ede

\section{Topological invariance of $\alpha$}

In this section, we shall prove that $\alpha$ is independent of the choice
of the Riemannian metric and admissible perturbation. 

Given any two metrics and admissible perturbations 
${{\cal {CSD}}^\prime}_{\mu_i}$, $i=-1,1$, we can connect them by a
path ${{\cal {CSD}}^\prime}_{\mu_t}$ $t\in [-1,1]$ for which Proposition 1.1.8 holds. We only need to consider the following two situations:
\begin{enumerate}
\item $D_{\mu_t}$ is invertible for all $t$.
\item $D_{\mu_t}$ is invertible for all $t$ but $t=0$.
\end{enumerate}
Here $D_{\mu_t}$ is the Dirac operator at the reducible point $[\theta_{\mu_t}]$.
In the first case, $\mbox{Index}D_X + \frac{1}{8}\mbox{Sign}(X)$ does not change, neither does $\chi$.  In fact, we have

\slm
Suppose two admissible perturbations $\mu_0$ and $\mu_1$ are connected by a
path $\mu_t$ which provides a partial cobordism $Z$ between part of ${{\cal M}}^{\ast}_0$ and part of 
${{\cal M}}^{\ast}_1$. If $\beta_0\in {{\cal M}}^{\ast}_0$ is
cobordant to $\beta_1\in {{\cal M}}^{\ast}_1$ via $Z$, then $SF(\beta_0,\beta_1)$
is even. If $\beta_0\in {{\cal M}}^{\ast}_0$ is
cobordant to $\beta_1\in {{\cal M}}^{\ast}_0$ via $Z$, then $SF(\beta_0,\beta_1)$
is odd. Here $SF(\beta_0,\beta_1)$ stands for the spectral flow between
$\nabla{s^\prime}_{\beta_0}$ and $\nabla{s^\prime}_{\beta_1}$.
\elm

\noindent{\bf Proof:} The lemma follows from the fact that the cobordism $Z$ can be
arranged so that the projection from $Z$ to $[0,1]$ is a Morse function. See
\cite{DK}, p.143.
\hfill ${\quad \Box}$

In the second case, $\mbox{Index}D_X + \frac{1}{8}\mbox{Sign}(X)$ changes by
$\pm 1$. We shall prove that $\chi$ also changes by $\pm 1$ which is 
compatible to the change of $\mbox{Index}D_X + \frac{1}{8}\mbox{Sign}(X)$ so that $\alpha$ remains unchanged. This is done by analyzing the Kuranishi model
near the reducible point at $t=0$.

Nonlinear Fredholm maps between Hilbert spaces admit local reductions to 
finite dimensional maps. Suppose $\Psi$: $X \rightarrow Y$ is a nonlinear 
Fredholm map satisfying $\Psi(0)=0$. Let $T=(d\Psi)_0$. Then there are 
splittings $X=\ker\;T\oplus (\ker\;T)^{\perp}$, $Y={Im}T \oplus {Coker}T$ and
a map $\psi : X \rightarrow {Coker}T$ so that $\Psi$ is equivalent to $T+\psi$
near $0$ by a diffeomorphism of $X$, and $\psi(0)=0$, $(d\psi)_0=0$. Moreover,
$\Psi^{-1}(0)$ is diffeomorphic to $\{\psi|_{\ker T}=0\}$ near $0$. If there
is a group action, the above can be made equivariant.

The detailed construction goes as follows. Let $\pi_k : X \rightarrow \ker T$,
$\pi_c : Y \rightarrow {Coker}T$ be the orthogonal projections. Then
$\chi : X \rightarrow X$ given by $\chi : x \rightarrow
\pi_k(x) + T^{-1}(1-\pi_c)(\Psi(x))$ is a local diffeomorphism at $0$.
Define $\psi(y)=\pi_c(\Psi(\chi^{-1}(y)))$. Then $\Psi \circ \chi^{-1}=
T+\psi$, and $\Psi^{-1}(0)=\{\psi|_{\ker T}=0\}$. See \cite{FU}.

Suppose two admissible perturbations $\mu_{-1}$ and $\mu_{1}$ are connected 
by a path $\mu_t$, $t\in [-1,1]$, in the sense of Proposition 1.1.8 and $D_{\mu_t}$ is invertible except for $t=0$. We will study the Kuranishi model near the reducible point at $t=0$ of the following family of Seiberg-Witten equations
$$
\left\{\begin{array}{c}
{\ast}_t dA + {\tau}_t (\psi,\psi)=0 \\ (D_{\mu_t}+A)\psi=0
\end{array} \right.
$$
where $A \in \ker\;d^{\ast}$. Here $d^\ast$ stands for $d^{\ast_t}$ at $t=0$.

Consider map $\Psi : {{\bf R}}\oplus L^2_1(\ker\;d^{\ast}\oplus W_0)
\rightarrow L^2(\ker\;d^{\ast}\oplus W_0)$ given by
$$
\Psi(t,A,\psi)=(\pi({\ast}_t dA + {\tau}_t (\psi,\psi)), (D_{\mu_t}+A)\psi)
$$
where $\pi : \Omega^1(Y)\otimes i{{\bf R}} \rightarrow \ker\;d^{\ast}$ is the $L^2$
orthogonal projection. Then $\ker\;(d\Psi)_0={{\bf R}}\oplus \ker\;D_0$,
${Coker}(d\Psi)_0=\ker\;D_0$. Here $D_0$ stands for $D_{\mu_0}$.
Write $\psi=\psi_0+\psi_1$ where $\psi_0\in
\ker\;D_0$ and $\psi_1\in (\ker\;D_0)^{\perp}$, then we have a local 
diffeomorphism $\chi : {{\bf R}}\oplus L^2_1(\ker\;d^{\ast}\oplus W_0) \rightarrow
{{\bf R}}\oplus L^2_1(\ker\;d^{\ast}\oplus W_0)$,
\begin{eqnarray*}
\chi : (t,A,\psi_0+\psi_1) & \rightarrow & (t,(\ast d)^{-1}(\pi({\ast}_t dA + {\tau}_t (\psi_0+\psi_1,\psi_0+\psi_1))),\\
                           &             & \psi_0+D^{-1}_0(1-\pi_k)
((D_{\mu_t}+A)(\psi_0+\psi_1))),
\end{eqnarray*}
and $\chi^{-1}(t,0,\psi_0)=(t,A,\psi_0+\psi_1)$ where $A=A(t,\psi_0)$,
$\psi_1=\psi_1(t,\psi_0)$ satisfy 
$$
\left\{\begin{array}{c}
A+(\pi{\ast_t}d)^{-1}(\pi\tau_t(\psi_0+\psi_1))=0 \\
\psi_1+D_0^{-1}(1-\pi_k)(D_{\mu_t}-D_0+A)(\psi_0+\psi_1)=0.
\end{array} \right.
$$

\slm
$(D_{\mu_t}+A(t,\psi_0))(\psi_0+\psi_1(t,\psi_0))\in \ker\;D_0$. If we write
$$
(D_{\mu_t}+A(t,\psi_0))(\psi_0+\psi_1(t,\psi_0))=a\psi_0+ib\psi_0
$$
where $a$, $b$ are real numbers, then $b=0$.
\elm

\noindent{\bf Proof:} 
For simplicity, denote $D_{\mu_t}+A(t,\psi_0)$ by $D$. Then $b\|\psi_0\|^2=\int_{Y}\langle ib\psi_0,i\psi_0\rangle_{Re}
=\int_{Y}\langle D(\psi_0+\psi_1)-a\psi_0,i\psi_0\rangle_{Re}
=\int_{Y}\langle D\psi_1,i\psi_0\rangle_{Re}
=-\int_{Y}\langle i\psi_1,D\psi_0\rangle_{Re}
=-\int_{Y}\langle i\psi_1,a\psi_0+ib\psi_0-D\psi_1\rangle_{Re}=0$.
\hfill ${\quad \Box}$

\slm
There exists a constant $C$ so that for small $s$, if $\|\psi_0\|_{L^2_1}
\leq s$, $t \leq s$, then
$$\|\psi_1(t,\psi_0)\|_{L^2_1}\leq Cs^2, \mbox{ and }
\|A(t,\psi_0)\|_{L^2_1}\leq Cs^2.
$$
\elm

\noindent{\bf Proof:}  We have continuous maps $L^2_1\times L^2_1
\rightarrow L^2$ and $(\ast d)^{-1}, D_0^{-1}: L^2\rightarrow L^2_1$.
Apply Banach lemma to the map
$$
B(A,\psi_1)=((\pi{\ast_t}d)^{-1}(\pi\tau_t(\psi_0+\psi_1)),D_0^{-1}(1-\pi_k)(D_{\mu_t}-D_0+A)(\psi_0+\psi_1)),
$$
which maps $\{\|A\|_{L^2_1}\leq Cs^2$, $\|\psi_1\|_{L^2_1}\leq Cs^2\}$ into itself when $t\leq s$ and $\|\psi_0\|_{L^2_1}\leq s$ for small $s$. The lemma
follows easily.
\hfill ${\quad \Box}$

Next we examine the finite dimensional reduction $\phi|_{\ker\;(d\Psi)_0}:
{{\bf R}}\oplus \ker\;D_0 \rightarrow \ker\;D_0$.
Let $\psi_0\in\ker\;D_0$, $\|\psi_0\|_{L^2}=1$. We have
$$
\phi|_{\ker\;(d\Psi)_0}(t,s\psi_0)=
\pi_k(D_{\mu_t}+A(t,s\psi_0))(s\psi_0+\psi_1(t,s\psi_0)).
$$
Without loss of generality, we assume that $s$ is real and positive. 
By Lemma 1.3.2, $\phi|_{\ker\;(d\Psi)_0}(t,s\psi_0)=0$ if and only if
$$
\int_{Y}\langle D_{\mu_t}(s\psi_0+\psi_1(t,s\psi_0)),s\psi_0\rangle_{Re}+
\int_{Y}\langle A(t,s\psi_0)(s\psi_0+\psi_1(t,s\psi_0)),s\psi_0\rangle_{Re}=0.
$$

\slm
Let $D_{\mu_t}\psi_t=\lambda_t\psi_t$, $\lambda_t(0)=0$, $\psi_t(0)=\psi_0$ 
as in Proposition 1.1.8.
Then
\begin{enumerate}
\item
$$
\int_{Y}\langle D_{\mu_t}(s\psi_0+\psi_1(t,s\psi_0)),s\psi_0\rangle_{Re}=
s^2(\lambda_t+O(st+t^2))
$$ 
as $t,s \rightarrow 0$.
\item
\begin{eqnarray*}
\int_{Y}\langle A(t,s\psi_0)(s\psi_0+\psi_1(t,s\psi_0)),s\psi_0\rangle_{Re}
& = &
2s^4(-\int_{Y}\langle(\ast d)^{-1}(\tau(\psi_0)),\tau(\psi_0)\rangle\\
&  & + O(s+t))
\end{eqnarray*}
as $t,s \rightarrow 0$.
\end{enumerate}
\elm

\noindent{\bf Proof:} Let $D_{\mu_t}\psi_t=\lambda_t\psi_t$, and $\psi_t=a_t\psi_0+b_t\psi_t^{\perp}$
where $\psi_t^{\perp}\in (\ker D_0)^{\perp}$, $\|\psi_t^{\perp}\|_{L^2}=1$, 
$a_t \rightarrow 1$, $b_t=O(t)$.
Then 
$$
\lambda_t=|a_t|^2(D_{\mu_t}\psi_0,\psi_0)+2|b_t|^2\lambda_t-|b_t|^2(D_{\mu_t}\psi_t^{\perp},\psi_t^{\perp}).
$$
Since $a_t \rightarrow 1$, $b_t=O(t)$, we have 
$(D_{\mu_t}\psi_0,\psi_0)=\lambda_t+O(t^2)$.

On the other hand, for any $\psi_2\in (\ker D_0)^{\perp}$, we have
$$
(D_{\mu_t}\psi_2,\psi_0)=
a_t^{-1}b_t(\lambda_t(\psi_t^{\perp},\psi_2)-(D_{\mu_t}
\psi_t^{\perp},\psi_2))=O(\|\psi_2\|\cdot t).
$$
So 
$$
\int_{Y}\langle D_{\mu_t}(s\psi_0+\psi_1(t,s\psi_0)),s\psi_0\rangle_{Re}=
s^2(\lambda_t+O(st+t^2))
$$ 
as $t,s \rightarrow 0$.

For the second assertion, we have
\begin{eqnarray*}
A(t,s\psi_0)
& = & -(\pi{\ast_t}d)^{-1}(\pi\tau_t(s\psi_0+\psi_1(t,s\psi_0)))\\
& = & -(\ast d)^{-1}(\tau(\psi_0))s^2 + O(ts^2+s^3).
\end{eqnarray*}
So
\begin{eqnarray*}
\int_{Y}\langle A(t,s\psi_0)(s\psi_0+\psi_1(t,s\psi_0)),s\psi_0\rangle_{Re}
& = &
2s^4(-\int_{Y}\langle(\ast d)^{-1}(\tau(\psi_0)),\tau(\psi_0)\rangle\\
&  & + O(s+t))
\end{eqnarray*}
as $t,s \rightarrow 0$.
\hfill ${\quad \Box}$

\scor
The equation $\phi|_{\ker\;(d\Psi)_0}(t,s\psi_0)=0$ has exactly one solution $s$
for and only for those $t$ such that $\lambda_t$ and
$\int_{Y}\langle(\ast d)^{-1}(\tau(\psi_0)),\tau(\psi_0)\rangle$ have the
same sign, if $\int_{Y}\langle(\ast d)^{-1}(\tau(\psi_0)),\tau(\psi_0)\rangle
\neq 0$. Moreover, we have $t\sim cs^2$ as $t,s \rightarrow 0$.
\ecor

\noindent{\bf Remark:}  $\int_{Y}\langle(\ast d)^{-1}(\tau(\psi_0)),\tau(\psi_0)\rangle
\neq 0$ is generically true by slightly perturbing $\mu_t$ near $t=0$,
observing that $\int_{Y}\langle(\ast d)^{-1}(\tau(\psi_0)),\tau(\psi_0,\phi)
\rangle=0$ for any $\phi$ implies that $\psi_0=0$, and also observing
that $\mu_t$ is transversal to the projection $\Pi$ (see Proposition 1.1.8).

\slm
Let $(A,\psi)$ be the solution to
$$
\left\{\begin{array}{c}
{\ast}_t dA + {\tau}_t (\psi,\psi)=0 \\ (D_{\mu_t}+A)\psi=0
\end{array} \right.
$$
near the reducible and $t=0$, then $SF({\cal K}_{(A,\psi)},{\cal K}_{\mu_t,s}(0,\psi_0))$
is odd as $t,s \rightarrow 0$.
\elm

\noindent{\bf Proof:} ${\cal K}_{(A,\psi)}$ is an analytic perturbation in $s=(\psi,\psi_0)$ of
$$
{\cal K}_0=\left(\begin{array}{clcr}
D_0 & 0 & 0\\
0 & \ast d & -d\\
0 & -d^{\ast} & 0
\end{array} \right).
$$

${\cal K}_0$ has three zero eigenvectors $E^1=\psi_0$, 
$E^2=\frac{1}{\sqrt{2}}(i\psi_0+i)$,
$E^3=\frac{1}{\sqrt{2}}(i\psi_0-i)$. Let ${\cal K}_{(A,\psi)}E^i_s=\lambda_s^i E_s^i$ where 
$E^i_s(0)=E^i$, $\lambda_s^i(0)=0$. Then
$$
\dot{\lambda_s^1}(0)=0, \mbox{ } \ddot{\lambda_s^1}(0)=-8\int_{Y}\langle(\ast d)^{-1}(\tau(\psi_0)),\tau(\psi_0)\rangle, \mbox{ } \dot{\lambda_s^2}(0)=1,
\mbox{ } \dot{\lambda_s^3}(0)=-1.
$$
So $\lambda_s^1 \sim \lambda s^2$, $\lambda_s^2 \sim s$ and $\lambda_s^3 \sim
-s$ where $\lambda$ has the same sign with $-\lambda_t$ (see Corollary 1.3.5).

On the other hand, by Lemma 1.2.2, ${\cal K}_{\mu_t,s}(0,\psi_0)$ has three small eigenvalues
$\lambda_t,\lambda_t, \lambda_1 s^2$ as $t\rightarrow 0$ and $s=o(t)$ where
$\lambda_1=-(D_{\mu_t}\tilde{\psi_0},\tilde{\psi_0})$ and 
$\tilde{\psi_0}=D_{\mu_t}^{-1}(i\psi_0)$. It is easy to see that
$\lambda_1$ has the same sign with $-(D_{\mu_t}{\psi_0},\psi_0)\sim -\lambda_t$
as $t \rightarrow 0$. So $SF({\cal K}_{(A,\psi)},{\cal K}_{\mu_t,s}(0,\psi_0))$ is odd
as $t,s \rightarrow 0$.
\hfill ${\quad \Box}$

\sthm
Let $Y$ be an oriented homology 3-sphere. Then 
\begin{enumerate}
\item $\alpha(Y)$ is a topological invariant of $Y$, and 
$\alpha(Y)+\alpha(-Y)=0$.
\item $\alpha(Y)\equiv \mu(Y) \pmod 2$, where $\mu(Y)$ is the Rohlin invariant
of $Y$.
\end{enumerate}
\ethm

\noindent{\bf Proof:}
There is a family of irreducible critical points disappearing or being
created when $t$ passes $0$. Call it $\beta_t$. It is easy to see from Lemma 1.3.6 that sign$(\beta_t)=$ sign$\lambda_t$. The rest of ${{\cal M}}_{\mu_t}^
{\ast}$ provides a cobordism between the rest of ${{\cal M}}_{\mu_{-1}}^{\ast}$ and ${{\cal M}}_{\mu_{1}}^{\ast}$. The sign convention fixed near the reducibles does
not change since ${\cal K}_{\mu_t,s}(0,\psi_0)$ has a spectral flow equal to $\pm 1$
when $t$ passes $0$ (the point is that $D_{\mu_t}$ is complex linear).  So
we have $\chi_{\mu_{-1}}-\chi_{\mu_1}=-SF(D_{\mu_{-1}},D_{\mu_1})$ and $\alpha$ remains unchanged. As for $\alpha(Y)+\alpha(-Y)=0$, it follows from Lemmas
1.2.5 and 1.2.6.

The second assertion is an easy consequence of the existence of $\sigma$-invariant admissible perturbations. We will construct them in the
next two sections.
\hfill ${\quad \Box}$ 

\noindent{\bf Remark:}
In \cite{H}, Hitchin studied a family of Riemannian metrics on $S^3$ which shows that the second term in the definition of $\alpha$ may take infinitely many different values.
Therefore we prove that even for the simplest manifold, $S^3$, the Seiberg-Witten invariant $\chi(S^3)$ takes infinitely many different values.

\section{Perturbations of Dirac operator}

In this section, we show that the perturbed Dirac operators $D_g+f$ are
invertible for generic pairs of $(g,f)$ and they admit a chamber structure.

Throughout this section, we assume that $Y$ is a closed oriented 3-manifold.
Given a metric $g$ on $Y$,
let $P_{SO}$ be the orthonormal tangent frame bundle of $Y$.
Let $H\subset GL(3,\bf R)$ be the subset of symmetric matrices with positive eigenvalues,
then $C^k (P_{SO} \times_{Ad} H)$ which is the set of $C^k$ sections of the 
associated fiber bundle $P_{SO} \times_{Ad} H$ parameterizes the $C^k$--smooth Riemannian metrics on $Y$.
We use the $C^k$--norm of  $C^k (P_{SO} \times_{Ad} H)$ to topologize it.
Let $h$ be a section of $P_{SO} \times_{Ad} H$,
$g^h$ be the corresponding metric,
and $P^h_{SO}$ be the orthonormal tangent frame bundle associated to $g^h$.
Let $\xi$ be a given spin structure on $Y$,
$\pi : P_{Spin(\xi)} \rightarrow P_{SO}$, 
$\pi : P_{Spin(\xi)}^h \rightarrow P_{SO}^h$ be the $Spin(3)$ bundles correspondent to the metrics $g$ and $g^h$,
then we have a lifting $\tilde{h}$
$$
\begin{array}{ccc}
P_{Spin(\xi)} & \stackrel{\tilde{h}}{\longrightarrow} & P_{Spin(\xi)}^h\\
\longdownarrow \pi & & \longdownarrow \pi\\
P_{SO} & \stackrel{{h}}{\longrightarrow} & P_{SO}^h.
\end{array}
$$
Note that if $h$ is not symmetric, we may not remain in the same spin 
structure.
Let $V=P_{Spin(\xi)}\times_{\rho}{{\bf C}^2}$,
$V^h=P_{Spin(\xi)}^h\times_{\rho}{{\bf C}^2}$ be the spinor bundles where
$\rho: Spin(3)\longrightarrow SU(2)$ is the standard representation. 
We have an isometry $\tilde{h}: V\longrightarrow V^h$ given by
 $\tilde{h}(\sigma,\theta)=({\tilde{h}}(\sigma),\theta)$.

Let ${{\cal D}}:\Gamma(V)\times C^k(P_{SO} \times_{Ad} H)\longrightarrow \Gamma(V)$ 
be the map defined by 
${{\cal D}}(\psi,h)={\tilde{h}}^{-1}\cdot D_{g^h}\cdot {\tilde{h}}(\psi)$
where $\psi\in \Gamma(V) $ and $h\in C^k(P_{SO} \times_{Ad} H)$.
Let $\sigma$ be a local frame of $P_{Spin(\xi)}$, 
$\pi(\sigma)=(e_1, e_2, e_3)$, and $(f_1,f_2,f_3)=(e_1,e_2,e_3)h$ which is the
local orthonormal frame with respect to the metric $g^h$. Write 
$\psi=(\sigma,\theta)$,$ h=(\pi(\sigma),(h_{ij}))$, then
\begin{eqnarray*}
{{\cal D}}(\psi,h)
 & = &\tilde{h}^{-1}\cdot{D}_{g^h}\cdot(\tilde{h}(\sigma),\theta) \\
 & = &\tilde{h}^{-1}\cdot(\tilde{h}(\sigma),
\sum_{i=1}^3(c_if_i(\theta)-\frac{1}{2} \sum_{k<j}\omega_{kj}^{i}(h)
c_ic_k c_j \theta))\\
 & = & (\sigma,\sum_{i=1}^{3} (c_{i} h_{si} e_{s} (\theta)
- \frac{1}{2} \sum_{k<j} \omega_{kj}^{i} (h) c_{i} c_{k} c_{j} \theta)) 
\end{eqnarray*}
where $\omega_{kj}^{i}(h)$ is the Levi-Civita connection 1-forms of the metric $g^{h}$ with respect to $(f_1,f_2,f_3)$, i.e., $\nabla_{f_i}^{h}{f_j}=
{f_k}{\omega_{kj}^{i}(h)}$, and
$$
c_1=\left(\begin{array}{clcr}
i & 0 \\
0 & {-i}
\end{array} \right),
c_2=\left(\begin{array}{clcr}
0 & {-1} \\
1 & 0
\end{array} \right),
c_3=\left(\begin{array}{clcr}
0 & i \\
i & 0
\end{array} \right).
$$
Direct calculation shows that
\begin{eqnarray*}
\omega_{kj}^{i}(h)
& = &\frac{1}{2}(h_{kr}^{-1}h_{li}h_{sj}+ h_{jr}^{-1}h_{lk}h_{si}
- h_{ir}^{-1}h_{lj}h_{sk})(\omega_{rs}^{l}- \omega_{rl}^{s}) \\
&   &\mbox{} + \frac{1}{2}h_{ks}^{-1}h_{li}e_{l}(h_{sj}) - 
\frac{1}{2}h_{kl}^{-1}h_{sj}e_{s}(h_{li}) + 
\frac{1}{2}h_{js}^{-1}h_{lk}e_{l}(h_{si}) \\ 
&   &\mbox{} -
\frac{1}{2}h_{jl}^{-1}h_{si}e_{s}(h_{lk}) -
\frac{1}{2}h_{is}^{-1}h_{lj}e_{l}(h_{sk}) +
\frac{1}{2}h_{il}^{-1}h_{sk}e_{s}(h_{lj}) 
\end{eqnarray*}
where $\nabla_{e_i}{e_j}={e_k}{\omega_{kj}^{i}}$,
$h_{ij}^{-1}h_{jk}=\delta_{ik}$. See \cite{LMi} and \cite{KN}.

\slm
${{\cal D}}(\cdot,h)$: ${\Gamma(V)}\longrightarrow {\Gamma(V)}$ is smooth in $h$.
Moreover, ${{\cal D}}(\cdot,h)$ is self-adjoint if $\det(h)=1$ pointwise on $Y$.
\elm

\noindent{\bf Proof:} That ${{\cal D}}(\cdot,h)$ is smooth in $h$ follows from the local expressions of ${{\cal D}}(\cdot,h)$ and $\omega_{kj}^{i}(h)$.
For the self-adjointness of ${{\cal D}}(\cdot,h)$, we have
\begin{eqnarray*}
\int_{Y}{\langle{{\cal D}}(\psi,h),\phi\rangle}_{g}{Vol_{g}}
& = &\int_{Y}{\langle{\tilde{h}}^{-1}\cdot{D}_{g^h}\cdot {\tilde{h}}(\psi),
\phi\rangle}_{g}{Vol_{g}} \\
& = &\int_{Y}{\langle{D}_{g^h}\cdot {\tilde{h}}(\psi),
{\tilde{h}}(\phi)\rangle}_{g^{h}}{Vol_{g^{h}}} \\
& = &\int_{Y}\langle\tilde{h}(\psi),{D}_{g^h}(\tilde{h}(\phi)\rangle_
{g^h}{Vol_{g^h}} \\
& = &\int_{Y}\langle\psi,\tilde{h}^{-1}\cdot{D}_{g^h}\cdot \tilde{h}
(\phi)\rangle_{g}{Vol_{g}} \\
& = &\int_{Y}\langle\psi,{{\cal D}}(\phi,h)\rangle_{g}{Vol_{g}}
\end{eqnarray*}
where $Vol_{g}=Vol_{g^h}$ since $\det(h)=1$ pointwise on $Y$.
\hfill ${\quad \Box}$

\slm
Given any metric $g$ on $Y$, let $(e_1,e_2,e_3)$ be an oriented local orthonormal frame in an open subset $A$ of $Y$. 
Let $f$ be a smooth real valued function on $Y$.
Suppose $\psi,\phi\in \ker\;({D}_g +f)$.
If
$$
\frac{d}{dt}(\int_{Y}{\langle{{\cal D}}(\psi,e^{tX}),\phi\rangle}_{g}{Vol_g})=0
$$ 
at $t=0$ for any symmetric matrix function $X$ compactly supported in $A$ 
satisfying ${tr}(X)=~0$, then in $A$ we have
$$
\langle{e_j}\nabla_{e_j}{\psi},\phi\rangle_{g} +
\langle \psi,e_j\nabla_{e_j}\phi\rangle_{g}=-\frac{2}{3}\langle f\psi,\phi
\rangle_g
$$
for $j=1,2,3 $, and
$$
\langle{e_j} \nabla_{e_i}{\psi},\phi \rangle_{g} +
\langle\psi,{e_j} \nabla_{e_i}{\phi}\rangle_{g}=
-\frac{1}{2}e_{k}({\langle\psi,\phi\rangle}_g)
$$
for any $i,j,k$ such that
$e_{i}\wedge{e_j}\wedge{e_k} = {e_1}\wedge{e_2}\wedge{e_3}$.
\elm

\noindent{\bf Remark:} The same conclusions hold with the hermitian product $\langle\cdot,\cdot\rangle_g$ replaced by its real part,
if $\langle\cdot,\cdot\rangle_g$
is replaced by its real part in the condition 
$\frac{d}{dt}(\int_{Y}{\langle{{\cal D}}(\psi,e^{tX}),
\phi\rangle}_{g}{Vol_g})=0$.

The proof of this lemma is a lengthy calculation which is given at the end of
this section.

Let $Met_0$ be the subspace of $Met=C^k(P_{SO} \times_{Ad} H)$ given by
$$
Met_{0}=\{h\in{Met}|\det(h)=1\}.
$$
Every metric in $Met$ is conformal to a metric in $Met_0$.

\newpage

\noindent{\bf The Proof of Proposition 1.1.6:}

Consider the real Hilbert bundle $E$ over the Banach manifold 
$B={Met_0}\times C^k(Y)\times(L_1^2(V)\setminus\{0\})$.
At $(h,f,\psi)\in B$,
$E_{(h,f,\psi)}=\{\phi \in L^2(V)| \phi$ is orthogonal to 
$i\psi,j\psi,k\psi \}$. Here $i,j,k\in \bf H$ satisfying 
$$
ij=k, \hspace{2mm} jk=i,\hspace{2mm} ki=j,\hspace{2mm} \mbox{and} \hspace{2mm}
i^2=j^2=k^2=-1.
$$
The map $L: (h,f,\psi)\longrightarrow {{\cal D}}(\psi,h)+f\psi$ defines a section of the bundle $E$ over the Banach manifold $B$. Suppose that $(h,f,\psi)\in L^{-1}(0)$,
then the differential of $L$ at $(h,f,\psi)$ is
$$\delta L_{(h,f,\psi)}(H,F,\Psi)=
{{\cal D}}(\Psi,h)+f\Psi+\delta{{\cal D}}(\psi,\cdot)(h)(H)+F\psi,
$$
from which it is easy to see that if 
$\phi \in ({Im}\delta L)^{\perp}$, then
$\phi \in \ker\;({{\cal D}}(\cdot,h)+f)$ and 
$\phi=a_1(i\psi)+a_2(j\psi)+a_3(k\psi)$ for some real functions 
$a_1,a_2,a_3$. Moreover, by Lemma 1.4.2,
$$\int_{Y} \langle\delta{{\cal D}}(\psi,\cdot)(h)(H),\phi \rangle_{Re}{Vol}=0$$
for any $H$ implies that 
$$
\langle e_i\nabla_{e_i}\psi,\phi \rangle_{Re} +
\langle \psi,e_i\nabla_{e_i}\phi \rangle_{Re}=-\frac{2}{3}\langle f\psi,\phi
\rangle_{Re}
$$
for $i=1,2,3$, and
$$
\langle e_j\nabla_{e_i}\psi,\phi \rangle_{Re} +
\langle \psi,e_j\nabla_{e_i}\phi \rangle_{Re} = -\frac{1}{2}e_k(
\langle \psi,\phi \rangle_{Re})
$$
for $i,j,k$ such that
$e_{i}\wedge{e_j}\wedge{e_k} = {e_1}\wedge{e_2}\wedge{e_3}$.
From this we obtain that
$$\langle \psi,e_s \cdot(e_l(a_1)(i\psi)+e_l(a_2)(j\psi)+
e_l(a_3)(k\psi))\rangle_{Re}=0
$$
for any $s,l=1,2,3$. Since $\psi$ is not identically zero, we have
$e_l(a_i)=0$ for any $l,i=1,2,3$. Hence $a_1,a_2,a_3$ are constant. So $L$
is transversal to the zero section of $E$ and $L^{-1}(0)$ is a Banach
submanifold in $B$. The projection 
$$P: L^{-1}(0)\longrightarrow {Met_0}\times C^k(Y)
$$
is a Fredholm map of index $3$. Note that $L(h,f,\cdot)={{\cal D}}(\cdot,h)+f$ 
is quaternionic, so by Sard-Smale theorem, for a generic pair $(h,f)
\in {Met_0}\times C^k(Y)$, $P^{-1}(h,f)$ is empty, i.e.,
${{\cal D}}(\cdot,h)+f$ is invertible. Any two such regular pairs $(h_0,f_0)$ and
$(h_1,f_1)$ can be connected by an analytic path $(h_t,f_t)$ which is 
transversal to the projection $P$. The operators ${{\cal D}}(\cdot,h_t)+f_t$
are invertible except for finitely many points $t_i\in(0,1)$,
$i=1,2,\ldots,n$. The fact that $\ker\;({{\cal D}}(\cdot,h_{t_i})+f_{t_i})=\bf H$ follows from index counting. Suppose that 
${{\cal D}}(\psi_t,h_t)+f_t\psi_t=\lambda_t\psi_t$ near $t_i$ with 
$\lambda_{t_i}=0$ and $\|\psi_t\|_{L^2}=1$, then
$$
\frac{d\lambda_t}{dt}(t_i)=
\int_{Y}\langle \frac{d}{dt}({{\cal D}}(\psi_{t_i},h_t)+f_t\psi_{t_i})(t_i),\psi_{t_i}
\rangle_{Re}.
$$
Since the path $(h_t,f_t)$ is transversal to the projection $P$, we have 
$$
\int_{Y}\langle \frac{d}{dt}({{\cal D}}(\psi_{t_i},h_t)+f_t\psi_{t_i})(t_i),\psi_{t_i}
\rangle_{Re}\neq 0.
$$

Suppose $h_1\in Met$ is conformal to $h\in Met_0$ and $g^{h_1}=e^{2u}g^h$.
Let $m:V^{h_1}\rightarrow V^h$ be the isometry. The Dirac operators are
related in the following way (see \cite{H} or \cite{LMi}):
$$
D_{g^h}=e^{2u}mD_{g^{h_1}}m^{-1}e^{-u}.
$$
It is easy to see from this that $D_{g^{h_1}}+f$ is invertible if and only if
$D_{g^h}+e^u f$ is. Similar arguments justify the chamber structure.
\hfill ${\quad \Box}$

\noindent{\bf The Proof of Lemma 1.4.2:}

Let $\psi=(\sigma,\theta)$, $\pi(\sigma)=(e_1,e_2,e_3)$, then
\begin{eqnarray*}
{{\cal D}}(\psi,h)
&=&(\sigma,c_{1}e_{1}(\theta)+c_{2}e_{2}(\theta)+c_{3}e_{3}(\theta)
-\frac{1}{2}((\omega_{12}^{2}(h)+\omega_{13}^{3}(h))c_{1}\theta \\
& &\mbox{}+(\omega_{23}^{3}(h)-\omega_{12}^{1}(h))c_{2}\theta-
(\omega_{13}^{1}(h)+\omega_{23}^{2}(h))c_{3}\theta \\
& &\mbox{}+(\omega_{12}^{3}(h)-\omega_{13}^{2}(h)+\omega_{23}^{1}(h))
c_{1}c_{2}c_{3}\theta)).
\end{eqnarray*}

For $h=e^{tX}$, where $X=\left(\begin{array}{clcr}
x & 0 & 0\\
0 & -x & 0\\
0 & 0 & 0
\end{array}\right)$, we have
\begin{eqnarray*}
\omega_{12}^{2}(h)+\omega_{13}^{3}(h)
&=&(\omega_{12}^{2}+\omega_{13}^{3})(1+tx)+{(1+tx)}^{2}{e_1}(1-tx)+O(t^2),\\
\omega_{23}^{3}(h)-\omega_{12}^{1}(h)
&=&(\omega_{23}^{3}-\omega_{12}^{1})(1-tx)+{(1-tx)}^{2}{e_2}(1+tx)+O(t^2),\\
\omega_{13}^{1}(h)+\omega_{23}^{2}(h)
&=&-(1-tx){e_3}(1+tx)-(1+tx){e_3}(1-tx)+O(t^2),\\
\omega_{12}^{3}(h)-\omega_{13}^{2}(h)+\omega_{23}^{1}(h)
&=&\frac{1}{2}({(1+tx)}^{2}(\omega_{23}^{1}+\omega_{12}^{3})+
{(1-tx)}^{2}(\omega_{12}^{3}-\omega_{13}^{2}))\\
& & +O(t^2).
\end{eqnarray*}
So we have
\begin{eqnarray*}
\frac{d}{dt}({{\cal D}}(\psi,h))(0)
&=&(\sigma,x{c_1}{e_1}(\theta)-x{c_2}{e_2}(\theta)
-\frac{1}{2}((x(\omega_{12}^{2}+\omega_{13}^{3})-{e_1}(x)){c_1}\theta\\
& &\mbox{}-(x(\omega_{23}^{3}-\omega_{12}^{1})-{e_2}(x)){c_2}\theta+
x(\omega_{23}^{1}+\omega_{13}^{2}){c_1}{c_2}{c_3}\theta)).
\end{eqnarray*}
If we write $\psi=(\sigma,\theta)$, $\phi=(\sigma,\xi)$, then 
\begin{eqnarray*}
\int_{Y}\langle \frac{d}{dt}({{\cal D}}(\psi,h))(0),\phi \rangle{Vol}
&=& \int_{Y}(\langle xc_1e_1(\theta),\xi \rangle -
\langle xc_2e_2(\theta),\xi \rangle \\
& &\mbox{}-\frac{1}{2}(x\omega_{12}^{2}+x\omega_{13}^{3}-e_1(x))
\langle c_1\theta,\xi \rangle \\
& &\mbox{}-\frac{1}{2}(x\omega_{23}^{3}-x\omega_{12}^{1}-e_2(x))
\langle c_2\theta,\xi \rangle \\
& &\mbox{}+\frac{1}{2}x(\omega_{23}^{1}+\omega_{13}^{2})
\langle \theta,\xi \rangle){Vol}.
\end{eqnarray*}
Let $(e^1,e^2,e^3)$ be the dual to $(e_1,e_2,e_3)$, then
\begin{eqnarray*}
d(x\langle c_1\theta,\xi \rangle \ast e^1)
&=& e_1(x)\langle c_1\theta,\xi \rangle e^1\wedge e^2\wedge e^3\\
& &\mbox{}+x(\langle c_1e_1(\theta),\xi \rangle +
\langle c_1\theta,e_1(\xi)\rangle)e^1\wedge e^2\wedge e^3\\
& &\mbox{}-x(\omega_{12}^{2}+\omega_{13}^{3})
\langle c_1\theta,\xi \rangle e^1\wedge e^2\wedge e^3.
\end{eqnarray*}
Integration by parts, we have
\begin{eqnarray*}
\int_{Y}e_1(x)\langle c_1\theta,\xi \rangle e^1\wedge e^2\wedge e^3
&=& -\int_{Y}x(\langle c_1e_1(\theta),\xi \rangle +
\langle c_1\theta,e_1(\xi)\rangle)e^1\wedge e^2\wedge e^3\\
& &\mbox{}+\int_{Y}x(\omega_{12}^{2}+\omega_{13}^{3})
\langle c_1\theta,\xi \rangle e^1\wedge e^2\wedge e^3.
\end{eqnarray*}
Similarly, we have
\begin{eqnarray*}
\int_{Y}e_2(x)\langle c_2\theta,\xi \rangle e^1\wedge e^2\wedge e^3
&=& -\int_{Y}x(\langle c_2e_2(\theta),\xi \rangle +
\langle c_2\theta,e_2(\xi)\rangle)e^1\wedge e^2\wedge e^3\\
& &\mbox{}+\int_{Y}x(\omega_{23}^{3}-\omega_{12}^{1})
\langle c_2\theta,\xi \rangle e^1\wedge e^2\wedge e^3.
\end{eqnarray*}
These give us
\begin{eqnarray*}
\int_{Y}\langle \frac{d}{dt}({{\cal D}}(\psi,h))(0),\phi \rangle{Vol}
&=& \frac{1}{2}\int_{Y}x(\langle e_1\nabla_{e_1}\psi,\phi \rangle +
\langle \psi,e_1\nabla_{e_1}\phi \rangle) e^1\wedge e^2\wedge e^3\\
& &\mbox{}-\frac{1}{2}\int_{Y}x(\langle e_2\nabla_{e_2}\psi,\phi \rangle +
\langle \psi,e_2\nabla_{e_2}\phi \rangle) e^1\wedge e^2\wedge e^3.
\end{eqnarray*}
Therefore, if 
$$
\int_{Y}\langle \frac{d}{dt}({{\cal D}}(\psi,h))(0),\phi \rangle{Vol}=0
$$ 
for all $h=e^{tX}$ where 
$X=\left(\begin{array}{clcr}
x & 0 & 0\\
0 & -x & 0\\
0 & 0 & 0
\end{array}\right)$, we have
$$
\langle e_1\nabla_{e_1}\psi,\phi \rangle +
\langle \psi,e_1\nabla_{e_1}\phi \rangle =\langle e_2\nabla_{e_2}\psi,\phi \rangle +\langle \psi,e_2\nabla_{e_2}\phi \rangle.
$$
Similarly, we have
$$\langle e_1\nabla_{e_1}\psi,\phi \rangle +
\langle \psi,e_1\nabla_{e_1}\phi \rangle =\langle e_3\nabla_{e_3}\psi,\phi \rangle +\langle \psi,e_3\nabla_{e_3}\phi \rangle.
$$
But $\psi,\phi\in \ker\;(D_g+f)$, we have
\begin{eqnarray*}
\sum_{i=1}^{3}(\langle e_i\nabla_{e_i}\psi,\phi \rangle +
\langle \psi,e_i\nabla_{e_i}\phi \rangle)
&=& \langle {D}_g\psi,\phi \rangle +\langle \psi,{D}_g\phi \rangle \\
&=& -2\langle f\psi,\phi\rangle.
\end{eqnarray*}
So we have
$$
\langle e_i\nabla_{e_i}\psi,\phi \rangle +
\langle \psi,e_i\nabla_{e_i}\phi \rangle=-\frac{2}{3}\langle f\psi,\phi\rangle
$$
for $i=1,2,3$. Similar computation with $X=\left(\begin{array}{clcr}
0 & x & 0\\
x & 0 & 0\\
0 & 0 & 0
\end{array}\right)$ yields
$$
\langle e_2\nabla_{e_1}\psi,\phi \rangle +
\langle \psi,e_2\nabla_{e_1}\phi \rangle +
\langle e_1\nabla_{e_2}\psi,\phi \rangle +
\langle \psi,e_1\nabla_{e_2}\phi \rangle =0.
$$
Combined with 
$$
(\langle e_2\nabla_{e_1}\psi,\phi \rangle +
\langle \psi,e_2\nabla_{e_1}\phi \rangle)-
(\langle e_1\nabla_{e_2}\psi,\phi \rangle +
\langle \psi,e_1\nabla_{e_2}\phi \rangle) =-e_3(\langle \psi,\phi \rangle),
$$
we have 
$$
\langle e_2\nabla_{e_1}\psi,\phi \rangle +
\langle \psi,e_2\nabla_{e_1}\phi \rangle = -\frac{1}{2}e_3(
\langle \psi,\phi \rangle).
$$ 
In general, we have
$$\langle e_j\nabla_{e_i}\psi,\phi \rangle +
\langle \psi,e_j\nabla_{e_i}\phi \rangle = -\frac{1}{2}e_k(
\langle \psi,\phi \rangle)
$$
for $i,j,k$ such that
$e_{i}\wedge{e_j}\wedge{e_k} = {e_1}\wedge{e_2}\wedge{e_3}$.
\hfill ${\quad \Box}$

\section{The $\sigma$-invariant perturbations}

In this section, we give the construction of the $\sigma$-invariant
admissible perturbations using holonomy along embedded loops. Assume that
$(g,f)$ is regular. Let $s^f$ denote the gradient of ${{\cs}}_f$ and ${{\M}}_f$
denote the set of critical points where 
$$
{\cs}_f={\cs} +\frac{1}{2}\int_{Y}f|\psi|_g^2{Vol}_g,\hspace{2mm}\mbox{and}\hspace{2mm}
s^f (A,\psi)=(\ast dA +\tau(\psi,\psi), D_A\psi+f\psi).
$$
The moduli space ${\M}_f$ is compact and can be represented by smooth sections.

\sde
A thickened loop is an embedding $\gamma : S^1\times D^2 \rightarrow Y$,
together with a bump function $\eta(y)$ on $D^2$ centered at $0\in D^2$, with
$\int_{D^2}\eta(y)dy=1$.
\ede

Given a thickened loop $\lambda=(\gamma,\eta)$, one can define a pair of 
$\sigma$-invariant functions $(p,q)_\lambda : {\B}\rightarrow [-1,1]\times
{{\bf R}}^{+}$ by 
$$
p_\lambda(A,\psi)=\int_{D^2}\cos(\theta_y)\eta(y)dy,
$$
where $e^{i\theta_y}$ is the holonomy of $A$ along the loop $\gamma_y=
S^1\times \{y\}$, and 
$$ 
q_\lambda(A,\psi)=\int_{D^2\times S^1}|\psi|^2\eta(y)dydt.
$$

\slm
The function $(p,q)$ is smooth on ${\A}$.
\elm

\noindent{\bf Proof:} The same arguments as in \cite{T2}. It is useful to know that
$$ 
dp_\lambda|_{(A,\psi)}(a,\phi)=\int_{D^2}i\sin(\theta_y)\eta(y)
(\int_{S^1\times \{y\}}\gamma_y^{\ast}a)dy
$$ 
and
$$ 
dq_\lambda|_{(A,\psi)}(a,\phi)=2\int_{D^2\times S^1}\langle\psi,\phi\rangle \eta(y)dydt.
$$
\hfill ${\quad \Box}$

For any set $\Lambda$ of finitely many thickened loops, we have a smooth map
$\Phi_\Lambda : {\B}^{\ast}\rightarrow \prod_{\lambda\in \Lambda}
([-1,1]\times {{\bf R}}^{+})_\lambda$ given by
$$
\Phi_\Lambda([(A,\psi)])=((p,q)_\lambda(A,\psi),\lambda\in \Lambda).
$$
The map $\Phi_\Lambda$ is $\sigma$-invariant and continuous on ${\B}$.

\slm
There is a set $\Lambda$ of finitely many thickened loops such that
\begin{enumerate}
\item $\ker\;\nabla s^f\bigcap \bigcap_{\lambda\in \Lambda}\ker\;(d(p,q)_\lambda)=
\{0\}$ at any $[(A,\psi)]\in {{\M}_f}^{\ast}$.
\item $\Phi_\Lambda$ is injective up to the $\sigma$ action on ${{\M}_f}$. Therefore
we can identify ${{\M}_f}/{\sigma}$ with a compact subset of $\prod_{\lambda\in \Lambda}([-1,1]\times {{\bf R}}^{+})_\lambda$.
\end{enumerate}
\elm

\noindent{\bf Proof:} Suppose $[(A,\psi)]\in {{\M}_f}^{\ast}$, and $(a,\phi)\in\ker\;\nabla s^f_{(A,\psi)}$,
i.e., $(a,\phi)$ satisfies 
$$\left\{\begin{array}{c}
D_A\phi+f\phi+a\psi =0\\
\ast da+2\tau(\psi,\phi)=0\\
-d^{\ast}a+i\langle i\psi,\phi\rangle_{Re}=0.
\end{array} \right.
$$
Since $A$ is not flat, if $(a,\phi)\in\ker\;(d(p,q)_\lambda)$ for all thickened
loops, then $\int_{S^1\times\{y\}}\gamma^{\ast}_{y}a=0$ for all $\gamma$.
So $da=0$. $da=0$ implies $\tau(\psi,\phi)=0$. So $\phi=v\psi$ for some function $v\in\Omega^0(Y)\otimes i{{\bf R}}$ wherever $\psi\neq 0$. This implies $dv+a=0$ and
$\int_{Y}(|dv|^2+|v|^2|\psi|^2)=0$ by plugging into the equations. Since 
$\psi$ is not identically zero, we have $(a,\phi)=0$.

So for each $[(A,\psi)]\in{{\M}_f}^{\ast}$, there is a set of finitely many thickened
loops such that the first assertion holds for $[(A,\psi)]$. Then the first assertion follows by the compactness of ${{\M}}^{\ast}_f$ and the smoothness of the
function $(p,q)$.

For the second assertion, suppose $[(A_1,\psi_1)],[(A_2,\psi_2)]\in{{\M}_f}^{\ast}$ 
such that $(p,q)_\lambda(A_1,\psi_1)=(p,q)_\lambda(A_2,\psi_2)$ for all loops.
Then $dA_1=\pm dA_2$, and $|\psi_1|^2=|\psi_2|^2$. Assume $dA_1=dA_2$, then
$\tau(\psi_1)=\tau(\psi_2)$. By writing in a local frame, it is easy to see
that $\psi_1=s\psi_2$ for some $s\in {Map}(Y,S^1)$. Then it is easy to see
that $[(A_1,\psi_1)]=[(A_2,\psi_2)]$. In the case of $dA_1=-dA_2$, apply $\sigma$.

Now for any $[(A_1,\psi_1)]\neq [(A_2,\psi_2)]$ in ${{\M}_f}^{\ast}/{\sigma}$, there is a thickened loop $\lambda$ separating them. By the compactness of
${{\M}_f}^{\ast}/{\sigma}$ and the smoothness of $(p,q)$, there exists a set of finitely many loops separating any two points in ${{\M}_f}^{\ast}/{\sigma}$ 
with 
distance greater then a fixed number. Combining with the first assertion, since
each point in ${{\M}_f}^{\ast}/{\sigma}$ has a neighborhood described by a Kuranishi model, the second assertion follows.

\hfill ${\quad \Box}$

For any smooth function $h$ on $\prod_{\lambda\in \Lambda}
([-1,1]\times {{\bf R}}^{+})_\lambda$, the composition $u=h\circ \Phi_\Lambda$ is a smooth function on ${\A}$. We will perturb ${{\cs}}_f$ by adding $u$, i.e.,
${\css}={{\cs}}_f+u$. Denote the gradient of ${\css}$ by ${\s}$. The 
following lemma is standard (see \cite{T2}).

\slm
\begin{enumerate}
\item $\nabla s^f$ and $\nabla{\s}$ are continuous families of Fredholm 
operators from bundle $T{\B}^{\ast}$ to ${\E}$ over ${\B}^\ast$, and 
$\nabla s^f - \nabla{\s}$ is compact.
\item ${{\M}}={\s}^{-1}(0)$ can be represented by smooth sections.
\item There exists a constant $\epsilon >0$ such that when $|dh|<\epsilon$,
${{\M}}$ is compact.
\item When $|dh|\rightarrow 0$, the distance between ${{\M}_f}$ and ${{\M}}$ goes to zero.
\end{enumerate}
\elm

Next we define a section $G$ of the bundle ${\E}$ over ${\B}^{\ast}\times V$
where $V$ is the dual of the vector space
$\prod_{\lambda\in \Lambda}({\bf R}\times {\bf R})_\lambda$:
$$ 
G((A,\psi),(v,w)_\lambda)=
s^f(A,\psi)+grad(\rho(\Phi_\Lambda)(\sum_{\lambda\in \Lambda}
(v_\lambda p_\lambda+w_\lambda q_\lambda)))(A,\psi).
$$
Here the set $\Lambda$ of thickened loops satisfies the conditions in
Lemma 1.5.3, and $\rho$ is a cutoff function on $\prod_{\lambda\in \Lambda}
({{\bf R}}\times {{\bf R}})_\lambda$ satisfying that $\rho\equiv 0$ in a neighborhood
of $\prod_{\lambda\in \Lambda}([-1,1]\times \{0\})_\lambda$ and 
$\Phi_\Lambda([(A,\psi)])$ where $[(A,\psi)]\in {\B}^{\ast}$ is a 
non-degenerate critical
point of ${{\cs}}_f$, and $\rho\equiv 1$ in a neighborhood of the rest of
$\Phi_\Lambda({{\M}_f}^{\ast})$.

\slm
There exists $\epsilon>0$ (depending on $\rho$) such that $G$ is transversal
to the zero section of ${\E}$ when restricted to ${\B}^{\ast}\times
B(\epsilon)$, where $B(\epsilon)$ is a ball of radius $\epsilon$ centered
at the origin in $V=(\prod_{\lambda\in \Lambda}({{\bf R}}\times {{\bf R}})_\lambda)^\ast$.
\elm

\noindent{\bf Proof:} $G$ is transversal to the zero section of ${\E}$ over 
${{\M}_f}^{\ast}\times \{0\}$ by the choice of the set $\Lambda$.
By continuity and Lemma 1.5.4 (4), this lemma is proved.
\hfill ${\quad \Box}$

\noindent{\bf The Proof of Proposition 1.1.7:}

Apply Sard-Smale theorem to the projection $\Pi
: G^{-1}(0)\rightarrow B(\epsilon)$. For a generic $(v,w)_\lambda\in B(\epsilon)$,
the perturbation ${\css}={\cs}_f+u$ is admissible where 
$$
u=\rho(\Phi_\Lambda)(\sum_{\lambda\in \Lambda}(v_\lambda p_\lambda+w_\lambda q_\lambda)).
$$
\hfill ${\quad \Box}$

\chapter{Seiberg-Witten Equations on Cylindrical End Manifolds}
 
Throughout this chapter, we assume that $Y$ is an oriented 3-manifold with boundary which is the complement of a tubular neighborhood of a knot in an integral homology 3-sphere (many results proved in this chapter hold for general
3-manifolds with toroidal boundary). 
Equip $Y$ with a Riemannian metric $g_0$ 
such that a neighborhood of $\partial Y=T^2$ 
is orientedly isometric to $(-1,0]\times{\bf R}/2\pi{\bf Z}\times{\bf R}/2\pi{\bf Z}$. We attach $[0,\infty)\times T^2$ to $Y$ and still denote it by $Y$. 

Given a spin structure of $Y$, there is a unique $SU(2)$
vector bundle $W$ over $Y$ such that the oriented volume form acts 
on $W$ as identity by the Clifford multiplication. The spinor bundle $W$ is cylindrical, i.e.
on $[0,\infty)\times T^2$, $W$ is isometric to 
the pull back $\pi^{\ast} W_0$ where $\pi:
[0,\infty)\times T^2 \rightarrow T^2$ is the projection and $W_0$ is the total spinor bundle on $T^2$ associated to the spin structure induced from $Y$.

\section{The Fredholm theory}

In this section, we set up the Fredholm theory for Seiberg-Witten
equations on $Y$. 
Throughout ${\hh}^1(T^2)$ stands for the space of harmonic 1-forms on $T^2$. We 
fix a cut-off function $\rho$ on $Y$ which equals to $0$ on $Y\setminus
[0,\infty)\times T^2$ and $1$ on $[1,\infty)\times T^2$. 

\sde
For $\delta>0$, let
\begin{eqnarray*}
{\A}_{\delta}
&=&\{(A,\psi)| A=B+\rho\pi^{\ast}a,\\
& & B\in \oo, \psi\in \w, a\in \h\},
\end{eqnarray*}
where $\pi:[0,\infty)\times T^2\rightarrow T^2$ is the projection.
\ede

Here $L_{k,\delta}^2$ denotes the weighted Sobolev spaces with weight $\delta$
(\cite{LM}).
${\A}_{\delta}$ is a Hilbert space (over real numbers) with the norm  
$$
\|(A,\psi)\|_{{\A}_{\delta}}=\|(B,\psi)\|_{L^2_{2,\delta}}+\|a\|_{L^2}.
$$
Note that the decomposition of $A$ as $B+\rho\pi^{\ast}a$ is unique. We define
a map $R:{\A}_{\delta}\rightarrow\h$ by $R(A,\psi)=a$.

\sde
The group of gauge transformations is 
\begin{eqnarray*}
{\G}_{\delta}
&=& \{s\in L_{3,loc}^2(Y,S^1)|s^{-1}ds=g+\rho\pi^{\ast}h,\\
& & g\in \oo, h\in \h\}.
\end{eqnarray*}
${\G}_{\delta}$ acts on ${\A}_{\delta}$ by the formula $s\cdot (A,\psi)=(A-s^{-1}ds,s\psi)$ for $s\in{\G}_{\delta}$ and $(A,\psi)\in{\A}_{\delta}$.
\ede

\slm
${\G}_{\delta}$ is an Abelian Hilbert Lie group acting smoothly on ${\A}_{\delta}$ with the Lie algebra $T{\G}_{\delta,id}=\tz\oplus i{\bf R}$.
Moreover, for $s\in {\G}_{\delta}$, if $s^{-1}ds$ is decomposed as $g+\rho\pi^{\ast}h$,
then $h$ has zero period along the longitude and periods in $2\pi i{\bf Z}$ along the
meridian.
\elm

\noindent{\bf Proof:} 
Suppose that $s\in {\G}_{\delta}$ is in the component of identity, then $s=e^f$ for some $f\in L_{3,loc}^2(\Lambda^0(Y)\otimes i{\bf R})$. 
By Definition 2.1.2,
$df=s^{-1}ds$ can be decomposed as $g+\rho\pi^{\ast}h$, from which it follows that $h=0$ and $df\in\oo$. By Taubes inequality (Lemma 5.2 in \cite{T1}), there exists an imaginary valued constant $f_0$ on $Y$ such that
$$
\int_Y|f-f_0|^2e^{2\delta t}\leq C(\delta)\int_Y|df|^2e^{2\delta t},
$$
which proves that the Lie algebra $T{\G}_{\delta,id}$ is $\tz \oplus i{\bf R}$.

Let $\gamma_1,\gamma_2$ be the longitude and meridian, and $F$ be the Seifert
surface that $\gamma_1$ bounds in $Y$. For $s\in {\G}_{\delta}$, if $s^{-1}ds$ is
decomposed as $g+\rho\pi^{\ast}h$, then we have 
$$
\int_{\gamma_1}h=\int_F d(s^{-1}ds)=0 \hspace{2mm} \mbox{and} \hspace{2mm} 
\int_{\gamma_2}h=\int_{\gamma_2}s^{-1}ds\in 2\pi i{\bf Z}.
$$
The rest follows easily from the Sobolev theorems for weighted
spaces.
\hfill ${\Box}$

\slm
The de Rham cohomology group $H^1_{DR}(Y)$ can be represented by the space of
``bounded'' harmonic forms
$$
{\hh}^1(Y)=\{a\in\Omega^1(Y)| da=d^{\ast}a=0, \|a\|_{C^0(Y)}<\infty, 
\lim_{t \rightarrow \infty}a(\frac{\partial}{\partial t})=0 \}.
$$
Moreover, each element $a\in{\hh}^1(Y)$ can be decomposed as $b+\rho\pi^{\ast}
a_\infty$ with $a_\infty\in {\hh}^1(T^2)$ and $b\in L_{k,\delta}^2$ for some $\delta > 0$. The map $R: {\hh}^1(Y)\rightarrow {\hh}^1(T^2)$ defined by
$R(a)=a_\infty$ represents the embedding $H^1_{DR}(Y)\rightarrow H^1_{DR}(T^2)$.
As a corollary, for any $\kappa\in H^1(Y,{\bf Z})$, there is an $s_\kappa
\in C^\infty(Y,S^1)$ such that $s_\kappa^{-1}ds_\kappa\in {\hh}^1(Y)\otimes i{\bf R}$ and $[s_\kappa^{-1}ds_\kappa]=2\pi i\kappa$. So $\pi_0({\G}_\delta)=H^1(Y,{\bf Z})={\bf Z}$.
\elm

\noindent{\bf Proof:}
The Laplacian $\dd^\ast\dd: L^2_2(\Lambda^0(T^2))\rightarrow
L^2(\Lambda^0(T^2))$ restricted to $(\ker\dd^\ast\dd)^\perp$
is invertible.  Let $G$ be the inverse.
Suppose a closed form $A\in \Omega^1(Y)$ is written as $A_0 dt+A_1$ on 
$[0,\infty)\times T^2$ with $A_1\in \Omega^1(T^2)$.
Then $f=G(\dd^\ast A_1)$ is a smooth function on $[0,\infty)\times T^2$.
We extend $f$ to the rest of $Y$ and 
still call it $f$. Let $B=A-df$. We can further modify $f$ by a function of $t$ so that $\int_{T^2} B_0=0$, where $B=B_0 dt+B_1$ on $[0,\infty)\times T^2$.
$(B_0,B_1)$ satisfies the following equations:
$$
\frac{\partial B_1}{\partial t}=\dd B_0, \hspace{2mm} \dd B_1=0 \hspace{2mm}
\mbox{and} \hspace{2mm} \dd^\ast B_1=0,
$$
which shows that $B_1$ is in ${\hh}^1(T^2)$ and constant in $t$ and $B_0=0$. Since $d^\ast B$ is compactly supported, there is a unique solution $g\in L_{k,\delta}^2(\Lambda^0(Y))$ to the equation
$d^\ast B=d^\ast dg$ (see Lemma 2.1.7 below). Let $C=B-dg$, then $C\in{\hh}^1(Y)$ and the cohomology classes $[A]$ and $[C]$ are equal in $H^1_{DR}(Y)$.

Suppose $a\in {\hh}^1(Y)$ and $a=a_0 dt + a_1$ on $[0,\infty)\times T^2$. Then the pair $(a_0,a_1)$ satisfies the following system of equations
$$
\left\{\begin{array}{c}

\frac{\partial a_1}{\partial t}-\dd a_0 =0\\
\frac{\partial a_0}{\partial t}-\dd^\ast a_1 =0\\
\dd a_1=0.
\end{array} \right.
$$
The operator $L=\left(\begin{array}{ll}
0 & \dd \\
\dd^\ast & 0
\end{array} \right )$
is formally self-adjoint and elliptic on $\ker d\oplus \Omega^0(T^2)$.
By expanding $(a_1,a_0)$ in terms of an orthonormal basis of eigenvectors of $L$, we see that $a=a_0 dt + a_1$ can be decomposed as $b+\rho\pi^\ast a_\infty$
where $a_\infty\in{\hh}^1(T^2)$ and $b\in L_{k,\delta}^2$ for some
$\delta >0$.

Assume $a_1, a_2\in{\hh}^1(Y)$, if $a_1-a_2=df$ for a smooth function $f$ on $Y$,
then $df\in L_{k,\delta}^2$, and by Taubes inequality and integration by parts,  $df=0$. Hence the map ${\hh}^1(Y)\rightarrow H^1_{DR}(Y)$ is also injective.
The rest of the lemma follows easily.
\hfill ${\Box}$

\sde
Let ${\B}_{\delta}={\A}_{\delta}/{\G}_{\delta}$ and ${\B}^\ast_{\delta}={\A}^\ast_{\delta} /{\G}_{\delta}$ where ${\A}^\ast_{\delta}={\A}_{\delta} \setminus 
\{\psi\equiv 0 \}$.
\ede
 
\slm
\begin{enumerate}
\item ${\B}^\ast_\delta$ is a Hilbert manifold with the slice at $(A,\psi)\in {\A}^\ast_\delta$
given by $T_{(A,\psi),\epsilon}=U\times V$ where 
\begin{eqnarray*}
U=(B,\psi)+\{(a,\phi) & \in & \oo\oplus\w | -d^\ast a \\ 
 &   & + i\langle i\psi,\phi\rangle_{Re}=0, \|(a,\phi)\|_{L^2_{2,\delta}}<\epsilon \}, \\
\end{eqnarray*}
$$
V=R(A,\psi)+\{a_\infty\in{\hh}^1(T^2)\otimes i{\bf R} | \|a_\infty\|_{L^2} <\epsilon \},
$$
where $A$ is decomposed into $B+\rho\pi^\ast R(A,\psi)$. The tangent space of 
${\B}^\ast_\delta$ at $(A,\psi)$ is 
\begin{eqnarray*}
T{\B}^\ast_{\delta,(A,\psi)} & = & \{(a,\phi)\in \oo\oplus\w | -d^\ast a \\
&  & + i\langle i\psi,\phi\rangle_{Re}=0\}\oplus\h.
\end{eqnarray*}
\item A neighborhood of $[(A,0)]$ in ${\B}_\delta$ is diffeomorphic to $T_{(A,0),\epsilon}/S^1$. $T_{(A,0),\epsilon}=U\times V$ and
$$\begin{array}{c}
U=(B,0)+\{(a,\phi)\in \oo\oplus\w | d^\ast a =0,
\|(a,\phi)\|_{L^2_{2,\delta}}<\epsilon \}, \\
V=R(A,\psi)+\{a_\infty\in{\hh}^1(T^2)\otimes i{\bf R} | \|a_\infty\|_{L^2} <\epsilon \},
\end{array}
$$
where $A$ is decomposed into $B+\rho\pi^\ast R(A,\psi)$.
The action of $S^1$ on $T_{(A,0),\epsilon}$ is given by the complex multiplication on the factor $\phi$.
\end{enumerate}
\elm

\slm
Let $L_1=d^\ast d$ and $L_2=d^\ast d + |\psi|^2$ where $\psi\in\w$. Then there
is $\delta_0>0$ such that
for $k\geq 2$ and any $\delta\in (0,\delta_0]$,  $L_1: 
L_{k,\delta}^2(\Lambda^0(Y))\rightarrow L_{k-2,\delta}^2(\Lambda^0(Y))$ is a
Fredholm operator of index $-1$. $\ker\;L_1=0$, and
the range of $L_1$ is the $L^2$-orthogonal complement of the space of constant functions. $L_2:L_{k,\delta}^2(\Lambda^0(Y))\oplus i{\bf R}\rightarrow
L_{k-2,\delta}^2(\Lambda^0(Y))$ is isomorphic if $\psi$ is not identically
zero.
\elm

\noindent{\bf Proof:}
The operator $L_1=
d^\ast d: L_{k,\delta}^2(\Lambda^0(Y))\rightarrow L_{k-2,\delta}^2(\Lambda^0(Y))$ is Fredholm of index $-1$ by Theorem 7.4 of 
\cite{LM}. $\ker\;L_1=0$ follows from integration by parts. From index counting it follows that the range of $L_1$ is the $L^2$-orthogonal complement of the space of constant functions. For $\psi\in\w$, $L_2: L_{k,\delta}^2(\Lambda^0(Y))\rightarrow L_{k-2,\delta}^2(\Lambda^0(Y))$ is a compact perturbation of $L_1$, so it is also
a Fredholm operator of index $-1$. So $L_2:L_{k,\delta}^2(\Lambda^0(Y))\oplus i{\bf R}\rightarrow L_{k-2,\delta}^2(\Lambda^0(Y))$ is an isomorphism if $\psi$ is not identically zero, since $\ker\;L_2=0$ and index $L_2=0$.
\hfill ${\Box}$

\noindent{\bf The Proof of Lemma 2.1.6:} 
\begin{enumerate}
\item The construction of a local slice is standard by applying the implicit function theorem. The key point is the properties of $L_2$ stated in Lemma 2.1.7.
To prove that ${\B}^\ast_\delta$ is Hausdorff and the local slice is embedded into ${\B}^\ast_\delta$,
the argument in \cite{FU} can be used, combined with Taubes inequality
(Lemma 5.2 in \cite{T1}).
\item Part $2$ of this lemma follows similarly with Lemma 2.1.7 understood.
\end{enumerate}
\hfill ${\Box}$

\sde
For $(A,\psi)\in{\A}_\delta$, we define
$$
{\E}_{\delta,(A,\psi)}=\{(a,\phi)\in \zo\oplus\wo | -d^\ast a + i\langle i\psi,\phi\rangle_{Re}=0 \}.
$$ 
${\E}_{\delta,(A,\psi)}$ is a closed subspace of $\zo\oplus\wo$.
\ede

\slm
${\E}_\delta=\{{\E}_{\delta,(A,\psi)}\}$ is a Hilbert bundle over ${\A}^\ast_\delta$ which descends to a Hilbert bundle over ${\B}^\ast_\delta$ (we still call it ${\E}_\delta$).
\elm

\noindent{\bf Proof:} 
For any $(a,\phi)\in\zo\oplus \wo$,
we can project $(a,\phi)$ into ${\E}_{\delta,(A,\psi)}$ by solving the following
equation
$$
-d^\ast(a-df) + i\langle i\psi,\phi +f\psi\rangle_{Re}=0 
$$
for $f\in\oz\oplus i{\bf R}$. By Lemma 2.1.7, the operator $L_2=d^\ast d+|\psi|^2$
is an isomorphism from $\oz\oplus i{\bf R}$ to
$L_{0,\delta}^2(\Lambda^0(Y)\otimes i{\bf R})$ since $(A,\psi)\in{\A}^\ast_\delta$.
So the above equation has a unique solution $f(a,\phi)$ for any $(a,\phi)$.
If $(a,\phi)\in{\E}_{\delta,(A_1,\psi_1)}$ with $(A_1,\psi_1)$ close enough to
$(A,\psi)$, one can easily show that the projection $(a,\phi)\rightarrow
(a-df,\phi+f\psi)$ is one to one and onto, again using the invertibility of $L_2$. This proves the local triviality of ${\E}_\delta$. The bundle ${\E}_\delta$ over ${\A}^\ast_\delta$ is ${\G}_\delta$-equivariant, so it descends to a Hilbert bundle over ${\B}^\ast_\delta$.
\hfill ${\Box}$

\sde
For $(A,\psi)\in{\A}^\ast_\delta$, we define
$$
s(A,\psi)=(\ast dA + \tau(\psi,\psi), D_A\psi).
$$
Here $D_A=D_{g_0}+A$ where $D_{g_0}$ is the Dirac operator associated to the metric $g_0$.
$s$ is a section of ${\E}_\delta$
over ${\A}^\ast_\delta$, which descends to a section of ${\E}_\delta$ over ${\B}^\ast_\delta$.

The covariant derivative of $s$ is a section of $End(T{\B}^\ast_\delta,{\E}_\delta)$ over ${\A}^\ast_\delta$ which descends to ${\B}^\ast_\delta$, defined by
$$
\nabla s_{(A,\psi)}(a,\phi)=
(\ast da+2\tau(\psi,\phi)-df(a,\phi), D_A\phi+a\psi+f(a,\phi)\psi)
$$
where $f(a,\phi)$ is the unique solution to the equation
$$
d^\ast df+f|\psi|^2=i\langle D_A\psi,i\phi\rangle_{Re}.
$$
The map $(a,\phi)\rightarrow (-df(a,\phi),f(a,\phi)\psi)$ from 
$T{\B}^\ast_{\delta,(A,\psi)}$ to
$\zo\oplus\wo$ is compact by the Sobolev theorems for weighted spaces.
\ede

\sde
\begin{enumerate}
\item For any $r>0$, let ${\hh}(r)=\h\setminus\bigcup_{p\in B}D(p,r)$ where $B$ is the lattice of ``bad'' points for the induced spin structure on $T^2$ (see Appendix A) and $D(p,r)$ is the closed disc of radius $r$ centered at $p$.
\item ${\A}_\delta(r)=R^{-1}({\hh}(r))$, ${\A}^\ast_\delta(r)={\A}_\delta(r)\bigcap{\A}^\ast_\delta$, ${\B}_\delta(r)={\A}_\delta(r)/{\G}_\delta$ and ${\B}^\ast_\delta(r)={\A}^\ast_\delta(r)/{\G}_\delta$, where $R:{\A}_\delta\rightarrow\h$ is given by $R(A,\psi)=a$ for $A=B+\rho\pi^\ast a$.
\end{enumerate}
\ede
Note that for any $a\in {\hh}(r)$, the twisted Dirac operator $D^{T^2}_a=D^{T^2}+a$ is invertible, where $D^{T^2}$ is the Dirac operator on $T^2$ (see Appendix A). 

\sprop
For any $r>0$, there exists a $\delta(r)>0$ such that for each
$\delta\in (0,\delta(r))$,
$\nabla s:T{\B}^\ast_\delta\rightarrow{\E}_\delta$ is a continuous family of Fredholm operators of index $1$
over ${\A}^\ast_\delta(r)$. (So $s$ is a Fredholm section of ${\E}_\delta$ over ${\B}^\ast_\delta(r)$).
\eprop

At $(A,\psi)\in{\A}^\ast_\delta$, we have a short exact sequence
$$
0\rightarrow T{\G}_{\delta,id} \stackrel{d_{(A,\psi)}}{\rightarrow} T{\A}^\ast_\delta\stackrel{\pi_\ast}{\rightarrow}T{\B}^\ast_\delta\rightarrow 0
$$
where $d_{(A,\psi)}(f)=(-df,f\psi)$ and $\pi:{\A}^\ast_\delta\rightarrow{\B}^\ast_\delta$ is the natural projection. This enables us to extend $\nabla s_{(A,\psi)}: 
T{\B}^\ast_{\delta,(A,\psi)}
\rightarrow{\E}_{\delta,(A,\psi)}$ to a ${\G}_\delta$-equivariant map ${\K}^\prime_{(A,\psi)}$ (see \cite{T2}), where
$$
{\K}^\prime_{(A,\psi)}=\left(\begin{array}{ccc}
\nabla s_{(A,\psi)} & 0 & 0\\
0 & 0 & d_{(A,\psi)}\\
0 & d^\ast_{(A,\psi)} & 0 
\end{array} \right).
$$
${\K}^\prime_{(A,\psi)}$ is from $T{\A}^\ast_\delta\oplus(\oz\oplus i{\bf R})$ to
$\zo\oplus\wo\oplus\zz$. Since the operator $\left(\begin{array}{ll}
0 & d_{(A,\psi)}\\
d^\ast_{(A,\psi)} & 0
\end{array} \right)$ is invertible, $\nabla s_{(A,\psi)}$ is Fredholm if and only if ${\K}^\prime_{(A,\psi)}$ is, and 
$\mbox{index}{\K}^\prime=\mbox{index}\nabla s$.

For $(A,\psi)\in{\A}_\delta$, we define a map ${\K}_{(A,\psi)}:{\A}_\delta\oplus(\oz\oplus i{\bf R})
\rightarrow\zo\oplus\wo\oplus\zz$ by 
$$
{\K}_{(A,\psi)}(a,\phi,f)=
(\ast da+2\tau(\psi,\phi)-df,D_A\phi+a\psi+f\psi,-d^\ast a+
i\langle i\psi,\phi\rangle_{Re}).
$$
Then ${\K}^\prime_{(A,\psi)}$ is a compact perturbation of ${\K}_{(A,\psi)}$ and Proposition 2.1.12 follows from 

\slm
For any $r>0$, there exists a $\delta(r)>0$ such that for each
$\delta\in (0,\delta(r))$, ${\K}$ is a continuous family of Fredholm maps on ${\A}^\ast_\delta(r)$ of index $1$.
\elm

\noindent{\bf Proof:} Consider the following commutative diagram
$$
\begin{array}{ccccccccc}
0 & \rightarrow & V{\A}_{\delta,1}\oplus L^2_{1,\delta} & \rightarrow & V{\A}_{\delta,1}\oplus L^2_{1,\delta} &
\rightarrow & 0 &  &  \\
  &  & \uparrow V{\K}_{(A,\psi)} &  & \uparrow{\K}_{(A,\psi)} &  &
\uparrow &  &  \\
0 & \rightarrow & V{\A}_\delta\oplus L^2_{2,\delta}& \rightarrow & {\A}_\delta\oplus(L^2_{2,\delta}\oplus i{\bf R}) & \rightarrow & \h\oplus i{\bf R} & \rightarrow & 0
\end{array}
$$
where $L^2_{1,\delta}=\zz$ and $L^2_{2,\delta}=\oz$.
$V{\A}_\delta$ is the fiber of map $R: {\A}_\delta\rightarrow\h$ and $V{\A}_{\delta,1}$ is its $L_{1,\delta}^2$-completion. $V{\K}_{(A,\psi)}$ is
the restriction of ${\K}_{(A,\psi)}$.
We have a long exact sequence (see \cite{M})
$$
\begin{array}{cccccccc}
  &  & CokerV{\K}_{(A,\psi)} & \rightarrow & Coker{\K}_{(A,\psi)} &%
\rightarrow & 0 & \\
0 & \rightarrow & KerV{\K}_{(A,\psi)} & \rightarrow & Ker{\K}_{(A,\psi)} & \rightarrow & \h\oplus i{\bf R} & \rightarrow. 
\end{array}
$$
The lemma follows from the claim that for
any $r>0$, there exists a $\delta(r)>0$ such that for each
$\delta\in (0,\delta(r))$,
$V{\K}_{(A,\psi)}:
V{\A}_\delta\oplus\oz \rightarrow V{\A}_{\delta,1}\oplus\zz $ is a Fredholm map of index $-2$ for $(A,\psi)\in{\A}^\ast_\delta(r)$.

$V{\K}_{(A,\psi)}$ is a compact perturbation of an operator of form $I(\frac{\partial}{\partial t} + B_a)$ on $[0,\infty)\times T^2$ where 
$$
I=\left(\begin{array}{clcr}
dt & 0 & 0 & 0 \\
0 & \ast_{T^2} & 0 & 0 \\
0 & 0 & 0 & -1 \\
0 & 0 & 1 & 0
\end{array} \right)
\hspace{2mm} \mbox{and} \hspace{2mm} 
B_a=\left(\begin{array}{clcr}
D^{T^2}_a & 0 & 0 & 0\\
0 & 0 & -\dd & \ast\dd \\
0 & -\dd^\ast & 0 & 0\\
0 & -\ast\dd & 0 & 0
\end{array} \right)
$$
acting on ${\Gamma}(W_0\oplus(\Lambda^1\oplus\Lambda^0\oplus\Lambda^0(T^2))\otimes
i{\bf R})$. Here $W_0$ is the total spinor bundle over $T^2$, and $D^{T^2}_a$ 
is the
twisted Dirac operator with $a=R(A,\psi)\in\h$.
For any $r>0$, let $\delta(r)=\min\{|u|: u\neq 0$ is an 
eigenvalue of $B_a$ for some $a\in{\hh}(r)\}$. Then
$V{\K}_{(A,\psi)}:
V{\A}_\delta\oplus\oz \rightarrow V{\A}_{\delta,1}\oplus\zz $ is a Fredholm map of index $-2$ for any $\delta\in (0,\delta(r))$ by Theorem 7.4 of \cite{LM}.
\hfill ${\Box}$

\section{Perturbation and transversality}

Fix a small $r>0$ and a $\delta\in (0,\delta(r))$ for the weight of the Sobolev spaces. For simplicity we omit the subscript $\delta$ in the discussion. Consider the following perturbations of $s$ over ${\A}^\ast(r)$:
$$
s_{\mu,f}(A,\psi)=(\ast dA+\tau(\psi,\psi)+\mu, D_A\psi+f\psi)
\hspace{2mm} \mbox{for} \hspace{2mm} (A,\psi)\in{\A}^\ast(r),
$$
where $\mu$ is a co-closed imaginary valued 1-form
and $f$ is a smooth real valued function on $Y$, both supported in
$Y\setminus [0,\infty)\times T^2$. The metric $g$ being used here is a perturbation of $g_0$ supported in $Y\setminus [0,\infty)\times T^2$.

\sde 
Define the Seiberg-Witten moduli spaces 
$$\begin{array}{c}
{\M}_{\mu,f}(r)=\{[(A,\psi)]\in{\B}(r) | (\ast dA+\tau(\psi,\psi)+\mu, D_A\psi+f\psi)=0\},\\
 {\M}_{\mu,f}^\ast(r)={\M}_{\mu,f}(r)\bigcap{\B}^\ast.
\end{array}
$$
\ede

Let $[R]:{\B}(r)\rightarrow{\hh}(r)/{\bf Z}$ be the map induced by $R:{\A}\rightarrow\h$. 

\sprop
The moduli spaces ${\M}_{\mu,f}(r)$ and ${\M}_{\mu,f}^\ast(r)$ have the following properties.
\begin{enumerate}
\item For a generic $\mu$, ${\M}_{\mu,f}^\ast(r)$ is a collection of 
1-dimensional smooth curves, and the map
$[R]:{\M}_{\mu,f}^\ast(r)\rightarrow{\hh}(r)/{\bf Z}$ is an immersion.
\item Given a set $S$ of immersed curves in ${\hh}(r)/{\bf Z}$, for a generic $\mu$, the map $[R]:{\M}_{\mu,f}^\ast(r)\rightarrow{\hh}(r)/{\bf Z}$ is transversal to $S$.
\item For a generic $(g,f)$, the $L^2$-closed
extension of the perturbed Dirac operator $D_g+f$ is invertible. Fix such a $(g,f)$, then for any small enough $\mu$, there exists a neighborhood $U_\mu$ of $[(a_\mu,0)]$ in ${\B}(r)$ such that $U_\mu\bigcap{\M}_{\mu,f}^\ast(r)=\emptyset$,
where $a_\mu$ is the unique
solution to $\ast da_\mu+\mu=0$, $d^\ast a_\mu=0$ such that 
$R(a_\mu)$ has zero period along the meridian.
\item ${\M}_{\mu,f}(r)\setminus{\M}_{\mu,f}^\ast(r)=\{(a_\mu+iA,0)|A\in{\hh}^1(Y)/{\hh}^1(Y,{\bf Z})\}\simeq S^1$. Note that for any $A\in {\hh}^1(Y)$, $R(A)$ is a multiple of $dy$ where $e^{iy}$ parameterizes the meridian.
\end{enumerate}
\eprop

For simplicity we omit the subscript $f$. Consider the section $\tilde{s}$ of ${\E}$ over ${\B}^\ast(r)\times \ker\;d^\ast$:
$$
\tilde{s}([(A,\psi)],\mu)=[(\ast dA+\tau(\psi,\psi)+\mu, D_A\psi+f\psi)].
$$
For any $([(A,\psi)],\mu)\in \tilde{s}^{-1}(0)$, we have the following commutative
diagram:
$$
\begin{array}{ccccccccc}
0 & \rightarrow & {\E} & \rightarrow & {\E} & \rightarrow & 0 &  & \\
  &  & \uparrow V\nabla \tilde{s} &  & \uparrow \nabla \tilde{s} &   &
\uparrow 0 &  &  \\
0 & \rightarrow & VT{\B}^\ast(r)\times \ker\;d^\ast & \rightarrow & T{\B}^\ast(r)\times \ker\;d^\ast & \stackrel{d[R]}{\rightarrow} & \h & \rightarrow & 0
\end{array}
$$
at $([(A,\psi)],\mu)$. This gives rise to a long exact sequence (see \cite{M})
$$
\begin{array}{cccccccc}
  &  & Coker(V\nabla \tilde{s}_{([(A,\psi)],\mu)}) & \rightarrow & Coker(\nabla \tilde{s}_{([(A,\psi)],\mu)}) 
& \rightarrow 0 &  & \\
0 & \rightarrow & Ker(V\nabla \tilde{s}_{([(A,\psi)],\mu)}) & \rightarrow & 
Ker(\nabla \tilde{s}_{([(A,\psi)],\mu)})
& \stackrel{d[R]}{\rightarrow} & \h & \rightarrow. 
\end{array}
$$

\slm
$Coker(V\nabla \tilde{s}_{([(A,\psi)],\mu)})=0$ for any $([(A,\psi)],\mu)\in \tilde{s}^{-1}(0)$.
\elm

\noindent{\bf Proof:}  
First observe that $V\nabla \tilde{s}_{([(A,\psi)],\mu)}$ is Fredholm as a map from the $L_{1,\delta}^2$-completion of $VT{\B}^{\ast}_{[(A,\psi)]}$ to the $L_{0,\delta}^2$-completion
of ${\E}_{[(A,\psi)]}$. So by regularity, it suffices to show that the $L_{0,\delta}^2$-orthogonal complement of the image of $V\nabla \tilde{s}_{([(A,\psi)],\mu)}$ is
zero dimensional.

Let $(a^\prime,\phi^\prime)
\in {\E}_{[(A,\psi)]}$ be $L^2_{0,\delta}$-orthogonal to the range of $V\nabla \tilde{s}_{([(A,\psi)],\mu)}$. Set $(a,\phi)=e^{2\delta t}(a^\prime,\phi^\prime)$, then 
$(a,\phi)$ is $L^2$-orthogonal to the range of $V\nabla \tilde{s}_{([(A,\psi)],\mu)}$ and
$e^{-2\delta t}(a,\phi)$ is in ${\E}_{[(A,\psi)]}$, i.e.
$
-d^\ast(e^{-2\delta t}a)+i\langle i\psi, e^{-2\delta t}\phi\rangle_{Re}=0.
$
Note that $(a,\phi)$ is in $L^2_{1,-\delta}$.

Observe that $\oo\oplus\w=VT{\B}^\ast(r)\oplus Im(d_{[(A,\psi)]})$ (recall
the map $d_{(A,\psi)}$ is defined by $f\rightarrow (-df,f\psi)$), and
$V\nabla \tilde{s}_{([(A,\psi)],\mu)}$ vanishes on $Im(d_{[(A,\psi)]})$. So if $(a,\phi)$
is $L^2$-orthogonal to the range of $V\nabla \tilde{s}_{([(A,\psi)],\mu)}$, then for any $\mu^\prime\in\ker\;d^\ast$ and $(b,\theta)\in\oo\oplus\w$, we have
$$
(\ast db+2\tau(\psi,\theta)+\mu^\prime, a)+ (D_A\theta+f\theta+b\psi,\phi)=0.
$$
This implies that $D_A\phi+f\phi+a\psi=0$ and $(b\psi,\phi)=0$ for any $b$.
Since $\psi$ is not identically zero, by the unique continuation theorem for
Dirac operators, we have $\phi=ih\psi$ for a real valued function $h$. Then
$D_A\phi+f\phi+a\psi=0$ implies that $idh+a=0$. Hence
$$
d^\ast(e^{-2\delta t}dh)+he^{-2\delta t}|\psi|^2=0.
$$
That $(a,\phi)\equiv 0$ follows from $h\equiv 0$, which follows by integration by parts from the claim that $e^{-\delta t}|h|$ is bounded on 
$[0,\infty)\times T^2$.

Next we prove that $e^{-\delta t}|h|$ is bounded on $[0,\infty)\times T^2$.
First of all, $idh+a=0$ implies that $\frac{\partial h}{\partial t}\in
L^2_{1,-\delta}$ and $\dd h\in L^2_{1,-\delta}$. Let $h_0(t)$ be the $L^2$-orthogonal projection of $h$ onto $\ker\dd^\ast\dd$, then
$\|h-h_0(t)\|_{L^2(T^2)}\leq c\|\dd h\|_{L^2(T^2)}$. So $h-h_0(t)\in L^2_{0,-\delta}$. On the other hand, $|\frac{dh_0}{dt}|\leq C\|
\frac{\partial h}{\partial t}\|_{L^2(T^2)}$ so that $\frac{dh_0}{dt}
e^{-\delta t}$ is bounded on $[0,\infty)\times T^2$. So
$$
|h_0(t)-h_0(0)|\leq \int_{0}^{t}|\frac{dh_0}{dt}|dt\leq \frac{C}{\delta}
e^{\delta t}.
$$
It follows easily that $e^{-\delta t}|h|$ is bounded on $[0,\infty)\times T^2$.
\hfill ${\Box}$

\slm
Given any spin structure on $Y$, for a generic $(g,f)$, the $L^2$-closed extension of the perturbed Dirac operator $D_g+f$ is invertible.
\elm

\noindent{\bf Proof:} For a perturbed metric $g$ of $g_0$ which is supported in $Y\setminus [0,\infty)\times T^2$, the Dirac operator $D_g$ on $Y$ takes the form of $dt(\frac{\partial}{\partial t}+D^{T^2})$ on the cylindrical end. Note that $D^{T^2}$ is invertible (see Appendix A). So the $L^2$-closed extension of $D_g+f$ is an essentially self-adjoint Fredholm operator ($f$ is vanishing on the cylindrical end). The argument for the proof of Proposition 1.1.6 can be applied to prove this lemma.
\hfill ${\Box}$

\slm
For small enough $\delta > 0$, the operator
$\ast d: L^2_{k,\delta}(\Lambda^1(Y))\bigcap \ker\;d^\ast\rightarrow
L^2_{k-1,\delta}(\Lambda^1(Y))\bigcap \ker\;d^\ast$
is Fredholm with $\dim\ker\ast d=0$ and $\dim Coker\ast d=1$. Moreover,
for any compactly supported co-closed 1-form $\mu$, there exists a unique
$a_\mu\in\ker\;d^\ast$ such that $i)$
$\ast da_\mu+\mu=0$; $ii)$ $a_\mu$ can be decomposed as $b+\rho\pi^\ast a_\infty$ where $b\in L^2_{k,\delta}$ and $a_\infty\in {\hh}^1(T^2)$ with zero period along the meridian. Furthermore,
$a_\mu$ satisfies the estimate:
$\|b\|_{L^2_{k,\delta}}+|a_\infty|\leq C\|\mu\|_{L^2_{k-1,\delta}}$.
\elm

\noindent{\bf Proof:}
The Fredholm property and the index calculation of $\ast d$ follows from a
similar argument as in Proposition 2.1.12. $\ker\;\ast d=0$ follows from
$H^1(Y,T^2)=0$. Given $\mu\in\ker\;d^\ast$ or equivalently $\ast\mu\in\ker\;d$,
since $H^2(Y,{\bf R})=0$, there exists an $A\in\Omega^1(Y)$ such that 
$dA+\ast\mu=0$ or equivalently $\ast dA+\mu=0$. If $\mu$ is compactly supported,
the argument for Lemma 2.1.4 can be used to modify $A$ with an
exact 1-form and a ``bounded'' harmonic form, and the resulting 1-form
$a_\mu$ has the claimed properties.
\hfill ${\Box}$

\noindent{\bf The Proof of Proposition 2.2.2:}

Since $Coker(\nabla \tilde{s}_{([(A,\psi)],\mu)})=0$ for any $([(A,\psi)],\mu)\in \tilde{s}^{-1}(0)$,
$\tilde{s}^{-1}(0)$ is a Banach manifold. The projection $\Pi:\tilde{s}^{-1}(0)\rightarrow
\ker\;d^\ast$ is a Fredholm map of index $1$ (Proposition 2.1.12). So
by Sard-Smale theorem, for a generic $\mu$, ${\M}^\ast_\mu(r)=\Pi^{-1}(\mu)$ is a collection of 1-dimensional smooth curves. In addition, $\ker\;(V\nabla \tilde{s})\bigcap\ker\;\Pi=0$
since $V\nabla \tilde{s}$ is formally self-adjoint on $VT{\B}^\ast(r)$. So $d[R]:
T{\M}^\ast_\mu(r)\rightarrow \h$ is injective.

Since $Coker(V\nabla \tilde{s}_{([(A,\psi)],\mu)})=0$ for any $([(A,\psi)],\mu)\in \tilde{s}^{-1}(0)$,
the map $[R]:\tilde{s}^{-1}(0)\rightarrow {\hh}(r)/{\bf Z}$ is a submersion. For any set $S$ of
immersed curves in ${\hh}(r)/{\bf Z}$, $[R]^{-1}(S)$ is a set of immersed submanifolds of co-dimension $1$ in $\tilde{s}^{-1}(0)$. If $\mu$ is a regular value of the
projection $\Pi:[R]^{-1}(S)\rightarrow\ker\;d^\ast$, then the map $[R]:
{\M}^\ast_\mu(r)\rightarrow {\hh}(r)/{\bf Z}$ is transversal to $S$.

Properties $3$, $4$ follow easily from Lemmas $2.2.4$, $2.2.5$ and $2.1.4$.

\section{The finite energy monopoles}

Fix a perturbation $(g,f,\mu)$  which is supported in $Y\setminus [0,\infty)\times T^2$.

\sde
$(A,\psi)\in \Omega^1(Y)\otimes i{\bf R}\oplus{\Gamma}(W)$ is said
to be {\bf a monopole of finite energy} if $(A,\psi)$ satisfies
\begin{itemize}
\item the Seiberg-Witten equations 
$$
\left\{\begin{array}{c}
\ast dA+\tau(\psi,\psi)+\mu=0\\ D_g\psi+ A\psi +f\psi=0;
\end{array} \right.
$$
\item the finite energy condition
$$
\int_Y(|\nabla_A\psi|^2+\frac{1}{2}|\psi|^4) < \infty.
$$
\end{itemize}
\ede

\noindent{\bf The exponential decay estimates}

\slm
(Lemma $4$ in \cite{KM1})

Let $X$ be a compact 3-manifold with boundary. Assume that $(A,\psi)\in
\Omega^1(X)\otimes i{\bf R}\oplus{\Gamma}(W)$ satisfies the Seiberg-Witten equations 
on $X$. Then there exists a gauge transformation $s\in C^\infty(X,S^1)$ such
that for any sub-domain $X^\prime$ with $\overline{X^\prime}\subset {int}X$,
$s\cdot (A,\psi)$ satisfies:
$$\begin{array}{c}
\|s\cdot (\psi)\|_{C^k(X^\prime)}\leq C(k,X,X^\prime)h_1(\|\psi\|_{L^4(X)}),\\
\|s\cdot (A)\|_{C^k(X^\prime)}\leq C(k,X,X^\prime)h_2(\|\psi\|_{L^4(X)})
\end{array}
$$
for a constant $C(k,X,X^\prime)$ and polynomials $h_1,h_2$ with $h_1(0)=0$.
\elm

\scor
For a finite energy monopole $(A,\psi)$,
$$
\|\psi\|_{C^0(T^2)}(t)
\leq C \|\psi\|_{L^4([t-1,t+1]\times T^2)}.$$
In particular, $\psi\rightarrow 0$ as $t\rightarrow \infty$. Moreover, there exists a constant $K$ depending only on the geometry of $Y$ and the norm of  $(\mu,f)$ such that $\|\psi\|_{C^0(Y)} < K$.
\ecor

\noindent{\bf Proof:}
It follows from Lemma 2.3.2, the Weitzenb\"ock formula and maximum principle.
\hfill ${\Box}$

Throughout this section, we use $a(t)$ to denote the harmonic component of $A_1$
where $A=A_0 dt+A_1$ on $[0,\infty)\times T^2$.
After a gauge transformation, any $(A,\psi)$ takes the {\bf standard form} on $[0,\infty)\times T^2$, i.e. $\dd^\ast A_1=0$ and $\int_{T^2}A_0=0$ (see Lemma 2.1.4).

\slm
Assume that the finite energy monopole $(A,\psi)$ is in the standard form. Then
the following holds for a constant $c$: 
\begin{itemize}
\item [{a)}] $\int_{T^2}(|A_0|^2+|\dd A_0|^2)\leq c\int_{T^2}|\psi|^4$;
\item [{b)}] $\int_{T^2}(|A_1-a(t)|^2+|\nabla^{T^2}(A_1-a(t))|^2)
\leq c\int_{T^2}|\psi|^4$;
\item [{c)}] $\int_{T^2}|\frac{\partial}{\partial t}(A_1-a(t))|^2
\leq c\int_{T^2}|\psi|^4$;
\item [{d)}] $\int_{T^2}|\frac{\partial}{\partial t}a(t)|^2
\leq c\int_{T^2}|\psi|^4$;
\item [{e)}] $\|\frac{\partial}{\partial t}A_0\|_{L^2_1(T^2)}
\leq c\|\frac{\partial\psi}{\partial t}\|_{L^2(T^2)}$.
\end{itemize}
\elm

\noindent{\bf Proof:}
On $[0,\infty)\times T^2$, the equation $\ast dA+\tau(\psi,\psi)=0$
reads as 
$$\begin{array}{c}
\dd(A_1-a(t))+q_1(\psi)=0,\\ 
\frac{\partial}{\partial t}(A_1-a(t))+
\frac{\partial}{\partial t}a(t)-\dd A_0 +q_2(\psi)=0
\end{array}
$$
for some quadratic forms $q_1,q_2$. Observe that $
\frac{\partial}{\partial t}(A_1-a(t)),
\frac{\partial}{\partial t}a(t)$ and $\dd A_0$ are $L^2$ orthogonal to each other.
The estimates $a),b),c),d)$ follow easily.

For $e)$, note that $\dd^\ast\dd(\frac{\partial}{\partial t}A_0)=
\dd^\ast(\frac{\partial}{\partial t}q_2(\psi))$. So we have 
$$
\|\frac{\partial}{\partial t}A_0\|_{L^2_1(T^2)}\leq c\|
\dd^\ast\dd(\frac{\partial}{\partial t}A_0)\|_{L^2_{-1}(T^2)}\leq c\|
\frac{\partial}{\partial t}q_2(\psi))\|_{L^2(T^2)}\leq c\|
\frac{\partial\psi}{\partial t}\|_{L^2(T^2)}
$$ 
since $\int_{T^2}A_0=0$ and $\|\psi\|_{C^0(Y)} < K$.
\hfill ${\Box}$

\slm
For any $r>0$, there exists a $c(r)>0$ such that $c(r)\int_{T^2}|\psi|^2\leq\int_{T^2}|D^{T^2}_a\psi|^2$ for any $a$ in the closure of ${\hh}(r)$ (see Definition 2.1.11 for ${\hh}(r)$).
\elm

\noindent{\bf Proof:}
Observe that both $\int_{T^2}|\psi|^2$ and $\int_{T^2}|D^{T^2}_a\psi|^2$ are
gauge invariant, so we can assume that $a$ is in the compact set $\overline{{\hh}(r)}/({\bf Z}\oplus{\bf Z})$. The lemma follows by taking $c(r)=\min\{u^2: u$ is an eigenvalue of $D^{T^2}_a$ for some $a$ in $\overline{{\hh}(r)}/({\bf Z}\oplus{\bf Z})\}$.
\hfill ${\Box}$

The following estimate turns out to be crucial.

\slm
There exists a constant $c_1$ with the following significance.
Let $(A,\psi)$ be a finite energy monopole.
For any $r>0$, if $a(t)$ is in ${\hh}(r)$ for $T_1<t<T_2$, then 
$$
\int_{T_1}^{T_2}\int_{T^2}|\psi|^2 \leq \frac{c_1}{c(r)}
\int_{T_1}^{T_2}\int_{T^2}(|\nabla_A\psi|^2+\frac{1}{2}|\psi|^4).
$$
\elm

\noindent{\bf Proof:} Assume that $(A,\psi)$ is in the standard form without
loss of generality. 
Since $a(t)$ is in ${\hh}(r)$ for $T_1<t<T_2$, by Lemma 2.3.5, for $T_1<t<T_2$, 
we have
\begin{eqnarray*}
c(r)\int_{T^2}|\psi|^2(t)
& &\leq\int_{T^2}|D^{T^2}_{a(t)}\psi|^2(t)\leq\int_{T^2}
|\nabla^{T^2}_{a(t)}\psi|^2(t)\\
& &\leq\int_{T^2}(|\nabla_A\psi|^2+|(A_1-a(t))\otimes\psi|^2)(t).
\end{eqnarray*}
But 
$\int_{T^2}|(A_1-a(t))\otimes\psi|^2
 \leq  K^2\int_{T^2}|(A_1-a(t))|^2 
 \leq C\int_{T^2}|\psi|^4$
by Corollary 2.3.3 and Lemma 2.3.4 b). So we have 
$$
\int_{T_1}^{T_2}\int_{T^2}|\psi|^2\leq \frac{c_1}{c(r)}
\int_{T_1}^{T_2}\int_{T^2}(|\nabla_A\psi|^2+\frac{1}{2}|\psi|^4)
$$
for a constant $c_1$.
\hfill ${\Box}$

\slm
Let $\gamma$ be a loop in $T^2$.
Then there exists a constant $c(\gamma)$ such that for any  $(A,\psi)$ 
satisfying the Seiberg-Witten equations on the cylindrical end, the following estimate holds for any $t_1<t_2$:
$$
\int_{t_1}^{t_2}\int_{\gamma}|\psi|^2\leq c(\gamma)(\int_{t_1}^{t_2}\int_{T^2}|\psi|^2+\int_{t_1-1}^{t_2+1}\int_{T^2}
(|\nabla_A\psi|^2+\frac{1}{2}|\psi|^4)).
$$
\elm

\noindent{\bf Proof:}
Note that it suffices to prove the estimate for $t_2=t_1+1$. Also note that both sides of the estimate are gauge invariant.
By the embedding $L^2_1(T^2)\rightarrow L^2(\gamma)$, we have
$$
\int_{\gamma}|\psi|^2\leq C\int_{T^2}(|\nabla^{T^2}\psi|^2+|\psi|^2)
\leq C\int_{T^2}(|\nabla_A \psi|^2+|A\otimes\psi|^2+|\psi|^2).
$$
On the other hand, in $U=[t_1-1,t_1+2]\times T^2$, $A$ can be decomposed into
$A=B+h$ in a Hodge gauge (Lemma 4 in \cite{KM1}) such that $\|B\|_{L^2_1(U)}
\leq C\|dA\|_{L^2(U)}\leq C_1\|\psi\|^2_{L^4(U)}$ and $h$ is harmonic with norm bounded by $K$. Hence 
$$\int_{t_1}^{t_1+1}\int_{T^2}|A\otimes\psi|^2\leq K^2
\int_{t_1}^{t_1+1}\int_{T^2}|\psi|^2 + \|B\|_{L^2_1(U)}^2\cdot\|\psi\|^2_{L^4(U)}.
$$ 
The lemma follows easily from these estimates.
\hfill ${\Box}$

\slm
There exists a constant $c_2$ such that the following estimate
\begin{eqnarray*}
|a(t_1)-a(t_2)| & \leq & c_2(\int_{t_1}^{t_2}\int_{T^2}|\psi|^2+
\int_{t_1-1}^{t_2+1}\int_{T^2}(|\nabla_A\psi|^2+\frac{1}{2}|\psi|^4) \\
 &   & +
(\int_{t_1-1}^{t_2+1}\int_{T^2}(|\nabla_A\psi|^2+\frac{1}{2}|\psi|^4))^{\frac{1}{2}})
\end{eqnarray*}
holds for any finite energy monopole $(A,\psi)$.
\elm

\noindent{\bf Proof:}
Without loss of generality, we assume $(A,\psi)$ is in the standard form. Then
\begin{eqnarray*}
|a(t_1)-a(t_2)|
& \leq & \sum_{i=1}^2|\int_{\gamma_{i}}a(t_1)-\int_{\gamma_{i}}a(t_2)|\\
&\leq & \sum_{i=1}^2(|\int_{\gamma_{i}}A_1(t_1)-\int_{\gamma_{i}}A_1(t_2)|\\
&  & +
\int_{\gamma_{i}}(|A_1(t_1)-a(t_1)|+|A_1(t_2)-a(t_2)|))\\
& \leq & \sum_{i=1}^2(\int_{t_1}^{t_2}\int_{\gamma_{i}}|dA|+
C(\|A_1(t_1)-a(t_1)\|_{L^2_1(T^2)}\\
&  &+\|A_1(t_2)-a(t_2)\|_{L^2_1(T^2)}))\\
& \leq & \sum_{i=1}^2 C(\int_{t_1}^{t_2}\int_{\gamma_{i}}|\psi|^2+(\int_{t_1\times T^2}|\psi|^4)^{\frac{1}{2}}+(\int_{t_2\times T^2}|\psi|^4)^{\frac{1}{2}}\\
& \leq & c_2(\int_{t_1}^{t_2}\int_{T^2}|\psi|^2+
\int_{t_1-1}^{t_2+1}\int_{T^2}(|\nabla_A\psi|^2+\frac{1}{2}|\psi|^4)\\
&  & +
(\int_{t_1-1}^{t_2+1}\int_{T^2}(|\nabla_A\psi|^2+\frac{1}{2}|\psi|^4))^{\frac{1}{2}})
\end{eqnarray*}
holds for any $t_1<t_2$. Here $\gamma_1$, $\gamma_2$ are the longitude and meridian.
\hfill ${\Box}$

\sde
Choose an increasing function ${\Gamma}(r)>0$ satisfying
$$
(\frac{c_1}{c(r)}+1){\Gamma}(r)+{\Gamma}(r)^{\frac{1}{2}}<c_2^{-1} r.
$$
A finite energy monopole $(A,\psi)$ is said to be
{\bf ``r-good''} if there are $t_0$ and $T$ with $T\leq t_0$ such that $a(t_0)\in\overline{{\hh}(2r)}$ and $\int_{T-1}^{\infty}\int_{T^2}(|\nabla_A \psi|^2+\frac{1}{2}|\psi|^4)<{\Gamma}(r)$.
\ede

The ``r-good'' monopoles of finite energy have the following good property.

\slm
Let $(A,\psi)$ be an ``r-good'' monopole of finite energy with $T$ as in the
Definition 1.3.9. Then for all $t\in [T,\infty)$, $a(t)$ is in ${\hh}(r)$.
Moreover, $a_\infty=\lim_{t\rightarrow\infty}a(t)$ exists in ${\hh}(r)$ 
and the following estimate holds for any $t\in [T,\infty)$:
\begin{eqnarray*}
|a(t)-a_\infty|
& \leq & c_2((\frac{c_1}{c(r)}+1)\int_{t-1}^{\infty}\int_{T^2}(|\nabla_A\psi|^2+
\frac{1}{2}|\psi|^4)\\
&      & +(\int_{t-1}^{\infty}\int_{T^2}(|\nabla_A\psi|^2+
\frac{1}{2}|\psi|^4))^{\frac{1}{2}}).
\end{eqnarray*}
\elm

\noindent{\bf Proof:}
It follows easily from the definition of ``r-goodness'' and Lemmas 2.3.6 and
2.3.8.
\hfill ${\Box}$

\slm
For any $r>0$, there exists a $\delta_0(r) >0$ with the following significance. For any $\epsilon >0$, there exists an $\epsilon_1 >0$ such that for
any ``r-good'' monopole $(A,\psi)$, when $\int_{t_0-1}^{\infty}\int_{T^2}(|\nabla_A\psi|^2+\frac{1}{2}|\psi|^4) <
\epsilon_1$ for $t_0\in [1,\infty)$, we have $\int_{2t_0}^{\infty}\int_{T^2}
|\psi|^2e^{2\delta t} < \epsilon$ for any $\delta\leq \delta_0(r)$.
\elm

\noindent{\bf Proof:} Without loss of generality, assume that $(A,\psi)$ is in the standard form.
By Lemma 2.3.10, there is a number $T>0$ such that for all $t\in [T,\infty)$, $a(t)$ is in ${\hh}(r)$ and $a_\infty=\lim_{t\rightarrow\infty}a(t)$ exists in ${\hh}(r)$. Moreover, the estimate
\begin{eqnarray*}
|a(t)-a_\infty|
& \leq & c_2((\frac{c_1}{c(r)}+1)\int_{t-1}^{\infty}\int_{T^2}(|\nabla_A\psi|^2+
\frac{1}{2}|\psi|^4)\\
&      & +(\int_{t-1}^{\infty}\int_{T^2}(|\nabla_A\psi|^2+
\frac{1}{2}|\psi|^4))^{\frac{1}{2}}).
\end{eqnarray*}
holds for $t\in [T,\infty)$. We can further apply a gauge transformation so that $a_\infty$ lies in
the compact set $\overline{{\hh}(r)}/({\bf Z}\oplus{\bf Z})$. There exists a $\delta_0(r) >0$ such that
for any $a\in\overline{{\hh}(r)}/({\bf Z}\oplus{\bf Z})$,
$$
\|D^{T^2}_a\psi\|_{L^2(T^2)}^2\geq 4\delta_0(r)\|\psi\|_{L^2_1(T^2)}^2
$$
for any $\psi\in{\Gamma}(W_0)$.
Set $u(t)=\int_{t\times T^2}|\psi|^2$. Then we have
$\frac{\partial^2}{\partial t^2}u(t)=2\int_{t\times T^2}(\langle
\frac{\partial^2}{\partial t^2}\psi,\psi\rangle+|\frac{\partial}{\partial t}
\psi|^2)$. But 
\begin{eqnarray*}
\int_{t\times T^2}\langle\frac{\partial^2}{\partial t^2}\psi,\psi\rangle
& = & \int_{t\times T^2}(|D^{T^2}_{A_1}\psi|^2-|A_0\psi|^2-
\langle\frac{\partial A_1}{\partial t}\psi+\frac{\partial A_0}{\partial t}\psi,
\psi\rangle\\
& \geq & \int_{t\times T^2}|D^{T^2}_{a_\infty}\psi|^2-\theta(t)
\|\psi\|_{L^2_1(T^2)}^2(t)
\end{eqnarray*}
where $\theta(t)\rightarrow 0$ as $\int_{t-1}^{\infty}\int_{T^2}(|\nabla_A\psi|^2+\frac{1}{2}|\psi|^4)
\rightarrow 0$
by the estimates in Lemma 2.3.4, Corollary 2.3.3, and Lemma 2.3.10. So when 
$
\int_{t_0-1}^{\infty}\int_{T^2}(|\nabla_A\psi|^2+\frac{1}{2}|\psi|^4)
$
is small enough, we have 
$$
\frac{\partial^2}{\partial t^2}u(t)\geq 4\delta_0(r) u(t) 
$$
for $t\geq t_0$. By the maximum principle, we have
$u(t)\leq e^{4\delta_0(r)(t_0-t)}u(t_0)$ for $t \geq t_0$. Hence 
$$
\int_{2t_0}^{\infty}u(t)e^{2\delta t}dt\leq C(\delta_0(r))u(t_0)
$$ 
holds for any $\delta\leq\delta_0(r)$.
The lemma follows easily.
\hfill ${\Box}$

\sprop
Assume that the ``r-good'' monopole of finite energy $(A,\psi)$ is in the standard form. Then $(A,\psi)$ is in ${\A}_\delta(r)$ for any weight $\delta\in 
(0,\min(\delta(r),\delta_0(r)))$. Moreover, the following estimate
$$
\|(A-a_\infty,\psi)\|_{L^2_{2,\delta}([T,\infty)\times T^2)}\leq c_3(\delta)
(\int_{T-1}^\infty\int_{T^2}|\psi|^2e^{2\delta t})^{\frac{1}{2}}
$$
holds for a constant $c_3(\delta)$ and any $T\in [1,\infty)$.
Here $\delta(r)$ is referred to Proposition 2.1.12.
\eprop

\noindent{\bf Proof:}
It follows from Lemma 2.3.11, the estimates in Lemma 2.3.4, Taubes inequality and standard elliptic estimates.
\hfill ${\Box}$

\noindent{\bf The convergence of ``r-good'' monopoles of finite energy}

\sprop
Let $(A_n,\psi_n)$ be a sequence of ``r-good'' monopoles of finite energy.
Then a subsequence of $(A_n,\psi_n)$ converges in ${\M}_\delta(\frac{r}{2})$ to a 
``$\frac{1}{2}r$-good'' monopole of finite energy $(A_0,\psi_0)$ for any 
weight $\delta\in (0,\min(\delta(\frac{r}{2}),\delta_0(\frac{r}{2})))$.
\eprop

\noindent{\bf Proof:}
The Weitzenb\"ock formula and maximum principle yield an upper bound
$K$ for the $C^0$ norm of the spinors (see Corollary 2.3.3).
Then the existence of a local Hodge gauge (\cite{KM1}) plus elliptic regularity and a
patching argument (\cite{U}) imply the existence of a ``weak'' limit. More 
precisely, there exist a finite energy monopole $(A_0,\psi_0)$, a subsequence of $(A_n,\psi_n)$ (still denoted by $(A_n,\psi_n)$) and a
sequence of gauge transformations $s_n$ such that $s_n\cdot (A_n,\psi_n)$ converges to $(A_0,\psi_0)$ in $C^\infty$ over any compact subset of $Y$. We can further assume that $s_n\cdot (A_n,\psi_n)$ are in the standard form and therefore the limit $(A_0,\psi_0)$ is also in the standard form. For simplicity we still use $(A_n,\psi_n)$ to denote $s_n\cdot (A_n,\psi_n)$.

Take $T_0$ large enough so that $\int_{T_0-1}^{\infty}\int_{T^2}(|\nabla_{A_0}\psi_0|^2+
\frac{1}{2}|\psi_0|^4)<{\Gamma}
(\frac{r}{2})$ (see Definition 2.3.9). Note that for any finite energy monopole $(A,\psi)$ the Weitzenb\"ock formula yields the following equation
$$
\int_{0}^{\infty}\int_{T^2}(|\nabla_A\psi|^2+\frac{1}{2}|\psi|^4)=
\int_{0\times T^2}\langle D^{T^2}_{A_1}\psi,\psi\rangle.
$$
It follows from the ``weak'' convergence of $(A_n,\psi_n)$ to $(A_0,\psi_0)$ that there is an $N$ such that when $n>N$ we have
$$
\int_{T_0-1}^{\infty}\int_{T^2}(|\nabla_{A_n}\psi_n|^2+
\frac{1}{2}|\psi_n|^4)<{\Gamma}
(\frac{r}{2})<{\Gamma}(r).
$$
Since $(A_n,\psi_n)$ are ``r-good'', by Lemma 2.3.10, $a_n(t)$ is in ${\hh}(r)$
for any $n>N$ and $t\in [T_0,\infty)$. From this it follows that $(A_0,\psi_0)$
is a ``$\frac{1}{2}r$-good'' monopole of finite energy, and therefore is in
${\A}_\delta(\frac{r}{2})$ for any weight $\delta\in (0,\min(\delta(\frac{r}{2}),\delta_0(\frac{r}{2}))$. It is also easy to see
that $a_{n,\infty}\rightarrow a_{0,\infty}$.

In order to prove that $(A_{n},\psi_{n})$ converges to $(A_0,\psi_0)$ in 
${\A}_\delta(\frac{r}{2})$ for any given weight $\delta\in (0,\min(\delta(\frac{r}{2}),\delta_0(\frac{r}{2}))$, it suffices to prove
that given any $\epsilon >0$, there exist $t_0\in [1,\infty)$ and $N$
such that when $n > N$,
$$
\|(A_n-a_{n,\infty},\psi)\|_{L^2_{2,\delta}([2t_0+1,\infty)\times T^2)}
<\epsilon.
$$
This goes as follows.
By Lemma 2.3.11, there exists an $\epsilon_1>0$ such that when
$$
\int_{t_0-1}^{\infty}\int_{T^2}(|\nabla_{A_n}\psi_n|^2+
\frac{1}{2}|\psi_n|^4) < \epsilon_1,
$$
we have 
$$\int_{2t_0}^{\infty}\int_{T^2}|\psi_n|^2 e^{2\delta t}< (c_3^{-1}(\delta)\epsilon)^2.
$$
Now take $t_0$ large enough so that
$$
\int_{t_0-1}^{\infty}\int_{T^2}(|\nabla_{A_0}\psi_0|^2+
\frac{1}{2}|\psi_0|^4) < \frac{\epsilon_1}{2}.
$$
Then there exists an $N$ such that when $n>N$ we have 
$$
\int_{t_0-1}^{\infty}\int_{T^2}(|\nabla_{A_n}\psi_n|^2+
\frac{1}{2}|\psi_n|^4) < \epsilon_1.
$$
Therefore we have 
$$
\|(A_n-a_{n,\infty},\psi)\|_{L^2_{2,\delta}([2t_0+1,\infty)\times T^2)}
<\epsilon
$$
by Proposition 2.3.12. Hence the proposition is proved.
\hfill ${\Box}$

\chapter{The Gluing Formula}

\section{The gluing of moduli spaces}

Assume that $Y_i$ $(i=1,2)$ are oriented cylindrical end 3-manifolds over $T_i^2$ where $Y_2$ is actually diffeomorphic 
to $D^2\times S^1$ carrying a metric whose scalar curvature is non-negative and somewhere positive. By the Weitzenb\"{o}ck formula, the moduli space 
${\M}(Y_2)$ actually only consists of reducible solutions.
Assume that there exists an orientation reversing isometry $h: T^2_1
\rightarrow T^2_2$ which is covered by the corresponding bundle maps. For 
any $L > 0$, we can form an oriented Riemannian 3-manifold $Y_L$ as follows:
$$
Y_L=Y_1\setminus [L+1,\infty)\times T^2_1\bigcup_h Y_2\setminus [L+1,\infty)
\times T^2_2,
$$
where $h:(L, L+1)\times T^2_1\rightarrow (L+1,L)\times T^2_2$ is given by
$h(L+t,x)=(L+1-t, h(x))$ for $t\in (0,1)$ and $x\in T^2_1$. Note that
the isometry $h:T^2_1\rightarrow T^2_2$ induces an isometry $h: {\hh}^1(T^2_1)\rightarrow {\hh}^1(T^2_2)$.

Throughout this section, we fix a small $r>0$ and a small
weight $\delta$. For simplicity, we omit the dependence of $r$ and $\delta$ in the discussion. We also omit the perturbation data since it is vanishing on the neck.

Let $\tilde{{\G}}(Y_i)$ be the normal subgroup of ${\G}(Y_i)$ which consists of 
elements in the component of identity. We define $\tilde{{\M}}(Y_i)$ to be the
space of $\tilde{{\G}}(Y_i)$-equivalence classes of the solutions to the Seiberg-Witten equations on $Y_i$ ($i=1,2$). Then $\tilde{{\M}}(Y_i)$ is a ${\bf Z}$-fold cover of ${\M}(Y_i)$:
$$
\begin{array}{ccccccccc}
 0 & \rightarrow & {\bf Z} & \rightarrow & \tilde{{\M}}(Y_i) & \rightarrow & {\M}(Y_i)
& \rightarrow & 0. 
\end{array}
$$
The irreducible part of $\tilde{{\M}}(Y_1)$ is denoted by $\tilde{{\M}}^\ast(Y_1)$, 
which is a ${\bf Z}$-fold cover of ${\M}^\ast(Y_1)$.

Let ${{\S}(Y_1,Y_2)}$ be the set of pairs $(\alpha_1,\alpha_2)\in\tilde{{\M}}^\ast(Y_1)\times\tilde{{\M}}(Y_2)$ such that
there are smooth representatives $(A_1,\psi_1)$ and $(A_2,\psi_2)$ satisfying
$hR_1(A_1,\psi_1)=R_2(A_2,\psi_2)$,
where $R_i:{\A}(Y_i)\rightarrow{\hh}^1(T^2_i)\otimes i{\bf R}$. 

\sde
$(\alpha_1,\alpha_2)\in {{\S}(Y_1,Y_2)}$ is said to be regular if 
\begin{enumerate}
\item the map $[R_1]: \tilde{{\M}}^\ast(Y_1)\rightarrow {\hh}^1(T^2_1)\otimes i{\bf R}$ 
is injective at $\alpha_1$.
\item $h[R_1](\tilde{{\M}}^\ast(Y_1))$ and $[R_2](\tilde{{\M}}(Y_2))$ intersect transversally at $[R_2](\alpha_2)$.
\end{enumerate}
\ede
Note that for a generic perturbation, ${{\S}(Y_1,Y_2)}$ consists of regular pairs (Proposition 2.2.2). We assume that ${{\S}(Y_1,Y_2)}$ is regular throughout this chapter. 

For any $L >0$, fix a cut-off function $\rho_L(t)$ which equals to one for
$t<L$ and equals to zero for $t>L+1$ with $|\rho^\prime_L|<2$. 
Given a regular pair $(\alpha_1,\alpha_2)$ with smooth representatives
$(A_1,\psi_1)$ and $(A_2,\psi_2)$,  we construct an
``almost'' monopole $(A_L,\psi_L)$ on $Y_L$ as follows:
$$
\begin{array}{cccc}
\psi_L & = & \rho_L\psi_1+(1-\rho_L)h^{-1}\psi_2  & \mbox{on} \hspace{2mm} [L,L+1]\times T^2_1\\
A_L & = & \rho_L(A_1-R_1(A_1))+R_1(A_1) & \\
    &   & +(1-\rho_L)h^{-1}(A_2-R_2(A_2))
& \mbox{on} \hspace{2mm} [L,L+1]\times T^2_1\\
\psi_L & = & \psi_i & \mbox{on} \hspace{2mm} Y_i\setminus [L,\infty)\times T^2_i\\
A_L & = & A_i &\mbox{on} \hspace{2mm} Y_i\setminus [L,\infty)\times T^2_i.
\end{array}
$$
Note that $\psi_L$ is compactly supported in $Y_1\setminus [L+1,\infty)\times
T^2$.

\sprop
Assume that $(\alpha_1,\alpha_2)$ is regular.
Then for large $L$, the ``almost'' monopole $(A_L,\psi_L)$ can be deformed to 
a non-degenerate monopole $T(A_L,\psi_L)$ satisfying 
$$
\|T(A_L,\psi_L)-(A_L,\psi_L)\|_{L^2_1(Y_L)}\leq CL^2e^{-\delta L}.
$$
Moreover, any monopole $(A,\psi)$ on $Y_L$ which is in the $L^2_1$-ball of
radius $K_1L^{-6}$ centered at $(A_L,\psi_L)$ is gauge equivalent to $T(A_L,\psi_L)$. In particular, there is a well-defined gluing map $T: {{\S}(Y_1,Y_2)}\rightarrow {\M}^\ast(Y_L)$ by
$T(\alpha_1,\alpha_2)=[T(A_L,\psi_L)]$.
\eprop

The following estimate on $(A_L,\psi_L)$ is straightforward.

\slm
$
\|(\ast dA_L +\tau(\psi_L,\psi_L), D_{A_L}\psi_L)\|_{L^2(Y_L)}\leq
Ce^{-\delta L}.
$
\elm

Next we estimate the lowest eigenvalue of $\Delta_L=d^\ast d+|\psi_L|^2$. Set 
$$
\lambda_L=\inf_{f\neq 0}\frac{\int_{Y_L}|\Delta_L f|^2}{\int_{Y_L}|f|^2}.
$$

\slm
Assume that one of $\psi_1$ or $\psi_2$ is not identically zero. For any function $\gamma(L)=o(L^{-4})$, there exists an $L_0>0$ such that when $L>L_0$,
we have $\lambda_L > \gamma(L)$.
\elm

The basic idea of the proof is the same as in Theorem 4 of Appendix B, but the argument is more
difficult. We postpone the proof to the end of this section. From now on,
we assume that one of $\psi_1$ or $\psi_2$ is not identically zero.

\scor
The norm of $\Delta_L^{-1}: L^2(Y_L)\rightarrow L^2_2(Y_L)$ is at most $L^3$ for large $L$.
\ecor

\noindent{\bf Proof:}
In Lemma 3.1.4, take $\gamma(L)=K^2L^{-6}$ with $K$ to be determined later.
There exists a constant $C$ independent of $L$ such that for any $f\in L^2_2(Y_L)$, we have
\begin{eqnarray*}
\|f\|_{L^2_2(Y_L)}
& \leq & C(\|\Delta_Lf\|_{L^2(Y_L)}+\|f\|_{L^2(Y_L)})\\
& \leq & C(\|\Delta_Lf\|_{L^2(Y_L)}+K^{-1}L^3\|\Delta_Lf\|_{L^2(Y_L)})\\
& \leq & L^3\|\Delta_Lf\|_{L^2(Y_L)}
\end{eqnarray*}
for large $L$ and a suitable choice of $K$. This proves the lemma.
\hfill ${\Box}$

Let $T{\B}^\ast_{(A_L,\psi_L)}$ be the tangent space of ${\B}^\ast(Y_L)$ at $(A_L,\psi_L)$, ${\E}_{(A_L,\psi_L)}$ be the $L^2$ completion of 
$T{\B}^\ast_{(A_L,\psi_L)}$. Then
$T{\B}^\ast_{(A_L,\psi_L)}=\{(a,\phi)\in{\A}(Y_L)|-d^\ast a+i\langle i\psi_L,\phi\rangle_{Re}=0\}$.

\slm
There exist constants $K_1$, $K_2$ with the following significance. When $L$ is
sufficiently large, for any $(A,\psi)\in{\A}(Y_L)$ satisfying
$$
\|(A,\psi)-(A_L,\psi_L)\|_{L^2_1(Y_L)}\leq K_1L^{-6},
$$
there exists a gauge transformation $s\in{\G}(Y_L)$ such that 
$s\cdot (A,\psi)-(A_L,\psi_L)\in T{\B}^\ast_{(A_L,\psi_L)}$ and
$$
\|s\cdot (A,\psi)-(A_L,\psi_L)\|_{L^2_1(Y_L)}\leq K_2L^3\|(A,\psi)-(A_L,\psi_L)\|_{L^2_1(Y_L)}.
$$
\elm

\noindent{\bf Proof:}
The point of this lemma is to have an estimate on the size of the local slice
at $(A_L,\psi_L)$ if an upper bound for the norm of $\Delta_L^{-1}$ is known
(Corollary 3.1.5).

Assume that $(A,\psi)$ is in ${\A}(Y_L)$. Set $(a,\phi)=(A,\psi)-(A_L,\psi_L)$ for
simplicity. To construct the local slice, we look for $f\in L^2_2$ such that
$$
-d^\ast(A-A_L-df)+i\langle i\psi_L,e^f\psi-\psi_L\rangle_{Re}=0.
$$
This can be written in terms of $(a,\phi)$ as
$$
(\Delta_L+\langle i\psi_L,i\phi\rangle_{Re})f+G(\phi,f)=d^\ast a-i\langle i\psi_L,\phi\rangle_{Re},
$$
where $G(\phi,f)=i\langle i\psi_L, (e^f-f-1)(\phi+\psi_L)\rangle_{Re}$
satisfying
$$
\|G(\phi,f_1)-G(\phi,f_2)\|_{L^2}\leq C\max(\|f_1\|_{L^2_2}, \|f_2\|_{L^2_2})
\|f_1-f_2\|_{L^2_2}
$$
for some constant $C$ and $\phi,f_i$ satisfying $\|\phi\|_{L^2_1}<1$ and 
$\|f_i\|_{L^2_2}<1$ for $i=1,2$.
The lemma follows by applying Banach lemma to the map
$$
B(f)=(\Delta_L+\langle i\psi_L,i\phi\rangle_{Re})^{-1}(d^\ast a-i\langle i\psi_L,\phi\rangle_{Re}-G(\phi,f))
$$
mapping an $L^2_2$-ball of radius $KL^{-3}$ into itself for some constant $K$.
\hfill ${\Box}$

Next we deform the ``almost'' monopole $(A_L,\psi_L)$ to a monopole.
Let $\Pi$ be the $L^2$ orthogonal projection onto ${\E}_{(A_L,\psi_L)}$. 
For any $(a,\phi)\in T{\B}^\ast_{(A_L,\psi_L)}$, we define
\begin{eqnarray*}
L(a,\phi)
& = & \Pi(\ast d(A_L+a)+\tau(\psi_L+\phi), D_{(A_L+a)}(\psi_L+\phi))\\
& = & (\ast dA_L+\tau(\psi_L), D_{A_L}\psi_L) + \nabla s_{(A_L,\psi_L)}
(a,\phi) + \Pi Q(a,\phi).
\end{eqnarray*}
Here $\nabla s_{(A_L,\psi_L)}: T{\B}^\ast_{(A_L,\psi_L)}\rightarrow
{\E}_{(A_L,\psi_L)}$ is given by 
$$
\nabla s_{(A_L,\psi_L)}(a,\phi)=
(\ast da+2\tau(\psi_L,\phi)-df(\phi), D_{A_L}\phi+a\psi_L+f(\phi)\psi_L)
$$
with $f(\phi)$ given by the equation
$$
\Delta_Lf=i\langle D_{A_L}\psi_L,i\phi\rangle_{Re}.
$$
$Q(a,\phi)=(\tau(\phi), a\phi)$ satisfies
$$
\|Q(a_1,\phi_1)-Q(a_2,\phi_2)\|_{L^2(Y_L)}\leq 
C(\|(a_1,\phi_1)\|_{L^2_1}+\|(a_2,\phi_2)\|_{L^2_1})
\|(a_1,\phi_1)-(a_2,\phi_2)\|_{L^2_1}.
$$

\slm
There exists a constant $K_3$ such that when $\|(a,\phi)\|_{L^2_1(Y_L)}\leq K_3L^{-3}$ for large enough $L$, $L(a,\phi)=0$ implies that 
$$
(\ast d(A_L+a)+\tau(\psi_L+\phi), D_{(A_L+a)}(\psi_L+\phi))=0.
$$
\elm

\noindent{\bf Proof:}
\begin{eqnarray*}
L(a,\phi)
& = & \Pi(\ast d(A_L+a)+\tau(\psi_L+\phi), D_{(A_L+a)}(\psi_L+\phi))\\
& = & (\ast d(A_L+a)+\tau(\psi_L+\phi)-dg(a,\phi), D_{(A_L+a)}(\psi_L+\phi)
+g(a,\phi)\psi_L)
\end{eqnarray*}
where $g(a,\phi)$ satisfies the equation
$$
\Delta_Lg=i\langle D_{(A_L+a)}(\psi_L+\phi),i\phi\rangle_{Re}.
$$
If $L(a,\phi)=0$, then $D_{(A_L+a)}(\psi_L+\phi)+g(a,\phi)\psi_L=0$.
So for large $L$, we have
\begin{eqnarray*}
\|g(a,\phi)\|_{L^2_2(Y_L)}
& \leq & L^3\|\Delta_Lg(a,\phi)\|_{L^2(Y_L)}\\
& \leq & L^3\|\langle g(a,\phi)\psi_L,i\phi\rangle_{Re}\|_{L^2(Y_L)}\\
& \leq & CL^3\|\psi_L\|_{C^0}\|\phi\|_{L^2_1(Y_L)}\|g(a,\phi)\|_{L^2_2(Y_L)}.
\end{eqnarray*}
If $\|(a,\phi)\|_{L^2_1(Y_L)}\leq K_3L^{-3}$ for a small enough constant $K_3$, we have $g(a,\phi)=0$, which proves the lemma.
\hfill ${\Box}$

\slm
Assume that $(\alpha_1,\alpha_2)$ is regular.
Then $\nabla s_{(A_L,\psi_L)}: T{\B}^\ast_{(A_L,\psi_L)}
\rightarrow {\E}_{(A_L,\psi_L)}$ is invertible for large $L$. Moreover, the norm
of $(\nabla s_{(A_L,\psi_L)})^{-1}$ is at most $L^2$ for large $L$.
\elm

\noindent{\bf Proof:}
First of all, $\nabla s_{(A_L,\psi_L)}: T{\B}^\ast_{(A_L,\psi_L)}
\rightarrow {\E}_{(A_L,\psi_L)}$ can be extended to an operator
$$
{\K}^\prime_{(A_L,\psi_L)}=\left(\begin{array}{ccc}
\nabla s_{(A_L,\psi_L)} & 0 & 0\\
0 & 0 & d_{(A_L,\psi_L)}\\
0 & d^\ast_{(A_L,\psi_L)}
\end{array} \right),
$$
where $d_{(A_L,\psi_L)}(f)=(-df, f\psi_L)$.
Secondly, by Theorem 4 in Appendix B, the lowest eigenvalue of
$$
{\K}_{(A_L,\psi_L)}=\left(\begin{array}{ccc}
D_{A_L} & \psi_L\cdot & \psi_L\cdot\\
2\tau(\psi_L,\cdot) & \ast d & -d\\
i\langle i\psi_L,\cdot\rangle_{Re} & -d^\ast & 0
\end{array} \right)
$$
is at least $KL^{-2}$ for any constant $K$ when $L$ is sufficiently large.
Here we essentially use the fact that $\psi_L$ is 
identically zero on the $Y_2$ side and not identically zero on the $Y_1$ side
and the assumption that $(\alpha_1,\alpha_2)$ is regular
so that the regularity and the
transversality conditions in Theorem 4 of Appendix B hold.
Finally, the difference between ${\K}^\prime_{(A_L,\psi_L)}$ and ${\K}_{(A_L,\psi_L)}$ can be ignored since the norm of 
${\K}^\prime_{(A_L,\psi_L)}-{\K}_{(A_L,\psi_L)}$
is bounded from above by $CL^6\|D_{A_L}\psi_L\|_{C^0}\leq cL^6
e^{-\delta L}$ by Lemma 3.1.3 and Corollary 3.1.5.
So the lowest eigenvalue of $\nabla s_{(A_L,\psi_L)}$
is bounded from below by $KL^{-2}$ for any constant $K$ 
when $L$ is large enough. The lemma follows easily.
\hfill ${\Box}$

\noindent{\bf The Proof of Proposition 3.1.2:}

In order to deform the ``almost'' monopole $(A_L,\psi_L)$ to a monopole, we
need to solve the equation $L(a,\phi)=0$ for $(a,\phi)\in T{\B}^\ast_{(A_L,\psi_L)}$. This equation can be
written as
$$
(a,\phi)=-(\nabla s_{(A_L,\psi_L)})^{-1}((\ast d{A_L}+\tau(\psi_L), D_{A_L}\psi_L)+\Pi Q(a,\phi)).
$$
Assuming that $(\alpha_1,\alpha_2)$ is regular, it then follows from Lemmas
3.1.3 and 3.1.8 that the map
$$
B(a,\phi)=-(\nabla s_{(A_L,\psi_L)})^{-1}((\ast d{A_L}+\tau(\psi_L), D_{A_L}\psi_L)+\Pi Q(a,\phi))
$$ 
maps an $L^2_1$-ball of radius $KL^{-2}$ in $T{\B}^\ast_{(A_L,\psi_L)}$ into itself and satisfies
$$
\|B(a_1,\phi_1)-B(a_2,\phi_2)\|_{L^2_1(Y_L)} <
\|(a_1,\phi_1)-(a_2,\phi_2)\|_{L^2_1(Y_L)}
$$
for large enough $L$ and some small enough constant $K$. 
Therefore the equation $L(a,\phi)=0$ has a unique 
solution $(a_L,\phi_L)$ in the $L^2_1$-ball of radius $KL^{-2}$. By Lemma 3.1.7, the resulting monopole is $T(A_L,\psi_L)=(A_L,\psi_L)+(a_L,\phi_L)$, and 
$$
\|T(A_L,\psi_L)-(A_L,\psi_L)\|_{L^2_1(Y_L)}\leq CL^2e^{-\delta L}.
$$
Suppose that monopole $(A,\psi)$ is in an $L^2_1$-ball of radius $K_1L^{-6}$ centered at $(A_L,\psi_L)$. By Lemma 3.1.6, $(A,\psi)$ is gauge equivalent to 
a monopole in the local slice with distance from $(A_L,\psi_L)$ less than $K_2K_1L^{-3}$, and it must be $T(A_L,\psi_L)$ by the uniqueness of
the solution $(a_L,\phi_L)$ to the equation $L(a,\phi)=0$.
In particular, $[T(A_L,\psi_L)]$ depends only on the homotopy class of $(A_L,\psi_L)$, which implies that there is a well-defined gluing map $T:{{\S}(Y_1,Y_2)}\rightarrow{\M}^\ast(Y_L)$. By Lemma 3.1.8 and the estimate 
$$
\|T(A_L,\psi_L)-(A_L,\psi_L)\|_{L^2_1(Y_L)}\leq CL^2e^{-\delta L},
$$
$[T(A_L,\psi_L)]$ is non-degenerate. Therefore the proposition is proved.
\hfill ${\Box}$

\noindent{\bf The Proof of Lemma 3.1.4:}

Suppose that there is a sequence of $L_n \rightarrow \infty$ such that
$$
\lambda_{L_n}\leq \gamma(L_n)=o(L_n^{-4}),
$$
then there exists a sequence of $c_n>0$, $f_n\neq 0$ such that 
$\Delta_{L_n}f_n=c_n^2f_n$ and $c_n=o(L_n^{-1})$.

\vspace{2mm}

\noindent{\bf Claim:}{\em \hspace{2mm}
There exist constants $M$ and $L_0$ with the following significance. The 
$f_n^\prime s$ can be chosen such that $\|f_n\|_{C^0(Y_{L_n})}\leq M$, and
one of the following is true:
\begin{enumerate}
\item either $\int_{Y_1(L_0)}|f_n|^2$ or $\int_{Y_2(L_0)}|f_n|^2$ is equal to one;
\item either $\|f_n\|_{L^2(T^2_1)}(L_0)$ or $\|f_n\|_{L^2(T^2_2)}(L_0)$ is greater than or equal to one.
\end{enumerate}
Here $Y_i(L_0)=Y_i\setminus (L_0,\infty)\times T^2_i$, $i=1,2$}.

\vspace{2mm}

Assuming the {\bf Claim}, Lemma 3.1.4 is proved as follows. By elliptic 
estimates, we can select a subsequence of $f_n$ which is convergent in 
$C^\infty$ to $f_i$ on $Y_i$ ($i=1,2$) on any compact subset. Moreover, $f_1,f_2$
satisfy the following conditions:
\begin{itemize}
\item [{a)}] $d^\ast df_i+|\psi_i|^2f_i=0$ on $Y_i$, $i=1,2$;
\item [{b)}] $\|f_i\|_{C^0(Y_i)}\leq M$, $i=1,2$;
\item [{c)}] one of $f_1$ or $f_2$ is not identically zero.
\end{itemize}
Lemma 3.1.4 is proved if we show that conditions $a)$ and $b)$ contradict 
condition $c)$. In fact, by condition $a)$, for any $t$, we have
\begin{eqnarray*}
0
& = & \int_{Y_i(t)}\langle d^\ast df_i+|\psi_i|^2f_i, f_i\rangle\\
& = & \int_{Y_i(t)}(|df_i|^2+|\psi_i|^2|f_i|^2)-\frac{1}{2}\frac{\partial}
{\partial t}(\int_{T^2_i}|f_i|^2)(t).
\end{eqnarray*}
Since $\xi_i(t)=\int_{T^2_i}|f_i|^2(t)$ is bounded by $b)$, there exists a
sequence of $t_n$ such that $\frac{\partial\xi_i}
{\partial t}(t_n)\rightarrow 0$ as $t_n\rightarrow \infty$. Therefore
$\int_{Y_i}(|df_i|^2+|\psi_i|^2|f_i|^2)=0$, $i=1,2$.
So $f_1$ and $f_2$ are constant functions and one of them is zero, since one
of $\psi_1$ or $\psi_2$ is not identically zero. But in the proof of the
{\bf Claim}, it is easy to see that $f_1=f_2$. So both of $f_1$ and $f_2$ are
zero, contradicting $c)$.

\vspace{2mm}

\noindent{\bf The Proof of the Claim:}

For simplicity, we omit the subscript $L_n$ or $n$ if no confusion is caused.
Write $f=g_1+g_2$ where $g_1\in {{\K}}er\;\dd^\ast\dd$ and $g_2\in
({{\K}}er\;\dd^\ast\dd)^\perp$.

Pick $L_0>0$ large enough, there are two possibilities:
\begin{itemize}
\item On $[L_0,2L+1-L_0]$, $\max|g_1|\leq\max\|g_2\|_{L^2(T^2_1)}$.
In this case, by the maximum principle, for large $L_0$, $\|g_2\|_{L^2(T^2_1)}$
can not reach its maximum in the interior of $[L_0,2L+1-L_0]$. The {\bf Claim} 
follows easily in this case.
\item On $[L_0,2L+1-L_0]$, $\max\|g_2\|_{L^2(T^2_1)}\leq
\max|g_1|$. In this case, we need to show that for large enough $L$ 
either $|g_1|$ reaches its maximum at the endpoints of $[L_0,2L+1-L_0]$, from which the {\bf Claim} follows, or either
$\max|g_1|\leq K|g_1(L_0)|$ or $\max|g_1|\leq K|g_1(2L+1-L_0)|$ holds for a
constant $K$ independent of $L$.
\end{itemize}

Assume that we are in the second case. On $[L_0,2L+1-L_0]$, we have
$$
-\frac{\partial^2}{\partial t^2}g_1(t)+h(t)=c^2g_1(t)
$$
where $h(t)$ is the $L^2$-projection of $|\psi_L|^2f$ into ${{\K}}er\;\dd^\ast\dd$.
$c=o(L^{-1})$ and $h(t)$ satisfies
$$
|h(t)|\leq Ke^{-2\delta t}\max|g_1| \hspace{2mm} on \hspace{2mm} [L_0, L+1) 
$$
and 
$$
|h(t)|\leq Ke^{-2\delta(2L+1-t)}\max|g_1| \hspace{2mm} on \hspace{2mm} 
(L, 2L+1-L_0].
$$
Set $g_3(t)=c^{-1}\frac{\partial}{\partial t}g_1(t)$, then we have
$$
\begin{array}{c}
\frac{\partial}{\partial t}g_1(t)=cg_3(t)\\
\frac{\partial}{\partial t}g_3(t)=-cg_1(t)+c^{-1}h(t).
\end{array}
$$
These equations can be written equivalently as
$$
\frac{\partial}{\partial t}\left[e^{Ct}\left(\begin{array}{c}
g_1(t)\\g_3(t)
\end{array} \right)\right]
=e^{Ct}\left(\begin{array}{c}
0\\c^{-1}h(t)
\end{array} \right)
$$
where $C=\left(\begin{array}{cc}
0 & -c\\
c & 0
\end{array} \right)$. 
Note that $e^{Ct}=\left(\begin{array}{cc}
\cos ct & -\sin ct\\
\sin ct & \cos ct
\end{array} \right)$.

Since $c=o(L^{-1})$, we have $\cos ct>\frac{1}{2}$ and
$$|\int_{L_0}^t\sin cs\cdot c^{-1}h(s)ds|\leq K\int_{L_0}^t|sh(s)|ds
\leq Ke^{-\delta L_0}\max|g_1|
$$
for large $L$ and any $t\in [L_0,L+1)$. On the other hand,
$$
g_1(t)\cos ct -g_3(t)\sin ct=-\int_{L_0}^t\sin cs\cdot c^{-1}h(s)ds
+ g_1(L_0)\cos cL_0 -g_3(L_0)\sin cL_0.
$$
So if $|g_1(t)|$ reaches its maximum in the interior of $[L_0,2L+1-L_0]$, 
without loss of generality, assuming that it is in $(L_0, L+1)$, then
$$
\frac{1}{2}\max|g_1|\leq Ke^{-\delta L_0}\max|g_1|+
|g_1(L_0)\cos cL_0 -g_3(L_0)\sin cL_0|.
$$
Hence for large $L_0$, we have 
$$
\max|g_1|\leq 4|g_1(L_0)\cos cL_0 -g_3(L_0)\sin cL_0|.
$$
On the other hand,
$$
g_3(t)=\int_{L_0}^t\cos c(s-t)\cdot c^{-1}h(s)ds +
g_3(L_0)\cos c(L_0-t) +g_1(L_0)\sin c(L_0-t).
$$
Assume that $|g_1(t)|$ has its maximum at $t_0\in (L_0, L+1)$. Then $g_3(t_0)=0$
and 
\begin{eqnarray*}
&      & |g_1(L_0)\sin c(L_0-t_0)+g_3(L_0)\cos c(L_0-t_0)|\\
& \leq & \int_{L_0}^{t_0}|c^{-1}h(s)|ds\\
& \leq & Kc^{-1}e^{-2\delta L_0}\max|g_1|\\
& \leq & 4Kc^{-1}e^{-2\delta L_0}|g_1(L_0)\cos cL_0 -g_3(L_0)\sin cL_0|.
\end{eqnarray*}
Since $c=o(L^{-1})$ as $L\rightarrow \infty$ and $L_0e^{-2\delta L_0}=o(1)$ as
$L_0\rightarrow \infty$, we have
$$
|g_3(L_0)|\leq 10Kc^{-1}e^{-2\delta L_0}|g_1(L_0)|+|g_1(L_0)|.
$$
Hence
\begin{eqnarray*}
\max|g_1|
& \leq & 4|g_1(L_0)\cos cL_0 - g_3(L_0)\sin cL_0|\\
& \leq & 4|g_1(L_0)|+(10Kc^{-1}e^{-2\delta L_0}|g_1(L_0)|+|g_1(L_0)|)
|\sin cL_0|\\
& \leq & 5|g_1(L_0)|+10KL_0e^{-2\delta L_0}|g_1(L_0)|\\
& \leq & K_1|g_1(L_0)|
\end{eqnarray*}
with $K_1$ independent of $L$ when $L_0$ is sufficiently large. This proves the {\bf Claim}.
(Observe that $|(g_1(t_0)\cos ct_0-g_3(t_0)\sin ct_0)-(g_1(L_0)\cos cL_0-g_3(L_0)\sin cL_0)|\leq Ke^{-\delta L_0}$ and $|(g_1(t_0)\cos ct_0-g_3(t_0)\sin ct_0)-(g_1(2L+1-L_0)\cos cL_0-g_3(2L+1-L_0)\sin cL_0)|\leq Ke^{-\delta L_0}$ where $t_0=L+\frac{1}{2}$, from which one sees $f_1=f_2$).

\section{Geometric limits}

Let $Y$ be an oriented integral homology 
3-sphere decomposed into $Y=Y_1\bigcup_{T^2}Y_2$, where $Y_1$ is the complement
of a tubular neighborhood of a knot and $Y_2$ is diffeomorphic to $D^2\times S^1$. Equip $Y$ with a Riemannian metric such that a collar neighborhood of $T^2$ is orientedly isometric to $(-1, 1)\times {\bf R}/2\pi{\bf Z}\times{\bf R}/2\pi{\bf Z}$ with $(-1,0)\times{\bf R}/2\pi{\bf Z}\times{\bf R}/2\pi{\bf Z}$ in $Y_1$ and $Y_2$ carries a non-negative, somewhere positive scalar curvature metric. We insert cylinders of lengths $2L+1$ and obtain a family of stretched versions $Y_L$ of $Y$. We also use $Y_1$ and $Y_2$ to denote the corresponding cylindrical end manifolds. Note that the finite energy monopoles on $Y_2$ are reducible by Weitzenb\"ock formula, so the moduli
space $\tilde{{\M}}(Y_2)$ is identified with the line ${\hh}^1(Y_2)\otimes i{\bf R}$ of imaginary valued ``bounded'' harmonic 1-forms on $Y_2$ which is embedded into $\h$ by the map $R_2$ (Lemma 2.1.4). 
With a small perturbation, we assume that  $\tilde{{\M}}(Y_2)$  misses the
lattice of ``bad'' points for the spin structure on $T^2$ where the twisted
Dirac operators are not invertible.

Let $(A_n,\psi_n)$ be a sequence of monopoles on $Y_{L_n}$, $L_n\rightarrow
\infty$. Weitzenb\"ock formula and the maximum principle yield an upper bound
$K$ for the $C^0$ norm of the spinors, which depends only on the scalar curvature of the manifolds. Then the existence of a local Hodge gauge (\cite{KM1}) plus elliptic regularity and a patching argument (\cite{U}) imply the existence of geometric limits as we stretch the neck.

\slm
There exists an $r>0$ with the following significance.
Let $(A_n,\psi_n)$ be a sequence of monopoles on $Y_{L_n}$, $L_n\rightarrow
\infty$. Then there exist a sequence of gauge transformations $s_n$ and a
pair of finite energy monopoles $(A_i,\psi_i)$ on $Y_i$ ($i=1,2$) such that
a subsequence of $s_n\cdot (A_n,\psi_n)$ converges in $C^\infty$ to $(A_i,\psi_i)$ on any compact subset of $Y_i$ ($i=1,2$). Moreover, 
$(A_1,\psi_1)$ and $(A_2,\psi_2)$ are ``r-good'' monopoles and
have the same limiting value, i.e. $R_1(A_1,\psi_1)=R_2(A_2,\psi_2)$.
\elm

\noindent{\bf Proof:}
Exhaust $Y_i$ by a sequence of compact subsets $K_{i,n}$ such that
$K_{i,n}\subset K_{i,n+1}$, $i=1,2$. There exist a subsequence of $(A_n,\psi_n)$ (still labeled by $n$), a sequence of gauge transformations $s_{i,n}$ defined on $K_{i,n}$, and monopoles $(A_i,\psi_i)$ on $Y_i$ such that $s_{i,n}\cdot (A_n,\psi_n)$ converges to $(A_i,\psi_i)$ in $C^\infty$
on any compact subset of $Y_i$. Note that $(A_i,\psi_i)$ are of finite
energy by Weitzenb\"ock formula.

First we show that there is an $r>0$ such that the geometric limits $(A_i,\psi_i)$ ($i=1,2$) are ``r-good''. Let $d$ be
the distance between the lattice of ``bad'' points and the line $R_2({\M}(Y_2))$ in $\h$. We simply take $r=\frac{d}{100}$.
Note that $\psi_2\equiv 0$ and $a_2(t)=a_{2,\infty}$ for all
$t\in [0,\infty)$, so $(A_2,\psi_2)$ is ``r-good''.
On the other hand, there are $t_0$ and $N$ such that
$|(D_{A_{1,1}}^{T^2}\psi_1,\psi_1)(t_0)|<\frac{1}{2}{\Gamma}(r)$ and $|(D_{A_{n,1}}^{T^2}\psi_n,\psi_n)(t_0)|<{\Gamma}(r)$ when $n>N$
(see Definition 2.3.9 for ${\Gamma}(r)$).
For large $n$, $a_{s_{2,n}\cdot (A_n)}(2L_n)$ is in 
${\hh}(4r)$, so is $a_n(2L_n)$. 
By Weitzenb\"ock formula, 
$
\int_{Y_2(t_0)}(|\nabla_{A_n}\psi_n|^2+\frac{1}{2}|\psi_n|^4)<
|(D_{A_{n,1}}^{T^2}\psi_n,\psi_n)(t_0)|<{\Gamma}(r)
$
when $n>N$, where $Y_2(t_0)=Y_2\setminus (2L_n+1-t_0,\infty)\times T^2$. 
So for large $n$, $a_n(t_0+1)$ is in ${\hh}(3r)$ (Lemma 2.3.10) and
so is $a_{s_{1,n}\cdot (A_n)}(t_0+1)$.  So $a_1(t_0+1)$ is in ${\hh}(2r)$. It follows easily that $(A_1,\psi_1)$ is ``r-good''.

 
Next we show that 
\begin{enumerate}
\item The $s_{i,n}'s$ can be chosen so that as $n\rightarrow \infty$, 
$s_{1,n}^{-1}|_{T^2}\cdot s_{2,n}|_{T^2}$ is in the component of identity of 
$Map(T^2,S^1)$. As a consequence, $s_{1,n}$ and $s_{2,n}$ can be 
extended to an $s_n\in {\G}(Y_{L_n})$.
\item $R_1(A_1,\psi_1)=R_2(A_2,\psi_2)$.
\end{enumerate}

Given any $\epsilon > 0$, pick $L_0$ large enough so that 
$$
\|(A_i,\psi_i)(L_0)-R_i(A_i,\psi_i)\|_{C^k(T^2)} < \epsilon, \hspace{2mm}
i=1,2.
$$
For large enough $n$, we have
$$
\|s_{1,n}\cdot (A_n,\psi_n)(L_0)-(A_1,\psi_1)(L_0)\|_{C^k(T^2)} < \epsilon 
$$
and
$$
\|s_{2,n}\cdot (A_n,\psi_n)(2L_n+1-L_0)-(A_2,\psi_2)(L_0)\|_{C^k(T^2)} < \epsilon.
$$
Let $\gamma$ be a generator of $H_1(T^2)$. Then we have
\begin{eqnarray*}
&      & |\int_{\gamma}(R_1(A_1,\psi_1)-R_2(A_2,\psi_2))|\\
& \leq & 3C\epsilon +
|\int_{\gamma}(s_{1,n}\cdot A_n(L_0)-s_{2,n}\cdot A_n(2L_n+1-L_0))|\\
& \leq & 3C\epsilon + 
|\int_{L_0}^{2L_n+1-L_0}\int_{\gamma} dA_n+2\pi i[s_{1,n}^{-1}|_{T^2}\cdot s_{2,n}|_{T^2}]|,
\end{eqnarray*}
where $[s_{1,n}^{-1}|_{T^2}\cdot s_{2,n}|_{T^2}]$ is the degree of the map
$s_{1,n}^{-1}|_{T^2}\cdot s_{2,n}|_{T^2}:\gamma\rightarrow S^1$.
On the other hand, we have estimates
\begin{eqnarray*}
|\int_{L_0}^{2L_n+1-L_0}\int_{\gamma} dA_n|
& \leq & C\int_{L_0}^{2L_n+1-L_0}\int_{\gamma} |\psi_n|^2\\
& \leq & C_1\int_{L_0-1}^{2L_n+2-L_0}\int_{T^2}(|\nabla_{A_n}\psi_n|^2+ \frac{1}{2}|\psi_n|^4)\\
& \leq & C_1(|(D_{A_{n,1}}^{T^2}\psi_n,\psi_n)(L_0-1)|\\
&      & +|(D_{A_{n,1}}^{T^2}\psi_n,\psi_n)(2L_n+2-L_0)|
\end{eqnarray*}
by Lemmas 2.3.6, 2.3.7 and Weitzenb\"ock formula, from which it follows that $R_1(A_1,\psi_1)=R_2(A_2,\psi_2)$ in ${\hh}^1(T^2)\otimes i{\bf R}
/({\bf Z}\oplus{\bf Z})$. It follows that $(A_i,\psi_i)$ can be modified by an element in ${\G}(Y_i)$ so that $R_1(A_1,\psi_1)=R_2(A_2,\psi_2)$
in ${\hh}^1(T^2)\otimes i{\bf R}$ (due to the fact that $Y$ is a homology 3-sphere) 
and $s_{1,n}^{-1}|_{T^2}\cdot s_{2,n}|_{T^2}$ is in the identity component of  $Map(T^2,S^1)$ for large $n$.
\hfill ${\Box}$

In the following discussion, we fix the number $r$ in Lemma 3.2.1, a weight $\delta$ small and a generic perturbation $(g,f,\mu)$ with $\mu$ sufficiently small. 

Recall the set ${{\S}(Y_1,Y_2)}$ of pairs $(\alpha_1,\alpha_2)\in\tilde{{\M}}^\ast(Y_1)\times\tilde{{\M}}(Y_2)$ such that
there are smooth representatives $(A_1,\psi_1)$ and $(A_2,\psi_2)$ satisfying
$R_1(A_1,\psi_1)=R_2(A_2,\psi_2)$. By Proposition 3.1.2 and Lemma 3.2.1, each 
pair $(\alpha_1,\alpha_2)$ in ${{\S}(Y_1,Y_2)}$ is ``r-good''. Therefore by Proposition 2.3.13 ( convergence of ``r-good''
monopoles ) and Proposition 2.2.2 (3), ${{\S}(Y_1,Y_2)}$ is compact and
hence consists of finitely many points. 

\sprop
For large enough $L$, the gluing map $T:{{\S}(Y_1,Y_2)}\rightarrow {\M}^\ast(Y_L)$
given by $T(\alpha_1,\alpha_2)=[T(A_L,\psi_L)]$ is one to one and onto.
\eprop 

Assume that a sequence of irreducible monopoles $(A_n,\psi_n)$ on $Y_{L_n}$ converges to geometric limits $(A_i,\psi_i)$ on $Y_i$ ($i=1,2$). 
Note that $\psi_2\equiv 0$ and
$(A_1,\psi_1)$ is irreducible since the (perturbed) Dirac operator at the 
reducible point on $Y_{L_n}$ is 
invertible for large $n$ and the norm is uniformly bounded from below
(Theorem B in \cite{CLM1}).
Our next goal is to show that for large enough $n$, $(A_n,\psi_n)$ is in the image of the gluing map $T$. This is done by showing that up to a gauge transformation the $L^2_1$ distance between $(A_n,\psi_n)$ and
the ``almost'' monopole $(A_{L_n},\psi_{L_n})$ is less than $K_1L_n^{-6}$
(see Proposition 3.1.2).

For simplicity we omit the subscript $n$ in the notation if no confusion is caused.
As in Lemma 2.3.11, there exists an $L_0>0$ such that for any $t\in [L_0,2L+1-L_0]$, we have
$$
\int_{T^2}|\psi|^2(t)\leq
e^{4\delta(L_0-t)}\int_{T^2}|\psi|^2(L_0)+
e^{4\delta(t-2L-1+L_0)}\int_{T^2}|\psi|^2(2L+1-L_0).
$$
For each $L>0$, fix a cut-off function $\rho_L$ which equals to one for $t\leq L$ and equals to zero for $t\geq L+1$. We construct an ``almost'' monopole $(\tilde{A_1},\tilde{\psi_1})$ on $Y_1$ as follows:
\begin{eqnarray*}
\tilde{A_1} & = & \rho_L(A-a(L+1))+a(L+1) \hspace{5mm}\mbox{on} \hspace{2mm} Y_1\setminus [L+1,\infty)\times T^2\\
\tilde{A_1} & = & a(L+1)\hspace{23mm} \mbox{on} \hspace{2mm} [L+1,\infty)\times T^2\\
\tilde{\psi_1} & = & \rho_L\psi \hspace{23mm} \mbox{on}\hspace{2mm}
Y_1\setminus [L+1,\infty)\times T^2\\
\tilde{\psi_1} & = & 0 \hspace{23mm} \mbox{on} \hspace{2mm} [L+1,\infty)\times T^2.
\end{eqnarray*}
(Note that we have omitted the subscript $n$ in the notation; here $A=A_n$ and $L=L_n$). Here $a(t)$ is the harmonic component of $A|_{\{t\}\times T^2}$.

The following estimate is straightforward.

\slm
$
\|(\ast d{\tilde{A_1}}+\tau(\tilde{\psi_1}),
D_{\tilde{A_1}}\tilde{\psi_1})\|_{L^2_{1,\delta}(Y_1)}\leq 
Ce^{-\delta L}
$
holds for $(\tilde{A_1},\tilde{\psi_1})$ on $Y_1$. 
\elm

Recall from Definition 2.1.10 that for $(A,\psi)\in{\A}^\ast$, 
$\nabla s_{(A,\psi)}:T{\B}^\ast_{(A,\psi)}\rightarrow {\E}_{(A,\psi)}$ is given by
$$
\nabla s_{(A,\psi)}(a,\phi)=
(\ast da+2\tau(\psi,\phi)-df(a,\phi), D_A\phi+a\psi+f(a,\phi)\psi)
$$
where $f(a,\phi)$ is the unique solution to the equation
$$
d^\ast df+f|\psi|^2=i\langle D_A\psi,i\phi\rangle_{Re}.
$$

\slm
For all sufficiently large $L$, $\nabla s_{(\tilde{A_1},\tilde{\psi_1})}: T{\B}^\ast_{(\tilde{A_1},\tilde{\psi_1})}\rightarrow{\E}_{(\tilde{A_1},\tilde{\psi_1})}$ is surjective. So there exists a bounded right inverse $P:{\E}_{(\tilde{A_1},\tilde{\psi_1})}\rightarrow T{\B}^\ast_{(\tilde{A_1},\tilde{\psi_1})}$ satisfying 
$$
\|P(a,\phi)\|_{{\A}}\leq K\|(a,\phi)\|_{L^2_{1,\delta}(Y_1)}
$$
for a constant $K$ independent of $L$ (see Definition 2.1.1 for the norm
$\|\hspace{1.5mm}\|_{{\A}}$).
\elm

\noindent{\bf Proof:}
Let $\Pi$ be the $L^2$-orthogonal projection onto ${\E}_{(A_1,\psi_1)}$, $\pi$ be the $L^2$-orthogonal
projection onto $T{\B}^\ast_{(\tilde{A_1},\tilde{\psi_1})}$ and $I$ be the
right inverse of $\nabla s_{(A_1,\psi_1)}$ ($I$ exists by the assumption that
${{\S}(Y_1,Y_2)}$ is regular).
For $(a,\phi)\in {\E}_{(\tilde{A_1},\tilde{\psi_1})}$, we have
$$
\nabla s_{(\tilde{A_1},\tilde{\psi_1})}\pi I \Pi(a,\phi)
=(a,\phi)+o(1)(a,\phi)
$$
as $L\rightarrow \infty$. Here the key point is that 
$d^\ast d + |\tilde{\psi_1}|^2$ is invertible and the norm of the inverse is
bounded uniformly in $L$ (Lemma 2.1.7).
\hfill ${\Box}$

Next we deform the ``almost'' monopoles $(\tilde{A_1},\tilde{\psi_1})$ to monopoles. Let $\Pi_1$ be the $L^2$ orthogonal projection onto
${\E}_{(\tilde{A_1},\tilde{\psi_1})}$. For any $(a,\phi)\in 
T{\B}^\ast_{(\tilde{A_1},\tilde{\psi_1})}$, we define
\begin{eqnarray*}
L(a,\phi)
& = & \Pi_1(\ast d{(\tilde{A_1}+a)}+\tau(\tilde{\psi_1}+\phi), D_{(\tilde{A_1}+a)}(\tilde{\psi_1}+\phi))\\
& = & (\ast d{\tilde{A_1}}+\tau(\tilde{\psi_1}), D_{\tilde{A_1}}\tilde{\psi_1}) + \nabla s_{(\tilde{A_1},\tilde{\psi_1})}(a,\phi) + \Pi_1 Q(a,\phi)
\end{eqnarray*}
where $Q(a,\phi)=(\tau(\phi),a\phi)$ satisfying
$$
\|Q(a_1,\phi_1)-Q(a_2,\phi_2)\|_{L^2_{1,\delta}}\leq 
C(\|(a_1,\phi_1)\|_{{\A}}+\|(a_2,\phi_2)\|_{{\A}})
\|(a_1,\phi_1)-(a_2,\phi_2)\|_{{\A}}.
$$

\slm
$L(a,\phi)=0$ implies that
$$
(\ast d{(\tilde{A_1}+a)}+\tau(\tilde{\psi_1}+\phi), D_{(\tilde{A_1}+a)}(\tilde{\psi_1}+\phi))=0
$$ when $\|(a,\phi)\|_{{\A}}$ is sufficiently small.
\elm

\noindent{\bf Proof:}
A similar argument as in the proof of Lemma 3.1.7. The key point is that 
$d^\ast d + |\tilde{\psi_1}|^2$ is invertible and the norm of the inverse is
bounded uniformly in $L$ (Lemma 2.1.7).
\hfill ${\Box}$

\slm
The ``almost'' monopole $(\tilde{A_1},\tilde{\psi_1})$ can be deformed to a
monopole $(\tilde{A_1}^\prime,\tilde{\psi_1}^\prime)$ 
such that 
$(\tilde{A_1}^\prime,\tilde{\psi_1}^\prime)- 
(\tilde{A_1},\tilde{\psi_1})\in T{\B}^\ast_{(\tilde{A_1},\tilde{\psi_1})}$ and
$
\|(\tilde{A_1}^\prime,\tilde{\psi_1}^\prime)- 
(\tilde{A_1},\tilde{\psi_1})\|_{{\A}}\leq Ce^{-\delta L}.
$
\elm

\noindent{\bf Proof:}
A similar argument as in the proof of Proposition 3.1.2. The fact that 
$d^\ast d + |\tilde{\psi_1}|^2$ is invertible and the norm of the inverse is bounded uniformly in $L$ (Lemma 2.1.7) is also used here
to get an estimate $\|\Pi_1 Q(a,\phi)\|_{L^2_{1,\delta}}\leq c
\|Q(a,\phi)\|_{L^2_{1,\delta}}$.
\hfill ${\Box}$


\noindent{\bf The Proof of Proposition 3.2.2:}

We need an estimate on the restriction of $(A,\psi)$ on the $Y_2$
side. Similarly we
construct ``almost'' monopoles $(\tilde{A_2},\tilde{\psi_2})$ on $Y_2$:
\begin{eqnarray*}
\tilde{A_2} & = &(1- \rho_L)(A-a(L))+a(L) \hspace{5mm}\mbox{on} \hspace{2mm} Y_2\setminus [L+1,\infty)\times T^2\\
\tilde{A_2} & = & a(L)\hspace{23mm} \mbox{on} \hspace{2mm} [L+1,\infty)\times T^2\\
\tilde{\psi_2} & = & (1-\rho_L)\psi \hspace{23mm} \mbox{on}\hspace{2mm}
Y_2\setminus [L+1,\infty)\times T^2\\
\tilde{\psi_2} & = & 0 \hspace{23mm} \mbox{on} \hspace{2mm} [L+1,\infty)\times T^2.
\end{eqnarray*}
By Weitzenb\"ock formula and the exponential decay estimate for the spinor, we have
$$
\int_{Y_2(L+1)}(|\nabla_{A}\psi|^2+\frac{1}{2}|\psi|^4)\leq
|(D^{T^2}_{A_{,1}}\psi,\psi)(L)|
\leq Ce^{-6\delta L}
$$
for a small $\delta>0$ where $Y_2(L+1)=Y_2\setminus (L+1,\infty)\times T^2$.
It then follows that
$$
\|\ast d\tilde{A_2}\|_{L^2_{\delta}}\leq Ce^{-\delta L}
\hspace{2mm} and \hspace{2mm}
\|\tilde{\psi_2}\|_{L^2_{1,\delta}}\leq Ce^{-\delta L}.
$$
Therefore the distance between $R_2(\tilde{A_2},0)$ and 
$[R_2](\tilde{{\M}}(Y_2))$ 
is controlled by $Ce^{-\delta L}$ (Lemmas 2.1.4, 2.2.5). On the other hand, the distance between $R_2(\tilde{A_2},0)$ and  $R_1(\tilde{A_1},\tilde{\psi_1})$ which is given by $|a(L+1)-a(L)|$ is also controlled by $Ce^{-\delta L}$ 
(Lemma 2.3.4 (d) and the exponential decay estimate for the spinor $\psi$).
So is the distance between 
$R_1(\tilde{A_1}^\prime,\tilde{\psi_1}^\prime)$ and $[R_2](\tilde{{\M}}(Y_2))$
by Lemma 3.2.6. By the 
assumption of transversality (Definition 3.1.1 (2)), we have
$$
R_1(T\tilde{{\M}}^\ast(Y_1)_{(A_1,\psi_1)})\bigcap
R_2(T\tilde{{\M}}(Y_2)_{(A_2,\psi_2)})=\{0\}.
$$
Then it follows that the distance between $R_1(\tilde{A_1}^\prime,\tilde{\psi_1}^\prime)$ and $R_1(A_1,\psi_1)$ is controlled by $Ce^{-\delta L}$. Since
$[R_1]: \tilde{{\M}}^\ast(Y_1)\rightarrow \h$ is an immersion at $[(A_1,\psi_1)]$
(Definition 3.1.1 (1)),
the distance between $[(\tilde{A_1}^\prime,\tilde{\psi_1}^\prime)]$ and $[(A_1,\psi_1)]$ is controlled by $Ce^{-\delta L}$. So is the distance
between $[(\tilde{A_1},\tilde{\psi_1})]$ and $[(A_1,\psi_1)]$ by Lemma 3.2.6.
The distance between $[(\tilde{A_2},\tilde{\psi_2})]$ and $[(A_2,\psi_2)]$ is also controlled by $Ce^{-\delta L}$.
Now it is easy to see that  
up to a gauge transformation $(A_n,\psi_n)$ is within an $L^2_1$ ball of
radius $Ce^{-\delta L_n}$ centered at the ``almost'' monopole 
$(A_{L_n},\psi_{L_n})$ for large enough $n$.
By Proposition 3.1.2, $(A_n,\psi_n)$ is in the image of the gluing map $T$.
On the other hand,  it follows from the ``weak''convergence of the gauge transformations that the gluing map $T: {{\S}(Y_1,Y_2)}\rightarrow
{\M}^\ast(Y_L)$ is also one to one. Hence the proposition is proved.

\section{Spectral flow, Maslov index and the gluing formula}

First we recall the basic relation between Maslov index and the spectral flow
of a one-parameter family of first-order, self-adjoint, elliptic differential operators of APS type on a stretched manifold. The basic references are \cite{CLM1} and \cite{CLM2}.

Let $M$ be a closed, oriented, smooth manifold that is decomposed into the union of two submanifolds $M_1, M_2$ by a co-dimension $1$, compact oriented submanifold ${\Sigma}$,
$$
M=M_1\bigcup M_2,\hspace{2mm} {\Sigma}=M_1\bigcap M_2=\partial M_1=\partial M_2.
$$
Equip $M$ with a Riemannian metric such that the hypersurface ${\Sigma}$ has a
collar neighborhood isometric to $(-1,1)\times{\Sigma}$, and
${\Sigma}={0}\times{\Sigma}, (-1,0)\times{\Sigma}\subset M_1$.
We stretch $M$ by inserting  cylinders
$[0,2L]\times{\Sigma}$ and obtain a family of manifolds
$M(L)$. Let $M_1(\infty), M_2(\infty)$ be the cylindrical end manifolds obtained by attaching $[0,\infty)\times{\Sigma}$ to $M_1$, and $(-\infty,0]\times{\Sigma}$ to
$M_2$. 

Let $D:{\Gamma}(E)\rightarrow{\Gamma}(E)$ be a first-order, self-adjoint, elliptic differential operator acting on the space of smooth sections of a real Riemannian vector bundle $E\rightarrow M$
which is of the APS type near ${\Sigma}$. More precisely, on $(-1,1)\times{\Sigma}$, $E$ is isometric to the
pull-back bundle $\pi^\ast E_0$ and $D$ can be written as
$$
D=\sigma(\frac{\partial}{\partial t}+D_0),
$$
where $\pi:(-1,1)\times{\Sigma}\rightarrow{\Sigma}$ is the projection,
$E_0\rightarrow {\Sigma}$ is a Riemannian vector bundle on ${\Sigma}$,
$\sigma: E_0\rightarrow E_0$ is a bundle isometry, and
$D_0$ is a first-order, self-adjoint, elliptic operator acting on ${\Gamma}(E_0)$.
Then $E$ and  $D$ naturally extend to a vector bundle
$E(L)$ and an operator $D(L)$ on the stretched manifold $M(L)$, and to
$E_j(\infty)$ and $D_j(\infty)$ on the cylindrical end manifold $M_j(\infty)$,
$j=1,2$.

Let $l_j$ be the space of limiting values of the extended $L^2$-solutions of
$D_j(\infty)=0$ over $M_j(\infty)$.
Denote ${{\K}}er\;D_0$ by ${\hh}$. Then we have (see \cite{CLM1})

\slm
\begin{enumerate}
\item ${\hh}$ is a symplectic vector space with the preferred symplectic form
$$
\{x,y\}=\int_{{\Sigma}}\langle x,\sigma y\rangle.
$$
\item $l_1,l_2$ are Lagrangian subspaces in ${\hh}$.
\end{enumerate}
\elm

We call $l_j$ the Lagrangian subspace associated to $D_j(\infty)$.

Let $E_1$, $E_2$ be the restriction of the vector bundle $E$ and $D_1$, $D_2$ be the restriction of the operator $D$ on the submanifolds $M_1$ and $M_2$ of $M$.
For any pair of Lagrangian subspaces $l_1,l_2$ of the symplectic vector space
${\hh}={{\K}}er\;D_0$, we have a pair of self-adjoint Fredholm operators
$D_1(l_1), D_2(l_2)$ defined with global boundary conditions:
$$
\begin{array}{c}
D_1(l_1):L^2_1(E_1,P_{+}\oplus l_1)\rightarrow L^2(E_1)\\
D_2(l_2):L^2_1(E_2,P_{-}\oplus l_2)\rightarrow L^2(E_2)
\end{array}
$$
where $P_{\pm}$ are the subspaces of $L^2(E_0)$ spanned by the eigenvectors of 
positive/negative eigenvalues of $D_0$, and the space $L^2_1(E_1,P_{+}\oplus l_1)$ is the $L^2_1$-Sobolev completion of smooth sections of bundle $E_1$
whose restrictions on ${\Sigma}$ lie in the space $P_{+}\oplus l_1$ and similarly
is the other space $L^2_1(E_2,P_{-}\oplus l_2)$ understood.

Each homotopy class (with fixed ends) of one-parameter families of pairs of
Lagrangian subspaces $(l_1(s),l_2(s)): a\leq s\leq b$ is associated 
with an integer which is called the Maslov index of $(l_1(s),l_2(s))$ and denoted by $Mas\{(l_1(s),l_2(s)): a\leq s\leq b\}$ (see \cite{CLM2},\cite{CLM3}).

The $(\epsilon_1,\epsilon_2)$-spectral flow is defined as follows.
Let $D(s):a\leq s \leq b$ be a family of real self-adjoint operators such that
for some fixed $\delta>0$ the total spectrum of $D(s)$ in the range of eigenvalues $\lambda$ with $|\lambda|<\delta$ is finite-dimensional and has
no essential spectrum. Furthermore, after taking into consideration of multiplicities, these eigenvalues $\lambda$ with $|\lambda|<\delta$ vary
continuously with respect to $s$. Let $\epsilon_1,\epsilon_2$ be real numbers
with $|\epsilon_1|<\delta, |\epsilon_2|<\delta$, such that $\epsilon_1$ is not
an eigenvalue of $D(a)$ and $\epsilon_2$ is not an eigenvalue of $D(b)$.
Then the $(\epsilon_1,\epsilon_2)$-spectral flow of $D(s):a\leq s\leq b$ is
equal to the number of times the eigenvalues $\lambda$ of $D(s)$ in the range $|\lambda|<\delta$ cross the line joining $(a,\epsilon_1)$ and $(b,\epsilon_2)$ from below, minus the number of times crossing from above (see \cite{CLM2} for details). The $(\epsilon,\epsilon)$-spectral flow will
be called briefly as $\epsilon$-spectral flow.

Let $D(s): a\leq s\leq b$ be a one-parameter family of first-order, self-adjoint, elliptic differential operators on $M$ which are of the APS type,
i.e. in the collar neighborhood $(-1,1)\times{\Sigma}$, 
$D(s)=\sigma(\frac{\partial}{\partial t}+D_0(s))$.
Furthermore, there exists a $\delta>0$ such that there are no eigenvalues of $D_0(s)$ in the range $(-\delta,0)$ and $(0,\delta)$, and ${\hh}={{\K}}er\;D_0(s)$
is a fixed symplectic vector space for $a\leq s\leq b$. A one-parameter family
of pairs of Lagrangian subspaces $(l_1(s),l_2(s)): a\leq s\leq b$ in ${\hh}$ is
said to satisfy the endpoint condition if $(l_1(s),l_2(s))$ is the pair of
Lagrangian subspaces associated to $(D_1(\infty)(s),D_2(\infty)(s))$ at the endpoints $s=a,b$.  

The basic relation between Maslov index and spectral flow is given by the 
following

\sthm
(Theorem C in \cite{CLM2})

There exists an $L_0>0$ such that for any choice of smoothly varying pairs of
Lagrangian subspaces $(l_1(s),l_2(s)): a\leq s\leq b$ satisfying the endpoint
condition, for all $L>L_0$, the $(L^{-2})$-spectral flow of $D(s)(L)$ on $M(L)$
for $a\leq s\leq b$ equals to 
$$
\sum_{j=1}^{2} SF_\epsilon \{D_j(s)(l_j(s)): a\leq s\leq b\} + Mas\{(l_1(s),l_2(s)): a\leq s\leq b\}
$$
where $SF_\epsilon \{D_j(s)(l_j(s)): a\leq s\leq b\}$ is the $\epsilon$-spectral flow of $D_j(s)(l_j(s)): a\leq s\leq b$. Here $\epsilon >0$ is chosen so that the eigenvalues of $D_j(s)(l_j(s))$ in the range $[-\epsilon,\epsilon]$ consist of at most zero eigenvalues for the endpoints $s=a,b$.
\ethm

Now let's go back to our own problem. Suppose that
$Y$ is an oriented integral homology 3-sphere that is decomposed as $Y=Y_1\bigcup_{T^2} Y_2$ with $Y_1$ being the complement of a tubular neighborhood of a knot and $Y_2=D^2\times S^1$. $Y$ carries a Riemannian metric such that a collar neighborhood of $T^2$ is orientedly isometric to $(-1,1)\times {\bf R}/2\pi{\bf Z}\times{\bf R}/2\pi{\bf Z}$ with $(-1,0)\times {\bf R}/2\pi{\bf Z}\times{\bf R}/2\pi{\bf Z}\subset Y_1$, where we assume that the first and second factors in ${\bf R}/2\pi{\bf Z}\times{\bf R}/2\pi{\bf Z}$ represent the longitude and meridian respectively (\cite{AM}),
and the metric on $Y_2$ has non-negative, somewhere positive scalar curvature.  By inserting cylinders $[0,2L+1]\times T^2$, we obtain a family of stretched versions $Y_L$ of $Y$. We also use $Y_1$ and $Y_2$ to denote the corresponding cylindrical end manifolds if no confusion occurs.

The basic result we've obtained so far (Proposition 3.2.2) is that for a
large enough $L$, the irreducible Seiberg-Witten moduli space ${\M}^\ast(Y_L)$
of $Y_L$ is identified via the gluing map $T$ with the set of ``intersection
points'' ${{\S}(Y_1,Y_2)}$. Here ${{\S}(Y_1,Y_2)}$ consists of the pairs $(\alpha_1,\alpha_2)\in\tilde{{\M}}^\ast(Y_1)\times\tilde{{\M}}(Y_2)$ such that there are smooth
representatives $(A_1,\psi_1)$ and $(A_2,\psi_2)$ having the same limiting
value, i.e. $R_1(A_1,\psi_1)=R_2(A_2,\psi_2)$. Our next goal is to orient $\tilde{{\M}}^\ast(Y_1)$ and $\tilde{{\M}}(Y_2)$ appropriately so that their 
``intersection number'' 
$\# {{\S}(Y_1,Y_2)}$ equals to the Seiberg-Witten invariant $\chi(Y_L)$ as
the oriented sum of the points in the moduli space ${\M}^\ast(Y_L)$. This is referred to as the gluing formula of $\chi$.

Fix a generic perturbation $(g,f,\mu)$ compactly supported on the $Y_1$ side according to Proposition 2.2.2 and thereafter omit it in the discussion for simplicity. Assume that $L_0$ is large enough so that Proposition 3.2.2 holds for $Y_{L_0}$. Pick a smooth section $\phi$ of the spinor bundle $W\rightarrow Y_{L_0}$ which is compactly supported in $Y_1\setminus [0,\infty)\times T^2$ and satisfies 
$
((D_g+f)^{-1}(i\phi),(i\phi))<0.
$
Then by Lemma 1.2.2, for small enough $t>0$, the self-adjoint operator
(on $Y_{L_0}$)
$$
{\K}_{(t,\phi)}=\left(\begin{array}{ccc}
D_g+f & 0 & 0\\
0 & \ast d & -d\\
0 & -d^\ast & 0
\end{array} \right)
+t\left(\begin{array}{ccc}
0 & \phi\cdot & \phi\cdot \\
2\tau(\phi,\cdot ) & 0 & 0\\
i\langle i\phi,\cdot \rangle_{Re} & 0 & 0
\end{array} \right)
$$
acting on ${\Gamma}(W\oplus (\Lambda^1\oplus\Lambda^0)\otimes i{\bf R})$
is invertible and has one small eigenvalue 
$$\lambda_t\sim 
-((D_g+f)^{-1}(i\phi),(i\phi))t^2>0.
$$ 
According to Definition 1.2.4, the Euler characteristic $\chi(Y_{L_0})$
is defined by
$$
\chi(Y_{L_0})=\sum_{\beta\in{\M}^\ast(Y_{L_0})} \mbox{sign}\beta, \hspace{2mm} \mbox{where}
\hspace{2mm}
\mbox{sign}\beta=(-1)^{SF({\K}_\beta,{\K}_{(t,\phi)})}
$$ 
for small $t>0$ ($SF$ denotes the spectral flow).
Here if $\beta$ is represented by $(A,\psi)$, then
$$
{\K}_\beta={\K}_{(A,\psi)}=\left(\begin{array}{ccc}
D_A & \psi\cdot & \psi\cdot\\
2\tau(\psi,\cdot ) & \ast d & -d\\
i\langle i\psi,\cdot \rangle_{Re} & -d^\ast & 0
\end{array} \right).
$$
Let $(A_{L_0},\psi_{L_0})$ be the ``almost'' monopole being deformed to $(A,\psi)$ under the gluing map $T$ (note that $\psi_{L_0}$ is compactly supported in $Y_1\setminus (L_0+1,\infty)\times T^2$). It is obvious that ${\K}_\beta$ can be replaced by ${\K}_{(A_{L_0},\psi_{L_0})}$ for the purpose of spectral flow calculation.
For any $L>0$, we insert cylinders of lengths $2L$ into $Y_{L_0}$ and obtain a family of manifolds $Y_{L_0,L}$ and operators ${\K}_{(A_{L_0},\psi_{L_0})}(L)$
on $Y_{L_0,L}$ from ${\K}_{(A_{L_0},\psi_{L_0})}$ in the obvious way.

\slm
For large enough $L_0$, ${\K}_{(A_{L_0},\psi_{L_0})}(L)$
are invertible for any $L>0$. In particular, the spectral flow between
${\K}_{(A_{L_0},\psi_{L_0})}$ and ${\K}_{(A_{L_0},\psi_{L_0})}(L)$ is zero for any
$L>0$.
\elm

\noindent{\bf Proof:}
For large enough $L_0>0$, the operators ${\K}_{(A_{L_0},\psi_{L_0})}(L)$ are
invertible for all $0<L<L_0$ by Theorem 4 in Appendix B. Suppose that
${\K}_{(A_{L_0},\psi_{L_0})}(L)$ has a non-zero kernel for some $L\geq L_0$, i.e.
there is an $x\neq 0$ such that ${\K}_{(A_{L_0},\psi_{L_0})}(L)x=0$. On the
inserted cylinder, ${\K}_{(A_{L_0},\psi_{L_0})}(L)$ has the form
$
I(\frac{\partial}{\partial t}+B)
$
where
$$
I=\left(\begin{array}{clcr}
dt & 0 & 0 & 0 \\
0 & \ast_{T^2} & 0 & 0 \\
0 & 0 & 0 & -1 \\
0 & 0 & 1 & 0
\end{array} \right)
\hspace{2mm} \mbox{and} \hspace{2mm} 
B=\left(\begin{array}{clcr}
D^{T^2}_a & 0 & 0 & 0\\
0 & 0 & -\dd & \ast\dd \\
0 & -\dd^\ast & 0 & 0\\
0 & -\ast\dd & 0 & 0
\end{array} \right)
$$
acting on ${\Gamma}(W_0\oplus(\Lambda^1\oplus\Lambda^0\oplus\Lambda^0(T^2))\otimes
i{\bf R})$. Here $W_0$ is the total spinor bundle over $T^2$, and $D^{T^2}_a$ is an
invertible twisted Dirac operator. It follows that
$x$ can be decomposed as 
$x=x_0+x_{+}+x_{-}$ with $x_0\in Ker\;B$ constant in $t$ and $x_{\pm}$ have
exponential decay to the right/left. Take a cut-off function $\gamma$ in the
middle of the inserted cylinder, define $y_{\pm}$ on the cylindrical end
manifolds $Y_1$/$Y_2$ by:
$$
y_{+} = \gamma(x-x_0)+x_0, \hspace{2mm}
y_{-} = (1-\gamma)(x-x_0)+x_0.
$$
Then it follows that 
$$
\|{\K}_{(A_{L_0},\psi_{L_0}),1}(\infty)y_{+}\|_{L^2_\delta}
\leq c e^{-\delta L}(\|y_{+}-x_0\|_{L^2_\delta}+\|y_{-}-x_0\|_{L^2_\delta}) 
$$ 
and 
$$
\|{\K}_{(A_{L_0},\psi_{L_0}),2}(\infty)y_{-}\|_{L^2_\delta}
\leq c e^{-\delta L}(\|y_{+}-x_0\|_{L^2_\delta}+\|y_{-}-x_0\|_{L^2_\delta})
$$ 
for some small $\delta>0$ and a constant $c$. Here ${\K}_{(A_{L_0},\psi_{L_0}),j}(\infty)$ 
is the corresponding operator on the cylindrical end manifold $Y_j$, $j=1,2$.
On the other hand, observe that $y_{+}$ and $y_{-}$ have the same limiting value
$x_0$ and for all large enough $L_0$, the Lagrangian subspaces associated to
${\K}_{(A_{L_0},\psi_{L_0}),j}(\infty)$ ($j=1,2$) are transversal to each other
with angles larger than a fixed number ( due to the fact that ${{\S}(Y_1,Y_2)}$ is regular). 
Then the above estimates yield 
$$
\|x_0\|\leq c_1 e^{-\delta L}(\|y_{+}-x_0\|_{L^2_\delta}+\|y_{-}-x_0\|_{L^2_\delta}).
$$
Since both of
${\K}_{(A_{L_0},\psi_{L_0}),1}(\infty)$ and ${\K}_{(A_{L_0},\psi_{L_0}),2}(\infty)$  have no $L^2$ kernels, we have estimates
$$
\|y_{\pm}-x_0\|_{L^2_\delta}\leq c_2 
e^{-\delta L}(\|y_{+}-x_0\|_{L^2_\delta}+\|y_{-}-x_0\|_{L^2_\delta})
$$
which imply that for large $L_0$ (therefore $L\geq L_0$ large) $y_{\pm}$ vanish
identically, contradicting the assumption that $x\neq 0$. Therefore the lemma
is proved.
\hfill ${\Box}$

The operators considered here have the APS form 
$
I(\frac{\partial}{\partial t}+B)
$
on the inserted cylinder where
$$
I=\left(\begin{array}{clcr}
dt & 0 & 0 & 0 \\
0 & \ast_{T^2} & 0 & 0 \\
0 & 0 & 0 & -1 \\
0 & 0 & 1 & 0
\end{array} \right)
\hspace{2mm} \mbox{and} \hspace{2mm} 
B=\left(\begin{array}{clcr}
D^{T^2}_a & 0 & 0 & 0\\
0 & 0 & -\dd & \ast\dd \\
0 & -\dd^\ast & 0 & 0\\
0 & -\ast\dd & 0 & 0
\end{array} \right).
$$
The symplectic vector space ${\hh}=Ker\;B$ is ${\hh}^1(T^2)\otimes i{\bf R}\oplus i{\bf R}\oplus i{\bf R}$. Let's fix the notation about ${\hh}$ first. Recall that $T^2$ is orientedly isometric to ${\bf R}/2\pi{\bf Z}\times {\bf R}/2\pi{\bf Z}$ (longitude, meridian). Let $(x,y)$ be the oriented coordinates, then we orient ${\hh}^1(T^2)\otimes i{\bf R}$ by $idx\wedge idy$. Furthermore, the 3rd component of ${\hh}$ corresponds
to the $dt$-component of the 1-forms and the 4th one is from the Lie algebra of the gauge group. 
   
Let ${\K}_0$ (acting on ${\Gamma}(W\oplus(\Lambda^1\oplus\Lambda^0)\otimes i{\bf R})$)
be the operator at the reducible point $(0,0)$:
$$
{\K}_0=\left(\begin{array}{ccc}
D_g+f & 0 & 0\\
0 & \ast d & -d\\
0 & -d^\ast & 0
\end{array} \right).
$$
The corresponding operators ${\K}_{0,j}(\infty)$ on the 
cylindrical end manifolds $Y_j$ have no $L^2$ kernels and the associated 
Lagrangian subspaces of ${\K}_{0,1}(\infty)$ and ${\K}_{0,2}(\infty)$ are spanned by
$(idy,(0,0,0,1))$ and $([R_2](T\tilde{{\M}}(Y_2)),(0,0,0,1))$ respectively. Note that $[R_2](T\tilde{{\M}}(Y_2))$ is transversal to $idy$ since $Y$ is a homology 3-sphere.

Now we are ready to orient the moduli spaces $\tilde{{\M}}^\ast(Y_1)$ and $\tilde{{\M}}(Y_2)$.
Assume that $\alpha_1\in\tilde{{\M}}^\ast(Y_1)$ is represented by $(A,\psi)$.
For any vector $V\in\h$ with positive $idx$-component which is not in $R_1(T\tilde{{\M}}(Y_1)_{(A,\psi)})$, let $v\in T\tilde{{\M}}(Y_1)_{(A,\psi)}$ such
that $V\wedge R_1(v)=idx\wedge idy$.
Pick an $L_0>0$ and cut down $(A,\psi)$ at $L_0+1$. Denote the result by 
$(A_{L_0},\psi_{L_0})$. We assume that $L_0$ is large enough so that the Lagrangian subspace associated to ${\K}_{(A_{L_0},\psi_{L_0}),1}(\infty)$ is 
transversal to the Lagrangian subspace spanned by $(V,(0,0,0,1))$.
Connect $(A_{L_0},\psi_{L_0})$ to the reducible point $(0,0)$ by a path $(A,\psi)_s$ which is constant in $t$ on $[L_0+1,\infty)\times T^2$. Choose a smooth path of Lagrangian subspaces $l_1(s)$ which equals to the Lagrangian subspace $(idy,(0,0,0,1))$ associated to ${\K}_{0,1}(\infty)$ or the Lagrangian subspace associated to ${\K}_{(A_{L_0},\psi_{L_0}),1}(\infty)$ at the endpoints of the path. 

\sde
\begin{enumerate}
\item The orientation of $\tilde{{\M}}^\ast(Y_1)$ at $\alpha_1=[(A,\psi)]$ is determined by the tangent vector $(-1)^{m} v$.
Here $m$ is the sum of the $(\epsilon)$-spectral flow of operators
${\K}_{(A,\psi)_s,1}(L_0+1)(l_1(s))$ (for a small $\epsilon > 0$) and the Maslov index Mas$\{(l_1(s),l_V)\}$, where $l_V$ is the Lagrangian subspace spanned by $(V,(0,0,0,1))$. 
\item The orientation of $\tilde{{\M}}(Y_2)$ is determined so that the positive direction of $[R_2](T\tilde{{\M}}(Y_2))$ has positive $idx$-component. Note that $[R_2](T\tilde{{\M}}(Y_2))$ is transversal to $idy$-axis since $Y$ is a homology 3-sphere.
\end{enumerate}
\ede

\slm
The orientation on $\tilde{{\M}}^\ast(Y_1)$ is well-defined, which induces an
orientation on ${\M}^\ast(Y_1)$ via the ${\bf Z}$-fold covering map $\tilde{{\M}}^\ast(Y_1)\rightarrow{\M}^\ast(Y_1)$.
\elm

\noindent{\bf Proof:}
We need to prove that the orientation of $\tilde{{\M}}^\ast(Y_1)$
is independent of the choice of $\alpha_1$ (and its representatives $(A,\psi)$), the vector $V\in\h$, the cut-off point $L_0$, the path $(A,\psi)_s$ and the path of Lagrangian subspaces $l_1(s)$.

First of all,
the independence on the choice of cut-off point $L_0$, the path $(A,\psi)_s$ and the path of Lagrangian subspaces $l_1(s)$ follows easily from Theorem 3.3.2.
Secondly, suppose that two monopoles $(A_1,\psi_1)$ and $(A_2,\psi_2)$ are in the same component. Join them by a path of monopoles $(A_s,\psi_s)$ and then cut down the path at $L_0+1$ for sufficiently large $L_0$ (still denote the path by $(A_s,\psi_s)$). Let $l(s)$ be the Lagrangian subspace associated to ${\K}_{(A_s,\psi_s),1}(\infty)$.
Then the $(\epsilon)$-spectral flow of ${\K}_{(A_s,\psi_s),1}(L_0+1)(l(s))$ is zero because $\tilde{{\M}}^\ast(Y_1)$ is immersed into $\h$ so that ${\K}_{(A_s,\psi_s),1}(\infty)$ have no $L^2$-kernels for large enough $L_0$.
On the other hand, since $(A_s,\psi_s)$ is irreducible so that the 3rd component of $l(s)$ is non-zero, Mas$\{(l(s),(V,(0,0,0,1))\}$ (mod 2) equals to
the sign change of $V\wedge R_1(v_s)$ where $v_s$ is a smooth tangent vector field in $T\tilde{{\M}}^\ast(Y_1)$ along the path $(A_s,\psi_s)$. So the orientation at $(A_1,\psi_1)$ and the orientation at $(A_2,\psi_2)$ are compatible. 
Finally, suppose that $V_1,V_2\in \h$ are two different vectors used in the 
definition. Then Mas$\{(l_1(s),(V_1,(0,0,0,1))\}-$Mas$\{(l_1(s),(V_2,(0,0,0,1))\}$
(mod 2) equals to the sign change from $V_1\wedge R_1(v)$ to $V_2\wedge R_1(v)$
for any $v\in T\tilde{{\M}}(Y_1)_{(A,\psi)}$, which implies that the orientation of $\tilde{{\M}}^\ast(Y_i)$ at $[(A,\psi)]$ is independent of the choice of the vector $V$. Therefore we have proved that the orientation of $\tilde{{\M}}^\ast(Y_1)$ is well-defined.


Next we prove that the ${\bf Z}$-fold covering map $\tilde{{\M}}^\ast(Y_1)\rightarrow{\M}^\ast(Y_1)$ induces an orientation on ${\M}^\ast(Y_1)$. Suppose that $(A_1,\psi_1)$ and $(A_2,\psi_2)$ are gauge equivalent by a gauge transformation $s_1$ not in the identity component of ${\G}(Y_1)$. Pick an $L_0$ large enough and cut down $(A_1,\psi_1)$ at $L_0+1$ 
and still denote it by $(A_1,\psi_1)$ (we can assume that $s_1$ is constant in $t$ on $[L_0+1,\infty)\times T^2$). Connect $(A_1,\psi_1)$ with the reducible 
point $(0,0)$ by a path $(A_s,\psi_s)$ (so $s_1\cdot (A_s,\psi_s)$ is a path joining $s_1\cdot (A_1,\psi_1)=(A_2,\psi_2)$ with $(-s_1^{-1}ds_1,0)$). Then 
the $(\epsilon)$-spectral flow of ${\K}_{(A_s,\psi_s),1}(L_0+1)(l_1(s))$ equals 
to that of ${\K}_{s_1\cdot (A_s,\psi_s),1}(L_0+1)(l_1(s))$ where $l_1(s)$ is a path of Lagrangian subspaces which equals to the associated 
Lagrangian subspace of ${\K}_{(A_s,\psi_s),1}(\infty)$ at the endpoints. On the other hand, the $(\epsilon)$-spectral flow of ${\K}_{(-us_1^{-1}ds_1,0),1}(L_0+1)(l_3): 0\leq u\leq 1$ is even (Dirac operators are complex linear) where the Lagrangian subspace $l_3$ is spanned by $(idy,(0,0,0,1))$. Now it is easy to see that the orientation at $(A_1,\psi_1)$ and $(A_2,\psi_2)$ are compatible.
So the lemma is proved.
\hfill ${\Box}$

Now we are ready to define the ``intersection number'' $\#{{\S}(Y_1,Y_2)}$ and prove the gluing formula.

\sde
\begin{enumerate}
\item For any $(\alpha_1,\alpha_2)\in{{\S}(Y_1,Y_2)}$, let $e_j$ be the positively oriented tangent vector of $\tilde{{\M}}(Y_j)$ at $\alpha_j$ ($j=1,2$). Then the sign of $(\alpha_1,\alpha_2)$ is the sign of $[R_1]e_1\wedge [R_2]e_2$ with respect to $idx\wedge idy$.
\item $\#{{\S}(Y_1,Y_2)}=\sum_{(\alpha_1,\alpha_2)\in{{\S}(Y_1,Y_2)}}{sign}(\alpha_1,\alpha_2)$.
\end{enumerate}
\ede

\sthm
(Gluing Formula)

$\chi(Y_L)=\#{{\S}(Y_1,Y_2)}$ for sufficiently large $L>0$.
\ethm

\noindent{\bf Proof:}
Let $(A_{L_0},\psi_{L_0})$ be the ``almost'' monopole being deformed to
$\beta\in{\M}^\ast(Y_{L_0})$. By Lemma 3.3.3, $\mbox{sign}\beta=(-1)^{m_1+1}$ where $m_1$ is the $L^{-2}$-spectral flow between ${\K}_{(A_{L_0},\psi_{L_0})}(L)$ and ${\K}_0$ for sufficiently large $L$. By Theorem 3.3.2, $m_1$ is equal to 
$$
\sum_{j=1}^{2} SF_\epsilon \{{\K}_{(A,\psi)_s,j}(L_0+1)(l_j(s))\} + Mas\{(l_1(s),l_2(s))\}
$$
for any choice of $(A,\psi)_s$ joining  $(A_{L_0},\psi_{L_0})$ with the reducible point $(0,0)$ and any choice of a path of Lagrangian subspaces $(l_1(s),l_2(s))$ satisfying the endpoint condition.
Here $SF_\epsilon \{{\K}_{(A,\psi)_s,j}(L_0+1)(l_j(s))\}$ is the $\epsilon$-spectral flow of ${\K}_{(A,\psi)_s,j}(L_0+1)(l_j(s))$ for some small $\epsilon>0$. We choose $(A,\psi)_s$ such that $\psi_s$ is identically zero on the $Y_2$ side, and choose $l_2(s)=l_2$ to be the Lagrangian subspace spanned by $([R_2](T\tilde{{\M}}(Y_2)),(0,0,0,1))$.
Then the $\epsilon$-spectral flow of ${\K}_{(A,\psi)_s,2}(L_0+1)(l_2)$ is even.
On the other hand, suppose $\beta=T(\alpha_1,\alpha_2)$. Let $e_j$ be the positively oriented tangent vector of $\tilde{{\M}}(Y_j)$ at $\alpha_j$ ($j=1,2$). Then by taking $V=[R_2]e_2$ in Definition 3.3.4, we have $\mbox{sign}(\alpha_1,\alpha_2)$ equals to the sign of $[R_1]e_1\wedge [R_2]e_2=(-1)^m [R_1]v\wedge [R_2]e_2=(-1)^{m+1}idx\wedge idy=(-1)^{m+1}$. Here $m$ and $v$ are referred to Definition 3.3.4.
The theorem follows from the relation $m\equiv m_1$ mod $2$.

\appendix

\chappendix{APPENDIX A}

\quad The purpose of this appendix is to find out for what $a\in\h$ the twisted Dirac operator $D_a=D+a$ is not invertible. Here $D$ is the Dirac operator on
$T^2$ associated to a given spin structure and the flat metric.

First of all, let's recall some basic facts about the spin structures on the torus $T^2$. There are two equivalent descriptions of spin
structures. Topologically, a spin structure on $T^2$ is a framing of
its stabilized tangent bundle $TT^2\oplus\epsilon$ (a homotopy equivalence class
of trivializations). There are four different spin structures on $T^2$ which
are parameterized by $H^1(T^2,{\bf Z}_2)={\bf Z}_2\oplus{\bf Z}_2$. It is well-known that 
among these four different spin structures, three of them are spin boundaries,
i.e.  spin structures induced from a spin 3-manifold bounded by the torus. The only one left which is not a spin boundary is usually called the
Lie group framing. Assume that $T^2={\bf R}/2\pi{\bf Z}\times{\bf R}/2\pi{\bf Z}$ and let $(\frac{\partial}{\partial x},\frac{\partial}{\partial y})$ be the tangent vectors of the circles. For $(k,l)=(0,0),(0,1),(1,0),(1,1)$, the following formula defines four different framings $\xi_{(k,l)}$
of the tangent bundle $TT^2$ which induce all the spin structures on $T^2$ (framings of $TT^2\oplus\epsilon$):
$$
\xi_{(k,l)}(x,y)=\left(\begin{array}{cc}
\cos(kx+ly) & -\sin(kx+ly)\\
\sin(kx+ly) & \cos(kx+ly)
\end{array} \right)
\left(\begin{array}{c}
\frac{\partial}{\partial x}\\\frac{\partial}{\partial y}
\end{array} \right), \hspace{4mm}
(x,y)\in T^2.
$$ 
The Lie group framing is $\xi_{(0,0)}$. See \cite{K1} for details.

The geometric aspect of spin structures is related to the groups $Spin(n)$.
The groups $Spin(n)$ sit inside the n-dimensional Clifford algebras $Cl(n)$
and double cover the groups $SO(n)$. Let $\pi: Spin(n)\rightarrow SO(n)$ be
the double covering map. Equip the torus $T^2$ with a Riemannian metric,
assuming that it is the product metric for simplicity. Let $P_{SO(2)}$ be the
$SO(2)$ principal bundle to which the tangent frame bundle of $T^2$ is reduced.
A spin structure on $T^2$ is then defined to be an equivalence class of liftings
of the principal bundle $P_{SO(2)}$ to a $Spin(2)$ principal bundle $P_{Spin(2)}$, i.e. $P_{Spin(2)}\stackrel{\pi}{\rightarrow} P_{SO(2)}$ such that $\pi$ restricts to the double covering map on each fiber. Two liftings
$P_{Spin(2)}^{(1)}\stackrel{\pi_1}{\rightarrow} P_{SO(2)}$ and
$P_{Spin(2)}^{(2)}\stackrel{\pi_2}{\rightarrow} P_{SO(2)}$ are said to be
equivalent if and only if there is a bundle isomorphism $i$ such that the
following diagram commutes:
$$
\begin{array}{ccc}
P_{Spin(2)}^{(1)} & \stackrel{i}{\rightarrow} & P_{Spin(2)}^{(2)}\\
\downarrow \pi_1 &  & \downarrow \pi_2 \\
P_{SO(2)} & \stackrel{identity}{\rightarrow} & P_{SO(2)}
\end{array}
$$
For each spin structure $P_{Spin(2)}\stackrel{\pi}{\rightarrow} P_{SO(2)}$, 
there is a canonically associated spinor bundle $W=W^+\oplus W^-$ on $T^2$, where $W^\pm=P_{Spin(2)}\times_{\varrho_\pm} {\bf C}$. The representations $\varrho_\pm: Spin(2)\rightarrow U(1)$ are distinguished by the conditions $\varrho_\pm(e_1e_2)=\mp i$ for any orthonormal basis $(e_1,e_2)$ of ${\bf R}^2$.

The topological and geometrical descriptions of spin structures on $T^2$ are
related in the following way. The spin structure induced by the trivialization
$\xi_{(k,l)}$ ($k,l=0,1$) corresponds to the unique equivalence class of 
liftings $P_{Spin(2)}^{(k,l)}\rightarrow P_{SO(2)}$ for which the trivialization
$\xi_{(k,l)}$ of $P_{SO(2)}$ can be lifted to a trivialization 
$\tilde{\xi_{(k,l)}}$ of $P_{Spin(2)}^{(k,l)}$, which further induces
trivializations for the spinor bundles $W^\pm$ and $W$.

{\bf Theorem:}{\em

Assume that $T^2={\bf R}/2\pi{\bf Z}\times{\bf R}/2\pi{\bf Z}$ carries the product metric
and $(\frac{\partial}{\partial x},\frac{\partial}{\partial y})$ is the oriented
orthonormal basis. For $(k,l)=(0,0),(0,1),(1,0),(1,1)$, define trivializations
$\xi_{(k,l)}$ of $TT^2$ by the following formula:
$$
\xi_{(k,l)}(x,y)=\left(\begin{array}{cc}
\cos(kx+ly) & -\sin(kx+ly)\\
\sin(kx+ly) & \cos(kx+ly)
\end{array} \right)
\left(\begin{array}{c}
\frac{\partial}{\partial x}\\\frac{\partial}{\partial y}
\end{array} \right), \hspace{4mm} (x,y)\in T^2.
$$
Then within the induced trivialization $\tilde{\xi_{(k,l)}}$ of
the spinor bundles associated to the spin structure given by $\xi_{(k,l)}$, the
Dirac operator $D^{(k,l)}$ is given by the following formula
$$
D^{(k,l)}\left(\begin{array}{c}
u\\v
\end{array} \right)=dx\left(\frac{\partial}{\partial x}+\frac{i}{2}\left(\begin{array}{cc}
k & 0\\
0 & -k
\end{array} \right)\right)\left(\begin{array}{c}
u\\v
\end{array} \right)+dy\left(\frac{\partial}{\partial y}+\frac{i}{2}\left(\begin{array}{cc}
l & 0\\
0 & -l
\end{array} \right)\right)\left(\begin{array}{c}
u\\v
\end{array} \right),
$$
where $u,v$ are complex valued functions on $T^2$.
As a consequence, for $a\in\h$, the twisted Dirac operator
$D^{(k,l)}_a=D^{(k,l)}+a$ is invertible unless $a=\frac{i}{2}(kdx+ldy)+
sidx+tidy$ for some integers $s$ and $t$, and $\dim_{\bf C}\ker\;D^{(k,l)}_a=2$ if
$a=\frac{i}{2}(kdx+ldy)+sidx+tidy$ for some integers $s$ and $t$}.

{\bf Remark:}
The lattice 
$$
B_{(k,l)}=\{a|a=\frac{i}{2}(kdx+ldy)+sidx+tidy, s,t\in{\bf Z}\}
$$
is called the lattice of ``bad'' points for the spin structure $\xi_{(k,l)}$.

{\bf Proof:}
In general, if $(e_1, e_2, ..., e_n)$ is an oriented local orthonormal frame,
then within the induced trivialization of the spinor bundles, the induced
connection is given by $-\frac{1}{2}\sum_{i<j}\omega_{ij}e_ie_j$, where
$\omega_{ij}$ is given by the formula $\nabla e_j=e_i\omega_{ij}$ (see \cite{LMi} for details). Back to our case of the torus, let $\xi_{(k,l)}(x,y)=(e_1,e_2)$, and 
$\nabla e_2=e_1\omega_{12}$, then $\omega_{12}=kdx+ldy$ by direct calculation. The theorem follows easily from this.

\chappendix{APPENDIX B}

\quad
The purpose of this appendix is to give an estimate on the lowest eigenvalue
of certain self-adjoint elliptic operators on a manifold containing long necks,
a technical result needed in Chapter 3. See \cite{C2}.
 
Let $X$ be an oriented Riemannian manifold with a cylindrical end modeled on
$Y$, i.e. there exists a compact subset $K$ such that $X\setminus K$ is 
isometric to $(-1,\infty)\times Y$.
Let $E$ be a cylindrical Riemannian vector bundle over $X$. By definition,
there is a Riemannian vector bundle $E_0$ over $Y$ such that $E$ is isometric
to $\pi^\ast E_0$ on the cylindrical end $(-1,\infty)\times Y$, where $\pi : (-1,\infty)\times Y\rightarrow Y$ is the natural projection. Assume that $D:
\Gamma(E)\rightarrow\Gamma(E)$ is
a first order formally self-adjoint elliptic operator on $X$, which takes the
following form on the cylindrical end $(-1,\infty)\times Y$
$$
D=I(\frac{\partial}{\partial t}+A)
$$
where $I$ is a bundle automorphism of $E_0$ which preserves its inner
product, and $A: \Gamma(E_0)\rightarrow\Gamma(E_0)$ is an elliptic operator
on $Y$ independent of $t$. The self-adjointness of $D$ implies that $I$ and
$A$ satisfy the following conditions:
$$
I^2=-1,\hspace{2mm} I^\ast=-I, \hspace{2mm} A^\ast=A,  \hspace{2mm}
IA+AI=0.
$$
Note that the spectrum of $A$ is symmetric about the origin and the 
automorphism $I$ maps $E_\lambda$ to $E_{-\lambda}$ where $E_\lambda$ is the
eigenspace correspondent to eigenvalue $\lambda$. See \cite{Mu}. We assume that $\ker\;A\neq 0$. Then the automorphism $I$ defines a 
complex structure on $\ker\;A$ which induces a symplectic structure on it.
In particular, the dimension of $\ker\;A$ is even. The operator $D$ as
described will be said cylindrical compatible.

\medskip

{\bf Definition 1}{\it

An exponentially small perturbation of a cylindrical compatible operator $D$
is a first order formally self-adjoint elliptic operator $D^\prime$ satisfying
the following conditions:
\begin{itemize}
\item [{a)}] $D^\prime$ is a zero order perturbation of $D$,
\item [{b)}] on the cylindrical end $(-1,\infty)\times Y$, $D^\prime=D+P(t)$
where $P(t): \Gamma(E_0)\rightarrow\Gamma(E_0)$ is a smooth family of zero
order self-adjoint operators satisfying the following exponential decay
conditions: there exist a small $\delta >0$, some $T_0 >0$ and a constant $C$
such that when $t > T_0$,
$$
\|P(t)\psi\|_{L^2(Y)}\leq Ce^{-\delta(t-T_0)}\|\psi\|_{L^2(Y)}\hspace{2mm} \mbox{and}
\hspace{2mm}
\|\frac{\partial P}{\partial t}\psi\|_{L^2(Y)}\leq Ce^{-\delta(t-T_0)}\|\psi\|_{L^2(Y)}
$$
for $\psi\in L^2(E_0)$.
\end{itemize}
}
\medskip

Let $D^\prime$ be an exponentially small perturbation of a cylindrical
compatible operator. The space of ``bounded'' harmonic sections of $D^\prime$
is denoted by $H_B(D^\prime)$, i.e.
$$
H_B(D^\prime)=\{\psi\in\Gamma(E) | D^\prime\psi=0, \|\psi\|_{C^0(X)} < \infty\}.
$$
The space of $L^2$ harmonic sections of $D^\prime$ is denoted by
$H_{L^2}(D^\prime)$, i.e.
$$
H_{L^2}(D^\prime)=\{\psi\in L^2(E) | D^\prime\psi=0\}.
$$
Let $\beta$ be a fixed cut-off function which is equal to one at $\infty$,
and $\pi: (-1,\infty)\times Y\rightarrow Y$ be the natural projection.

{\bf Lemma 2}{\it

There exists a small $\delta_1 >0$ such that for any $\psi\in H_B(D^\prime)$,
there exists a unique limiting value $r(\psi)\in\ker\;A$ such that
$$
\|\psi-\beta\pi^\ast r(\psi)\|_{L^2_{\delta_1}(E)} < \infty.
$$
In particular, $\psi\in H_{L^2}(D^\prime)$ if and only if $r(\psi)=0$.
Moreover, 
$$
dim H_B(D^\prime)-dim H_{L^2}(D^\prime)=\frac{1}{2}dim\ker\;A.
$$
}
\medskip

Now consider a pair of triples $(X_i,E_i,D_i^\prime)$ for $i=1,2$. Suppose
that there is an orientation reversing isometry $h: Y_1\rightarrow Y_2$ which
is covered by corresponding bundle maps which identify $A_1$ with $A_2$ in a
suitable way so that for any $L > 0$, we can form a triple 
$(X_L,E_L,D_L^\prime)$ where $X_L=X_1\setminus [L+1,\infty)\times Y_1
\bigcup_h X_2\setminus [L+1,\infty)\times Y_2$ with $h: (L,L+1)\times Y_1
\rightarrow (L+1,L)\times Y_2$ given by $h(L+t,y)=(L+1-t,h(y))$, $E_L=E_1
\bigcup_h E_2$, $D_L=D_1\bigcup_h D_2$ and $P_L=\beta_L P_1+(1-\beta_L)h^\ast P_2$ for 
some cut-off function $\beta_L$ supported in $(L,L+1)\times Y_1$ with $|\nabla\beta|\leq 2$, and 
$D_L^\prime=D_L+P_L$. Set
$$
\lambda_L=\inf_{\psi\neq 0}\frac{\int_{X_L}|D_L^\prime\psi|^2}
{\int_{X_L}|\psi|^2}.
$$
The purpose of this appendix is to investigate the behavior of $\lambda_L$ as
$L\rightarrow \infty$.


{\bf Definition 3}{\it

Suppose $D^\prime$, $D^\prime_1$ and $D^\prime_2$ are exponentially small
perturbations of cylindrical compatible operators.
\begin{itemize}
\item [{a)}] $D^\prime$ is said to be regular if $H_{L^2}(D^\prime)=0$.
\item [{b)}] $(D^\prime_1,D^\prime_2)$ is said to be a transversal pair if
$$
r(H_B(D^\prime_1))\bigcap h^\ast(r(H_B(D^\prime_2)))=\{0\}.
$$ 
\end{itemize}
}
\medskip

Here is the main result.

{\bf Theorem 4}{\it

\begin{itemize}
\item [{1)}] $\lambda_L=O(\frac{1}{L^2})$ as $L\rightarrow \infty$,
\item [{2)}] if $(D^\prime_1,D^\prime_2)$ is a regular transversal pair, then
for any function $\gamma(L)=o(\frac{1}{L^2})$ as $L\rightarrow\infty$, there
exists $L_0 > 0$ such that when $L > L_0$, we have
$$
\lambda_L > \gamma(L).
$$
In particular, $D^\prime_L$ is invertible for large $L$.
\end{itemize}
}

\medskip

We first introduce some notation. Let $\lambda_i$, $i\in{\bf Z}$ denote the 
eigenvalues of the operator $A$, and $u_i$ denote the corresponding eigensections.
Set $\mu=\inf_{\lambda_i\neq 0}|\lambda_i|$. For simplicity, we omit the 
subscript $L$ if no confusion is caused.

{\bf Lemma 5}{\it

There exist $L_0 > 0$ and $M > 1$  with the following significance. Assume that $\psi$ and $c$ satisfy $D^\prime\psi=c\psi$ with $\psi\neq 0$ and
$|c| < \delta(\mu)$ for some small $\delta(\mu)$, then 
$\psi$ can be rescaled so that $\|\psi\|_{C^0(X_L)} < M$ and 
one of the following conditions holds:
\begin{itemize}
\item either $\int_{X_1(L_0)}|\psi|^2$ or $\int_{X_2(L_0)}|\psi|^2$  is equal
to one,
\item either $\|\psi\|_{L^2(Y_1)}(L_0)$ or $\|\psi\|_{L^2(Y_2)}(L_0)$ 
 is greater than or equal to one.
\end{itemize}
Here $X_i(L_0)=X_i\setminus (L_0,\infty)\times Y_i$, $i=1,2$. 
}
\medskip

{\bf Proof:} 
Let $\Pi_1$, $\Pi_2$ be the $L^2$-orthogonal projection onto $\ker\;A$ and
$(\ker\;A)^\perp$. On the cylindrical neck of $X_L$, write $\psi=f_1+f_2$
where $f_1\in \ker\;A$ and $f_2\in (\ker\;A)^\perp$. Set
$\xi(t)=\int_{Y}|f_2|^2$.

Direct computation shows that
\begin{eqnarray*}
\frac{\partial f_1}{\partial t}
& = & I\Pi_1 P\psi - cI(f_1)\\
\frac{\partial f_2}{\partial t}
& = & -Af_2 + I\Pi_2 P\psi - cI(f_2)\\
\frac{\partial^2 f_2}{\partial t^2}
& = & (A^2-c^2)f_2 + IA\Pi_2 P\psi + I\Pi_2\frac{\partial P}{\partial t}\psi
+ I\Pi_2 P\frac{\partial\psi}{\partial t} + c\Pi_2 P\psi.
\end{eqnarray*}
For any $\epsilon > 0$, there exists $L_0 > 0$ such that on the neck
$[L_0, 2L+1-L_0]\times Y_1$ we have
\begin{eqnarray*}
\frac{\partial^2\xi}{\partial t^2}
& \geq & 2\int_{Y}(\frac{\partial^2 f_2}{\partial t^2}, f_2)\\
& \geq & K(\mu^2\|f_2\|^2_{L^2_1(Y)} - \epsilon\|f_2\|_{L^2_1(Y)} 
(\|f_1\|_{L^2(Y)} + \|f_2\|_{L^2(Y)}))
\end{eqnarray*}
for some constant $K$. Here $|c| < \delta(\mu)$ for some small $\delta(\mu)$.
If $\xi(t)$ reaches its maximum in an interior point
$t_0\in (L_0, 2L+1-L_0)$, then on the neck, we have
$$
\max\|f_1\|_{L^2(Y)}\geq\|f_1\|_{L^2(Y)}(t_0)\geq\frac{\mu^2-\epsilon}{\epsilon}
\max\|f_2\|_{L^2(Y)}.
$$
Otherwise, $\xi(t)=\|f_2\|_{L^2(Y)}^2$ reaches its maximum at the end points.

On the other hand, we have on the neck that
$$
\begin{array}{c}
\frac{\partial f_1}{\partial t}+cI(f_1)=I\Pi_1 P\psi\\
\frac{\partial (If_1)}{\partial t}-c(f_1)= -\Pi_1 P\psi.
\end{array} $$
Set $C=\left(\begin{array}{cc}
0 & c\\
-c & 0
\end{array} \right)$, then
$$
\left(\begin{array}{c}
f_1\\If_1
\end{array} \right)(t)=e^{-Ct}\int_{L_0}^{t}e^{Cs}\left(\begin{array}{c}
I\Pi_1 P\psi\\-\Pi_1 P\psi
\end{array} \right)ds + e^{-C(t-L_0)}\left(\begin{array}{c}
f_1(L_0)\\If_1(L_0)
\end{array} \right).
$$
This implies that on the interval $[L_0, 2L+1-L_0]$
$$
\|f_1\|_{L^2(Y)}(t)\leq c_{1}e^{-\delta (L_0-T_0)}
(\max\|f_1\|_{L^2(Y)} + \max\|f_2\|_{L^2(Y)}) + \|f_1(L_0)\|_{L^2(Y)}.
$$
If $\|f_2\|_{L^2(Y)}$ reaches its maximum in the interior, then
$$
\max\|f_1\|_{L^2(Y)}\leq 2\|f_1(L_0)\|_{L^2(Y)}
$$
for large enough $L_0$. If $\|f_2\|_{L^2(Y)}$ reaches its maximum at the end
points, assuming that it is the left end point without loss of generality, we have
$$
\max(\|f_1\|_{L^2(Y)} + \|f_2\|_{L^2(Y)})\leq 2
(\|f_1(L_0)\|_{L^2(Y)} + \|f_2(L_0)\|_{L^2(Y)})
$$
for large enough $L_0$. Lemma 5 follows easily from these estimates.
\hfill $\Box$

\noindent{\bf The Proof of Theorem 4:}

1). Pick $\phi\in\ker\;A$ with $\|\phi\|_{L^2(Y)}=1$. Let $\rho_L$ be a cut-off function which equals to one on $[\frac{3L}{4}, \frac{3L}{4}+
\frac{L}{2} +1]$ and equals to zero outside $[\frac{L}{2},\frac{L}{2}+L+1]$
with $|\nabla\rho_L|=O(\frac{1}{L})$. Then
$$
\int_{X_L}|D^\prime_L(\rho_L\phi)|^2\leq
\int_{X_L}|\nabla\rho_L|^2|\phi|^2+\int_{X_L}|P_L(\rho_L\phi)|^2
=O(\frac{1}{L}),\mbox{and} \int_{X_L}|\rho_L\phi|^2\geq \frac{L}{10}.
$$
So $\lambda_{L}=O(\frac{1}{L^2})$ as $L\rightarrow\infty$.

2). Suppose that there exists a sequence of $L_n\rightarrow\infty$ such that
$\lambda_{L_n}\leq\gamma(L_n)$. Then there exist $\psi_n$, $c_n$ such that
$D^\prime_{L_n}\psi_n=c_n\psi_n$ with $c_n^2=\lambda_{L_n}$. By Lemma 5,
there exist $\psi_1\in H_B(D^\prime_1)$, $\psi_2\in H_B(D^\prime_2)$ such that
a subsequence of $\psi_n$ converges to $\psi_1$ over $X_1$ and $\psi_2$ over
$X_2$ in $C^k$ norm on any compact subset. Note that one of $\psi_1$ and 
$\psi_2$ is nonzero. Part 2 of Theorem ~4 follows if we show that
$r(\psi_1)=h^\ast r(\psi_2)$. But this follows from the fact that if we write
$\psi=f_1 + f_2$ as in Lemma 5, 
\begin{eqnarray*}
\|f_1(t)-f_1(2L+1-t)\|_{L^2(Y)}
& \leq &
C(e^{-\delta t} + |\cos(c(2L+1-2t))-1|\\
&      & + |\sin(c(2L+1-2t))|),
\end{eqnarray*}
for large enough $t$ and $L$. $C$ is some constant independent of $t$ and $L$.


\noindent{\bf The Proof of Lemma 2:}

Suppose $\psi\in\Gamma(E)$ and $D^\prime\psi=0$. On the cylindrical end
$(T_0,\infty)\times Y$, write $\psi=\sum_{i}f_i u_i$ where $u_i$ are the
eigensections of the operator $A$ corresponding to eigenvalues $\lambda_i$, and
$f_i$ are smooth functions in $t$. Then we have
$$
\frac{\partial f_i}{\partial t} + \lambda_i f_i =
(IP(t)\psi, u_i).
$$
Set $g_i=(IP(t)\psi, u_i)$, then $\sum_{i}g_i^2=\|P\psi\|_{L^2(Y)}^2$ and
$$
f_i(t)=\int_{T_0}^{t}e^{-\lambda_i(t-s)}g_i(s)ds 
+ f_i(T_0)e^{-\lambda_i (t-T_0)}.
$$
Now assume that $\psi\in L^2_{-\gamma}$ for any small enough $\gamma > 0$.
Assume that $\delta_1 < \min( \frac{\delta}{2}, \frac{\mu}{4})$ where
$\mu=\inf_{\lambda_i\neq 0}|\lambda_i|$. 
\begin{itemize}
\item For $\lambda_i=0$, we have for any $t^\prime > t$,
$$
e^{\delta_1 t}|f_i(t^\prime) - f_i(t)|\leq
C(\int_{t}^{t^\prime}\int_{Y}e^{-\frac{\delta}{10}s}|\psi|^2{Vol_Y}ds)^
{\frac{1}{2}},
$$
so $f_i(\infty)=\lim_{t\rightarrow\infty}f_i(t)$ exists and
$f_i-f_i(\infty)\in L^2_{\delta_1}$.
\item For $\lambda_i > 0$, we have for some constant $C(\mu)$ that
$$
e^{2\delta_1 t}(\sum_{i}f_i^2(t))\leq
C(\mu)\int_{T_0}^{\infty}e^{2\delta_1 s}(\sum_{i}g_i^2(s))ds
+ (\sum_{i}f_i^2(T_0))e^{2\delta_1 T_0}.
$$
\item For $\lambda_i < 0$. First of all, we have 
$$
f_i(t)=-e^{-\lambda_i t}\int_{t}^{\infty}e^{\lambda_i s}g_i(s)ds
$$
since $\psi\in L^2_{-\gamma}$ for any small enough $\gamma > 0$.
On the other hand, for some constant $C(\mu)$, we have
$$
e^{2\delta_1 t}(\sum_{i}f_i^2(t))\leq
C(\mu)\int_{t}^{\infty}e^{2\delta_1 s}(\sum_{i}g_i^2(s))ds.
$$
\end{itemize}
Take $r(\psi)=\sum_{i}f_i(\infty)u_i$ where $u_i\in \ker\:A$, then
$$
\|\psi-\beta\pi^\ast r(\psi)\|_{L^2_{\delta_1}(E)} < \infty
$$
where $\beta$ is a fixed cut-off function which is equal to one at $\infty$,
and $\pi: (-1,\infty)\times Y\rightarrow Y$ is the natural projection.
As for $\mbox{dim}H_B(D^\prime)-\mbox{dim} H_{L^2}(D^\prime)=\frac{1}{2}\mbox{dim}\ker\;A$, it follows from Theorem 7.4 in \cite{LM}.


\end{document}